\documentclass[11pt]{article}

\usepackage{arxiv}

\usepackage{amsmath,amsfonts,latexsym,fullpage,graphicx,vmargin, mathrsfs}
\usepackage{pstricks}
\usepackage{comment}
\usepackage[inline]{enumitem}
\usepackage{cleveref}

\usepackage[ruled]{algorithm2e}
\usepackage{algpseudocode}

\RequirePackage[OT1]{fontenc}

\usepackage{amsmath}
\usepackage{amsfonts}
\usepackage{caption}
\usepackage{subcaption}
\usepackage{multirow}
\usepackage{tikz}
\usetikzlibrary{graphs}
\usetikzlibrary{cd}
\usepackage{bm}
\usepackage{url}
\usepackage{placeins}

\RequirePackage{natbib}

\newtheorem{teo}{Theorem}
\newtheorem{prop}{Proposition}
\newtheorem{defi}{Definition}
\newtheorem{rmk}{Remark}
\newtheorem{lem}{Lemma}
\newtheorem{cor}{Corollary}

\DeclareMathOperator*{\im}{Im}
\DeclareMathOperator*{\Top}{Top}

\DeclareMathOperator*{\Sets}{Sets}

\DeclareMathOperator*{\Mapp}{Mapp}

\DeclareMathOperator*{\len}{len}
\DeclareMathOperator*{\lvl}{lvl}

\usepackage{scalerel}
\DeclareMathOperator*{\Bigcdot}{\scalerel*{\cdot}{\bigodot}}

\newcommand{\X}{X_{\Bigcdot}}
\newcommand{\Y}{Y_{\Bigcdot}}

\newcommand{\T}{\pi_0(\X)}
\newcommand{\G}{\pi_0(\Y)}

\newcommand{\R}{\mathbb{R}}

\newcommand{\virgolette}[1]{``#1''}

\newcommand{\Gcal}{\mathcal{G}}

\firstpageno{1}

\begin{document}

\title{A Persistence-Driven Edit Distance for Trees with Abstract Weights}

\author{Matteo Pegoraro\thanks{Department of Mathematical Sciences, Aalborg University}}

%\author{\name Matteo Pegoraro \email matteo.pegoraro@polimi.it \\
%       \addr MOX --- Department of Mathematics \\
%       Politecnico di Milano \\
%       Milano, Italy
%       \AND
%       \name Mario Beraha \email mario.beraha@polimi.it \\
%       \addr Department of Computer Science \\
%       Universit\`a degli Studi di Bologna \\
%       Bologna, Italy}
%
\maketitle

\begin{abstract}
In this work we define a novel edit distance for trees considered with some abstract weights on the edges. 
The metric is driven by the idea of considering trees as topological summaries in the context of persistence and topological data analysis. Several examples related to persistent sets are presented.
The metric can be computed with a dynamical binary linear programming approach. This framework is applied and further studied in other works focused on merge trees, where the problems of stability and merge trees estimation are also assessed.
\end{abstract}



\section{Introduction}

Graphs are widely used across scientific disciplines. Whenever a system is discretized, or when relations among points, statistical units, or abstract objects need to be formalized, graphs naturally emerge as a fundamental modeling tool.

This ubiquity has led to the development of a broad array of mathematical techniques for analyzing and comparing graphs, ranging from differential calculus on discrete structures and graph metrics, to specialized algebraic frameworks. Over time, this has given rise to a rich and diverse landscape of graph-theoretic methods.

Within this broader context, trees occupy a particularly important role. Their structural simplicity, combined with considerable expressive power, offers an appealing balance between interpretability and computational tractability. As a result, trees have become central in many areas of science. Two especially prominent examples are phylogenetic trees \cite{billera2001geometry, phylo, garba2021information, lueg2024foundations} and clustering dendrograms \cite{dendro_1, dendro_2}, which are closely related structures. Both describe how a set of elements evolves or merges under some notion of similarity or agglomeration, and both illustrate the versatility of trees as a means of encoding hierarchical information.

In this paper, we approach trees from the perspective of topological data analysis (TDA) \citep{PH_survey}, interpreting them as carriers of topological information. Our goal is to develop an edit distance tailored to analyzing populations of tree-shaped topological summaries. 

TDA is a relatively recent framework for data analysis that leverages tools from algebraic topology \citep{hatcher} to extract shape-related features from data at the level of individual statistical units.

A typical TDA pipeline begins with a topological space and constructs a nested sequence of subspaces—known as a \emph{filtration}—from which topological information is extracted using functors such as simplicial homology \citep{hatcher}. This standard approach has been generalized in several directions, including filtrations indexed by more general partially ordered sets (e.g., $\mathbb{R}^n$ for $n \geq 2$) \citep{botnan2022introduction}, and the use of other mathematical objects such as Laplacians \citep{memoli2022persistent}, path-connected components \citep{curry2021decorated, curry2021trees, pegoraro2024finitely, cavinato2022imaging}, zigzag homology \citep{carlsson2010zigzag}, cosheaves \citep{de2016categorified} among others.

\subsection*{Previous Works on Edit Distances for Trees in TDA}

The topological perspective we assume in this work stems from  \emph{persistent sets}, as defined in \cite{carlsson2013classifying}, which are functors  
 of the form $F:(P,\leq)\rightarrow \Top \xrightarrow{\pi_0} \Sets$ (with $(P,\leq)$ being a poset, $\Top$ being the category of topological spaces and $\pi_0$ being the functor of path-connected components). The most typical example being $(P,\leq)$ equal to $\{f^{-1}((-\infty,t])\}_{t\in \R}$ for some $f:X\rightarrow \R$, which amounts to using \emph{merge trees} of real valued functions \citep{merge_interl, curry2021trees}, while, if $(P,\leq)$ is the poset of open intervals in $\R$, up to some gluing conditions, we recover \emph{Reeb cosheaves} \citep{de2016categorified}.

Relations between posets and graphs have been extensively studied: for instance it is well-known that any poset can be represented with a directed acyclic graph (DAG) \citep{transitive_reduction}. Thus, it is no surprise that a common way to represent instances of $F:(P,\leq)\rightarrow \Sets$ is via some kind of graphs, like merge trees or Reeb graphs.
Accordingly, the use of graphs (possibly labeled \cite{yan2019structural}) as summaries of persistent sets has produced a series of original definitions and results 
\citep{bauer2014measuring,
de2016categorified,
 di2016edit, bauer2020reeb, merge_interl, merge_intrins, merge_wass, merge_farlocca, merge_farlocca_2, sridharamurthy2021comparative,
cardona2021universal,
merge_frechet,
 curry2021decorated} which are all driven by these novel scenarios.

A significant number of these works focus on defining suitable metric structures for merge trees and Reeb graphs through variants of edit distances. In particular, 
\citep{ di2016edit, bauer2020reeb} are focused on sup-norm kinds of distances, where the difference between (Reeb) graphs is determined by the cost of the biggest modification needed to turn one graph into another. While \cite{merge_wass, merge_farlocca, merge_farlocca_2, sridharamurthy2021comparative} consider edit distances for merge trees which are more closely related to classical tree edit distances \cite{Tai}, where the objective is to minimize the total cost of a sequence of modifications converting one merge tree into another.

In this work, we propose an edit distance for \emph{weighted} trees that aligns more closely with the second family of approaches mentioned above, but it differs from previous methods in two key ways:
\begin{itemize}
    \item Stability Properties:  in \cite{pegoraro2024finitely} and \cite{pegoraro2024functional} it is shown that, when adapted to merge trees, the metric we define here produces an edit distance  $d_E$ satisfying the following stability relation :
    \[
    d_I(T,T')\leq d_E(T,T')\leq 2(\text{size}(T)+\text{size}(T'))d_I(T,T')
    \]
    where $d_I$ is interleaving distance between merge trees \citep{merge_interl} and $\text{size}(T)$ is the number of edges in $T$. Which is analogous to the relation  between the bottleneck distance ($d_B$) and the $1$-Wasserstein metric ($W_1$) between persistence diagrams (PDs):
    \[
    d_B(D,D')\leq W_1(D,D')\leq (\#D + \#D')d_B(D,D').
    \]
    See \cite{edelsbrunner2022computational}. Among the edit distances presented in \cite{merge_wass, merge_farlocca, merge_farlocca_2, sridharamurthy2021comparative}, only the one in \cite{merge_farlocca_2} was shown—in \cite{pegoraro2024finitely}—to posses some stability properties.       
    This stability result also clarifies the theoretical distinction between our approach and those in \cite{bauer2014measuring, di2016edit, bauer2020reeb}, all of which yield distances that are bounded from above by $d_I$. Since $d_B$ is the \emph{universal metric} between PDs and $d_I$ is the universal one between merge trees, by analogy, one can further characterize the difference between $d_E$ and \cite{di2016edit, bauer2020reeb}, via the easier and well understood difference between $d_B$ and $W_1$.
    \item Abstract Weights: we investigate generalizations in which trees are enriched not only with edge lengths, but with weights taking values in more general metric spaces. This allows topological summaries to capture and represent a wider variety of information.
\end{itemize}

\subsection*{Main Contributions}

As anticipated, the main contribution of this work is the introduction of a novel edit distance designed to compare trees interpreted as topological summaries. Increasingly, such trees are enriched with additional information, typically associated with the vertices or edges of the graph, as seen in various recent works \citep{ushizima2012augmented, curry2021decorated, pegoraro2024finitelyfunc, curry2023convergence}. To remain as general as possible, we allow abstract weights to be assigned to the edges of the trees.

Given the wide range of settings that emerge from this enriched perspective (as discussed in \Cref{sec:examples} and reflected in the diverse works built on the framework we establish) this paper focuses on the following core contributions:
\begin{enumerate}
    \item Defining a general edit distance for trees with abstract weights;
    \item Formalizing the abstract properties that edge weights must satisfy for the metric to be well-defined;
	 \item	Developing a computational framework for efficiently computing the metric;
	 \item	Beginning to explore the variety of edge weights that can be used, discussing their potential applications and how can they support diverse data analysis scenarios. We leave a deeper investigation into the interpretability of the metric in specific scenarios and the study of its stability properties to other works.
\end{enumerate}

\begin{figure}
    \centering
	\begin{subfigure}[c]{0.30\textwidth}
    	\centering
    	\includegraphics[width = \textwidth]{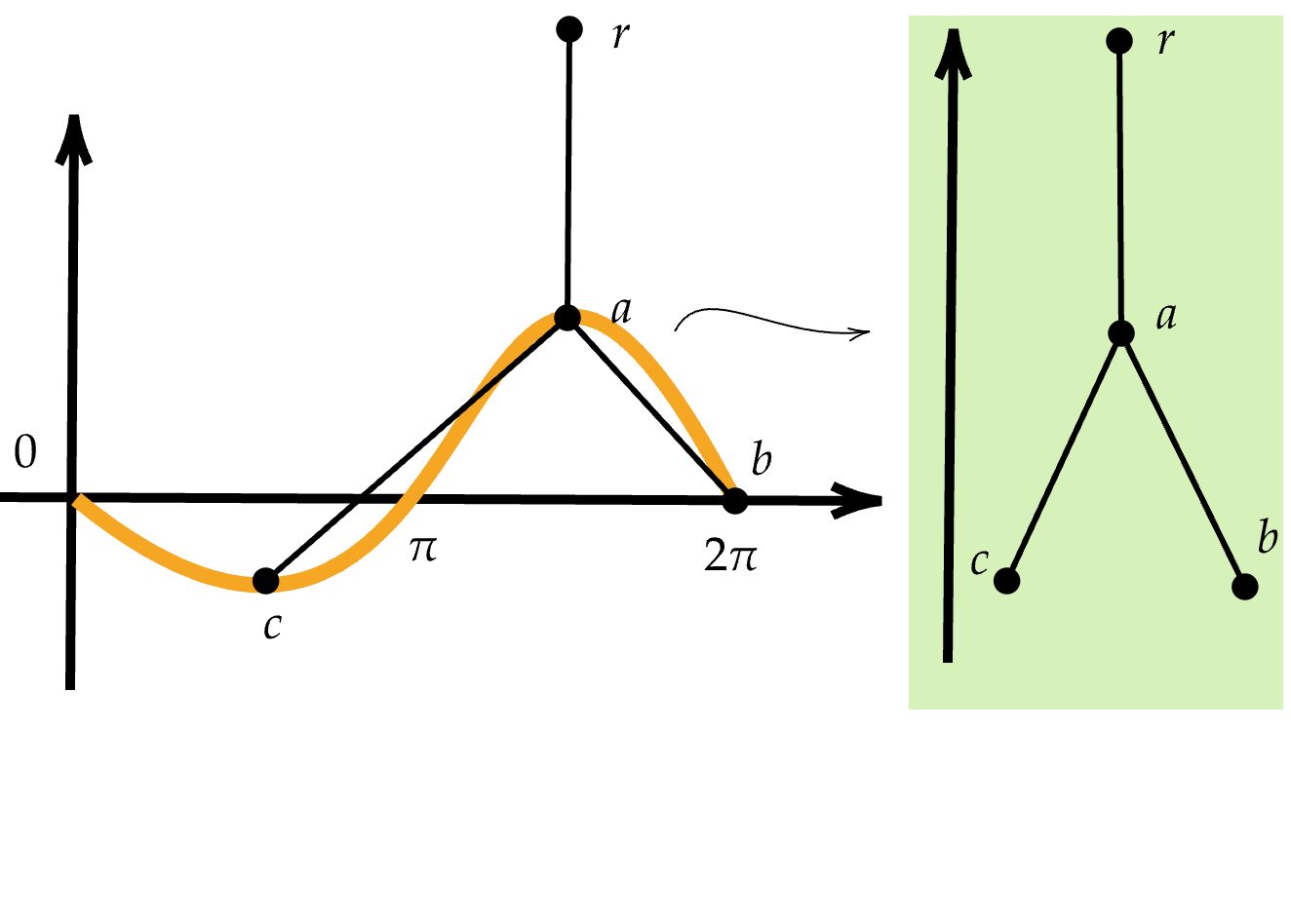}
    	\caption{A function (orange) with its associated merge tree.}
    	\label{fig:func_ex_0}
    \end{subfigure}
   	\begin{subfigure}[c]{0.30\textwidth}
    	\centering
    	\includegraphics[width = \textwidth]{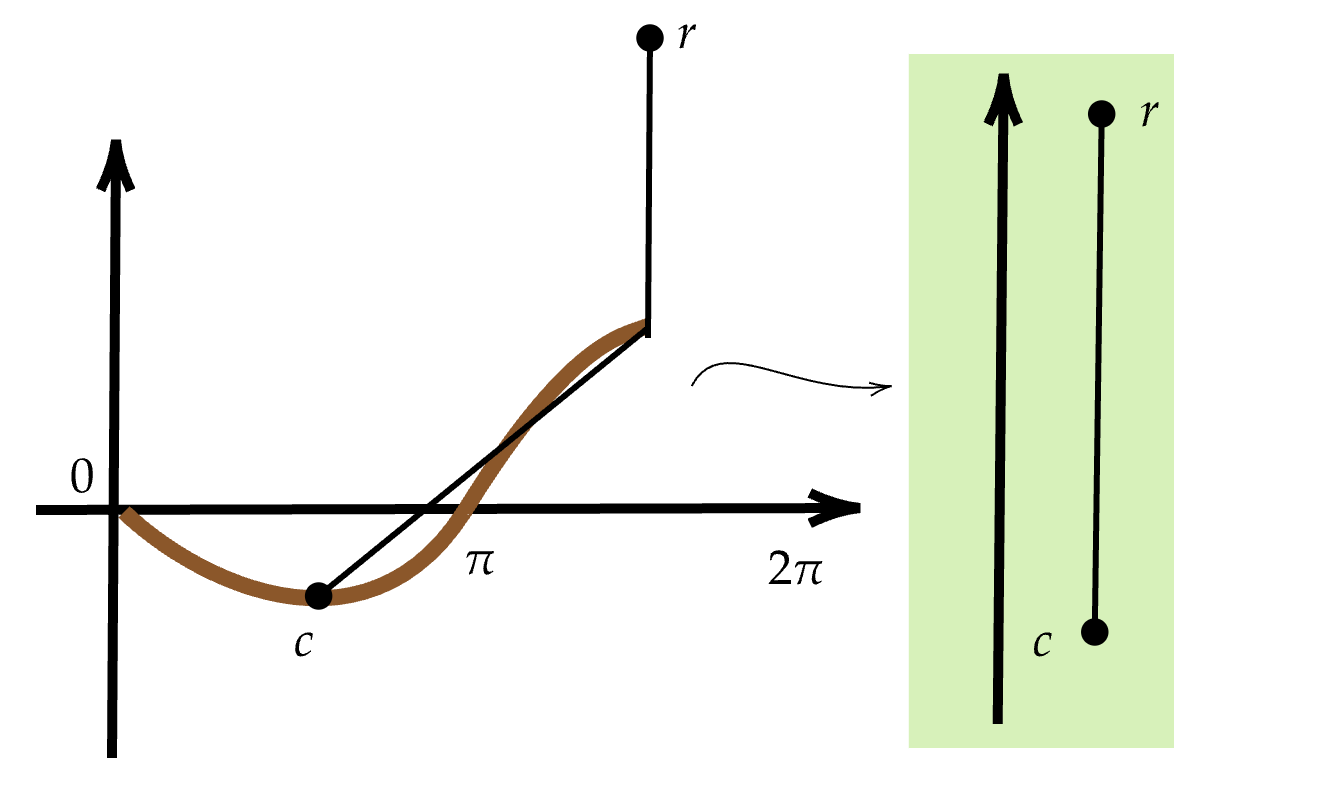}
    	\caption{Another function (brown) with its associated merge tree in black.}
    \label{fig:func_ex_1}
    \end{subfigure}
   	\begin{subfigure}[c]{0.30\textwidth}
    	\centering
    	\includegraphics[width = \textwidth]{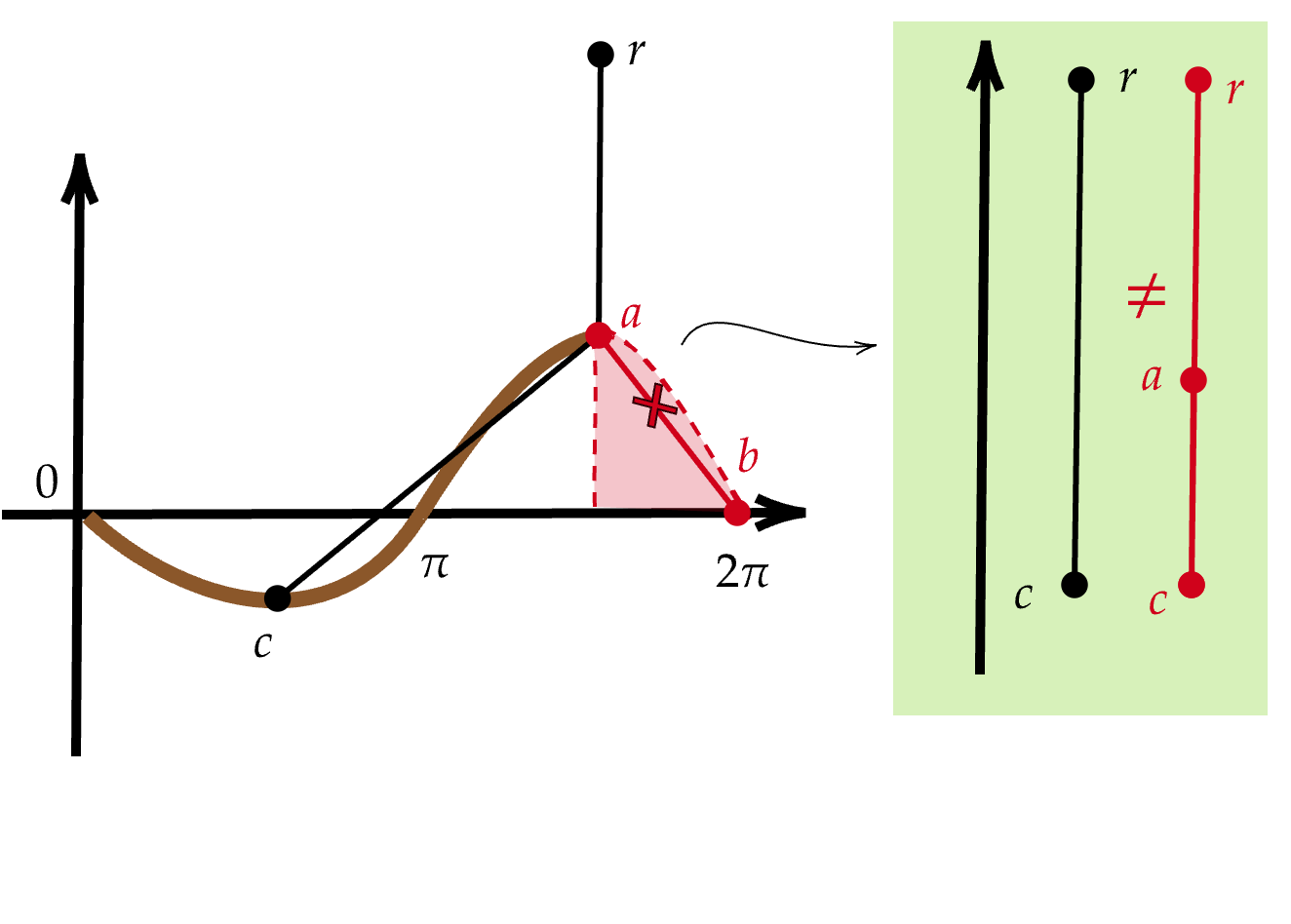}
    	\caption{Comparison bewtween the merge tree in \Cref{fig:func_ex_1}, and the tree obtained from the one in \Cref{fig:func_ex_0}, removing the edge $(a,b)$.}
    \label{fig:func_ex_2}
    \end{subfigure}
%    \hfill
\caption{Plots related to the example presented in the Introduction.}
\label{fig:func_intro}
\end{figure}

The edit distance we define stems from classical edit distances \citep{levenshtein1966binary, Tai, survey_ted, gao2010survey, lerouge2016exact}, but add some deep modifications, which can be intuitively motivated by the following example, represented in \Cref{fig:func_intro}. 

Consider the following function: $f(x)=x\cdot \sin( x)$ for $x\in [0,2\pi]$ - \Cref{fig:func_ex_0}. Then take the sublevel set filtration $F(t)=f^{-1}((-\infty,t])$ and the persistent set $\pi_0\circ F$. The information of this persistent set can be represented with a tree with four vertices and three edges, as in \Cref{fig:func_ex_0}, along with a real-valued height function defined on its vertices - with the root being at height $+\infty$. We can associate to any edge its length, obtained as 
the difference in height between its endpoints. Suppose that we were to remove from the (merge) tree the edge $(a,b)$, we would obtain the red (merge) tree in \Cref{fig:func_ex_2}. Vertex $a$ is now a degree $2$ vertex which represents a merging between path connected components (i.e. a local maximum in the function), which is no more happening. 
On a function level, this can be visualized as removing the local minimum in $2\pi$, for instance restricting the function on $[0,3/2\pi]$ - \Cref{fig:func_ex_1}. Usual edit distances would regard the red and black weighted  trees in the green box in \Cref{fig:func_ex_2} as very different trees, having a different number of edges and very different weights associated to the edges. As a consequence also the black trees in the green boxes in \Cref{fig:func_ex_0} and in \Cref{fig:func_ex_1} would be very far apart. In our case however, we would like our metric to be much less sensitive and, for instance, not to distinguish between the black  and red trees in  trees \Cref{fig:func_ex_2}. 
For this reason, differently from other edit distances, we add the possibility to remove the vertex $a$ and join the two adjacent edges.

\subsection*{Interconnected Works}

The framework introduced in this work lays the foundation for a series of developments explored in separate studies: the edit distance we define is applied across different contexts, each accompanied by dedicated theoretical results that both justify the specific application and reinforce the foundational framework developed in this manuscript. 

We summarize such relations as follows:

\begin{itemize}
\item All the works we mention in the following points rely heavily on \Cref{sec:mappings} and \Cref{teo:main_thm} for developing theoretic results, and on Algorithm \ref{alg:bottomup} for applications and case studies.
\item Merge trees: this general framework is used by \cite{pegoraro2024finitely} to define a metric for merge trees, which can be computed with Algorithm \ref{alg:bottomup}. With the help of some results contained in \cite{pegoraro2024functional}, which is more focused on the statistical aspects of such metric, \cite{pegoraro2024finitely} shows that the approach we pursue here induces a metric with stability properties which are analogous to the ones of the $1$-Wasserstein distance between persistence diagrams. 
\cite{pegoraro2024functional} and \cite{cavinato2022imaging} 
present real data applications to show the practical effectiveness of the metric.
\item Functions defined on merge trees: \cite{pegoraro2024finitelyfunc} exploits the present work,  exploring the case where weights are functions defined on the edges of the trees, obtaining another stability result. 
\end{itemize}

%A simplified but similar definition has already been considered in \cite{bio}, but it 
%is just cited in few lines as a possibility without a real motivation, which lacks 
%any kind of investigation (even whether or not it defines a proper metric).

\subsection*{Outline}

\Cref{sec:edit_space} introduces the idea of \emph{editable spaces}, which are the spaces we want the weights of our graphs to live in. Their properties are exploited in \Cref{sec:trees} to define \emph{dendrograms}.
\Cref{sec:edits} develops the novel edit distance and proves that it is well defined for trees considered up to degree $2$ vertices. \Cref{sec:examples} presents a series of general examples regarding persistent sets.
\Cref{sec:algo} obtains an algorithm to compute the distance via linear binary optimization.
\Cref{sec:conclusions} concludes the paper by highlighting the pros and cons of the approach along with possible further developments, with a brief comparison with graph edit distances.

\section{Editable Spaces}
\label{sec:edit_space}

As outlined in the introduction, the goal of this manuscript is to introduce a novel metric structure for working with weighted trees, where the term \emph{weighted}  is interpreted broadly: edge weights are defined via a map that assigns each edge a value in a metric space.
To define this metric meaningfully, we require the metric spaces to satisfy certain properties. These conditions are essential to guarantee well-behaved transformations when modifying trees.

\begin{defi}
A set $X$ is called editable if the following conditions are satisfied: 
\begin{itemize}
\item[(P1)] $(X,d_X)$ is a metric space
\item[(P2)] $(X,\odot,0_X)$ is a commutative monoid (that is $X$ has an associative operation $\odot$ with zero element $0_X$)
\item[(P3)] the map $d_X(0_X,\cdot):X\rightarrow \mathbb{R}$
is a map of monoids between $(X,\odot)$ and $(\mathbb{R},+)$:
$d_X(0_X,x\odot y)= d_X(0_X,x)+d_X(0_X,y)$.
\item[(P4)] $d_X$ is $\odot$ invariant, that is: $d_X(x,y)=d_X(z\odot x,z\odot y)=d_X(x\odot z,y\odot z)$
\end{itemize}
\end{defi}
Note that in property (P3), $d_X(x\odot y,0_X)=d_X(x,0_X)+d_X(y,0_X)$, implies that $x\odot y\neq 0_X$. Moreover (P3)-(P4)
imply that the points $0_X$, $x$, $y$ and $x\odot y $ form a rectangle which can be isometrically embedded in an Euclidean plane with the Manhattan geometry (that is, with the norm $\mid \mid \cdot\mid \mid _1$): 
$d_X(x,x\odot y)=d_X(0_X,y)$, $d_X(y,x\odot y)=d_X(0_X,x)$
and $d_X(x\odot y,0_X)= d_X(0_X,x)+d(0_X,y)$.

Exploiting the monoidal structure of an editable space, we can turn it into a partially ordered set.

\begin{prop}\label{prop:editable_poset}
An editable space $X$ admits a poset structure given by:
\[
x\leq y \text{ if and only if } 
\exists r\in X \text{ such that }
x\odot r=y.
\]
\begin{proof}
\begin{itemize}
\item Reflexivity: $x\odot 0=x$ and thus $x\leq x$.
\item Antisymmetry: suppose $x\leq y$ and $x\geq y$. Then $x\odot r=y$ and $y\odot q=x$ for some $r,q\in X$.
But then:
\begin{align*}
 d_X(y,0)=&d_X(x\odot r,0)=d_X(x,0) + d_X(r,0) \\
 =&d_X(y\odot q,0)+ d_X(r,0)=d_X(y,0)+ d_X(q,0)+d_X(r,0),
\end{align*}
 
which implies $r=q=0$ and so $x=y$.
\item Transitivity follows from associativity of $\odot$.
\end{itemize}
\end{proof}
\end{prop}

If the editable space $X$ is contained in a group $V$, the partial order can be further characterized. 

\begin{cor}\label{cor:editable_poset}
Consider $X$ editable space. If $X\subset V$ with $V$ being a group and the inclusion being a map of monoids, then for every $x\leq y$ in $X$, there exist one and only one $r\in X$ such that $x+r=y$. 
\end{cor}

\subsection{Examples of Editable Spaces}
\label{sec:edit_examples}

We give some examples of editable spaces.

\subsubsection{Positive Real Numbers}\label{sec:positive_editable}
The set $(\mathbb{R}_{\geq 0}, +, \mid \cdot \mid )$ is an editable space, as well as its subsets which are monoids, like $\mathbb{N}$.
This editable space can be used to recover the classical notion of (positively) weighted graphs.

\subsubsection{$1$-Wasserstein Space of Persistence Diagrams}

Persistence diagrams (PDs) are among the most used topological summaries in TDA.
They are used to represent sequences of homology groups, often called \emph{persistence modules}, arising by applying homology functors to filtrations of topological spaces. PDs are frequently used due to their interpretable nature, being finite set of points in the plane  
whose $x$ and $y$ coordinates - $y\geq x$ - represent the birth and death of homological features along the filtration.

There are many possible ways to introduce persistence diagrams and Wasserstein metrics, like \cite{edelsbrunner2022computational},
\cite{skraba2020wasserstein}, or 
\cite{divol2021understanding}. 
In particular, \cite{edelsbrunner2022computational}
is suitable for interested readers not familiar with the topic.
For notational convenience here we follow the approach in \cite{bubenik2022universality}, which we use in the proof of an upcoming proposition.

\begin{defi}[adapted from \cite{bubenik2022universality}]\label{def:PD}
Define $D(\R^2_\Delta)$ as the set of finite formal sums of elements in:
\[
\R^2_\Delta =\{(b,d)\in R^2\mid 0\leq b \leq d\}. 
\]
Similarly $D(\Delta)$ is the set of formal sums of elements in
 $\Delta=\{(b,b)\in \R^2\mid 0\leq b\}$. In both sets of formal sums we have an operation given by the sum of formal sums, which is commutative and has the empty formal sum as zero element.
Given $D_1,D_2\in D(\R^2_\Delta)$ we say that $D_1\sim D_2$ if there exist $D_3,D_4\in D(\Delta)$ such that $D_1+D_3= D_2+D_4$. 
We collect 
the equivalence classes of such relations in $D(\R^2_\Delta,\Delta)$ and call them persistence diagrams, and indicate such an equivalence class with the notation $[D]$. The sum operation and the zero element are coherently defined on the quotient.
\end{defi}

As a consequence of \Cref{def:PD} we can write every diagram $[D]$ as a sum of points in $\R^2_\Delta$ each with coefficient $1$. This \virgolette{decomposition} is clearly non-unique as it varies between elements in the same equivalence class. With an abuse of notation we may say that $x\in [D]$, meaning that $x=(b,d)$ is an off-diagonal point of $\R^2_\Delta$ (i.e. $d>b$) appearing in such decomposition of an element (in fact, in any element) in the equivalence class of $D$.

\begin{defi}[adapted from \cite{bubenik2022universality}]
\label{defi:PD}
The $1$-Wasserstein metric between persistence diagrams is:
\[
W_1([D],[D'])=\inf_{\gamma:\widetilde{D}\rightarrow \widetilde{D}'} \parallel x-\gamma(x)\parallel_1,
\]
where $\gamma$ is any bijection between some $\widetilde{D}$ such that $\widetilde{D} \sim D$, and some $\widetilde{D}'$ such that $\widetilde{D}'\sim D'$. 
The infimum in the definition of $W_1([D],[D'])$ is always attained by some $\gamma$, and thus, it is always a minimum.
\end{defi}

\begin{prop}\label{prop:wass_editable}
The space $(D(\R^2_\Delta,\Delta),W_1,+,0)$ is editable.
\end{prop}

\begin{rmk}
We assumed \Cref{defi:PD} as it makes \Cref{prop:wass_editable} easier to prove, compared to other formulations. However, for the sake of simplicity, we drop the notation $[D]$ in the remaining of the manuscript referring to a diagram simply as $D$, chosing the representative with no points of the form $(b,b)$.
\end{rmk}

\subsubsection{Functions in Editable Spaces}\label{sec:fun_weights}

Consider an editable space $(X,d_X,\odot_X,0_X)$ and a measure space $(M,\mu)$. Then, the space of 
functions:
\[
L_1(M,\mu,X):=\lbrace f:M\rightarrow X \mid \int_{M} d_X(0_X,f(p))d\mu(p)<\infty \rbrace/\sim,
\]
 with $\sim$ being the equivalence relation  of functions up to zero measure subsets and $\odot_W$ being defined pointwise - is editable. 
In fact, the function $p\mapsto d_X(0,f(p))$ is always non negative, so 
if properties (P3) and (P4) hold for a fixed $p\in M$, then 
they hold also for integrals. Call $W:=L_1(M,X)$, we verify 
(P3) as follows: 
\begin{align*}
d_W(f\odot_W g,0)=& \int_M d_X(f(p)\odot_X g(p),0)d\mu(p) \\
=&\int_M d_X(f(p),0)+ d_X(g(p),0)d\mu(p) = d_W(f,0)+d_W(g,0).
\end{align*}
And similarly for (P4).

\subsubsection{Finite Products of Editable Spaces}
\label{sec:product_spaces}
Consider  two editable spaces $X$ and $X'$, that is $(X,\odot,0_X)$ and $(X',\diamond,0_{X'})$ satisfying properties (P1)-(P4). Then 
$(X\times X',\star, (0_X,0_{X'}))$ 
is an editable space, with $\star$ being the component-wise operations $\odot$ and $\diamond$, and the metric $d$ on $X\times X'$ being the (possibly weighted) sum of the component-wise metrics of $X$ and $X'$.
For instance, $\R^n_{\geq 0}$ with the $1$-norm is an editable space.

\section{Dendrograms}
\label{sec:trees}

Now we introduce the objects we will work with: trees,  with a function associating to each edge a weight in a fixed editable space.

\begin{defi}\label{def:tree_struct}
A tree structure $T$ is given by a connected rooted acyclic graph $(V_T,E_T)$.  We indicate the root of the tree with  $r_T$. The degree of a vertex $v\in V_T$ is the number of edges which have that vertex as one of the endpoints.
Any vertex with an edge connecting it to the root is its child and the root is its parent: this is the first step of a recursion which defines the parent  and children relation  for all vertices in \(V_T.\)
%In this way we recursively define parent  and children (possibly none) for any vertex on the tree. 
The vertices with no children are called leaves and are collected in the set $L_T$. The relation $child < parent  $ generates a partial order on $V_T$ which induces an orientation on the edges. The edges in $E_T$ are identified in the form of ordered pairs $(a,b)$ meaning $a<b$ and $a\rightarrow b$.
A subtree of a vertex $v$, called $sub_T(v)$, is the tree structure whose set of vertices is $\{x \in V_T\mid x\leq v\}$. 
\end{defi}

A key fact is that given a tree structure $T$, identifying an edge $(v,v')$ with its lower vertex $v$, gives a bijection between $V_T-\{r_T\}$ and $E_T$, that is $E_T\cong V_T-\{r_T\}$ as sets. 
Given this bijection, we often use $E_T$ to indicate the vertices $v\in V_T-\{r_T\}$,  to simplify the notation. 

To avoid overlapping some terminology with already existing notions of trees and for ease of notation, we refer to a weighted tree, with the weight function taking values in an editable space, by calling it dendrogram, as formalized by the following definition. 

\begin{defi}
A tree structure $T$ with a weight function $\varphi_T:E_T\rightarrow X-\{0_X\}$, with $X$ editable space, is called dendrogram. The space of $X$-valued dendrograms is indicated as $(\mathcal{T},X)$.
\end{defi}

\subsection{Finite Posets, DAGs and Transitive Reductions}
\label{sec:trans_red}

In this section, we justify our focus on tree-shaped objects by outlining, at a high level, the various contexts in which such graphs arise in TDA. This also yields a conceptual blueprint for potential extensions of the present work.

We start from the definition of a directed graph.

\begin{defi}
A directed graph is made by a pair $G=(V,E)$ with $V$ being the set of its vertices and $E\subseteq V\times V$ the set of its directed edges. 
\end{defi}

Consider a finite partially ordered set $(P,\leq)$. 
It is well known
that one can associate to $P$ a unique (up to directed graph isomorphism) directed acyclic graph (DAG) $G_P$ so that $G_P$ is in bijection with the objects in $P$
and $(a,b)\in E_{G_P}$ if and only if $a\leq b$. Moreover, the graph $G_P$ can be simplified without losing ordering information. We know formalize these constructions.

\begin{defi}[\cite{transitive_reduction}]
A directed graph $G=(V,E)$ is called transitive if $(v,v')\in E$, whenever there is a directed path from $v$ to $v'$. For every directed graph $G=(V,E)$ we define its transitive closure as the graph $G^T:= (V,E_{G^T})$, with:
\[
E_{G^T} = \bigcap_{\{E'\subset V\times V \mid E\subset E', (V,E') \text{ transitive graph}\}} E'. 
\]
\end{defi}

\begin{teo}[adapted from \cite{transitive_reduction}]
For every $G=(V,E)$ DAG  there is a unique graph $G^t=(V,E_{G^t})$ such that:
\begin{itemize}
\item $(G^t)^T=G^T$;
\item for every $E'\subsetneq E_{G^t}$, $(V,E')^T\neq G^T$.
\end{itemize}
Such graph is called the transitive reduction of $G$ and it is given by:
\[
E_{G^t}:= \bigcap_{\{E'\subset V\times V \mid (V,E')^T=G^T\}} E'.
\]
\end{teo}

Given a finite poset $P$, we call $\Gcal(P)$ the transitive reduction of the associated DAG.

\subsection{From Persistent Sets to Dendrograms}
\label{sec:finite_pos}

We now establish a connection between functors, posets, and tree. To this end, we introduce a general procedure that applies beyond the specific case of functors $S: (\R,\leq)\rightarrow \Sets$, called \emph{persistent sets} \citep{carlsson2013classify, curry2018fiber, curry2021decorated}, which constitutes the primary focus of this work. Such objects are in fact of great interest in TDA and represent a (non-equivalent) alternative to persistence modules/diagrams when working with information collected from path-connected components of filtrations. A more thorough exploration of the broader scenarios is left for future research.

\begin{figure}
    \begin{subfigure}[c]{0.49\textwidth}
    	\centering
    	\includegraphics[width = \textwidth]{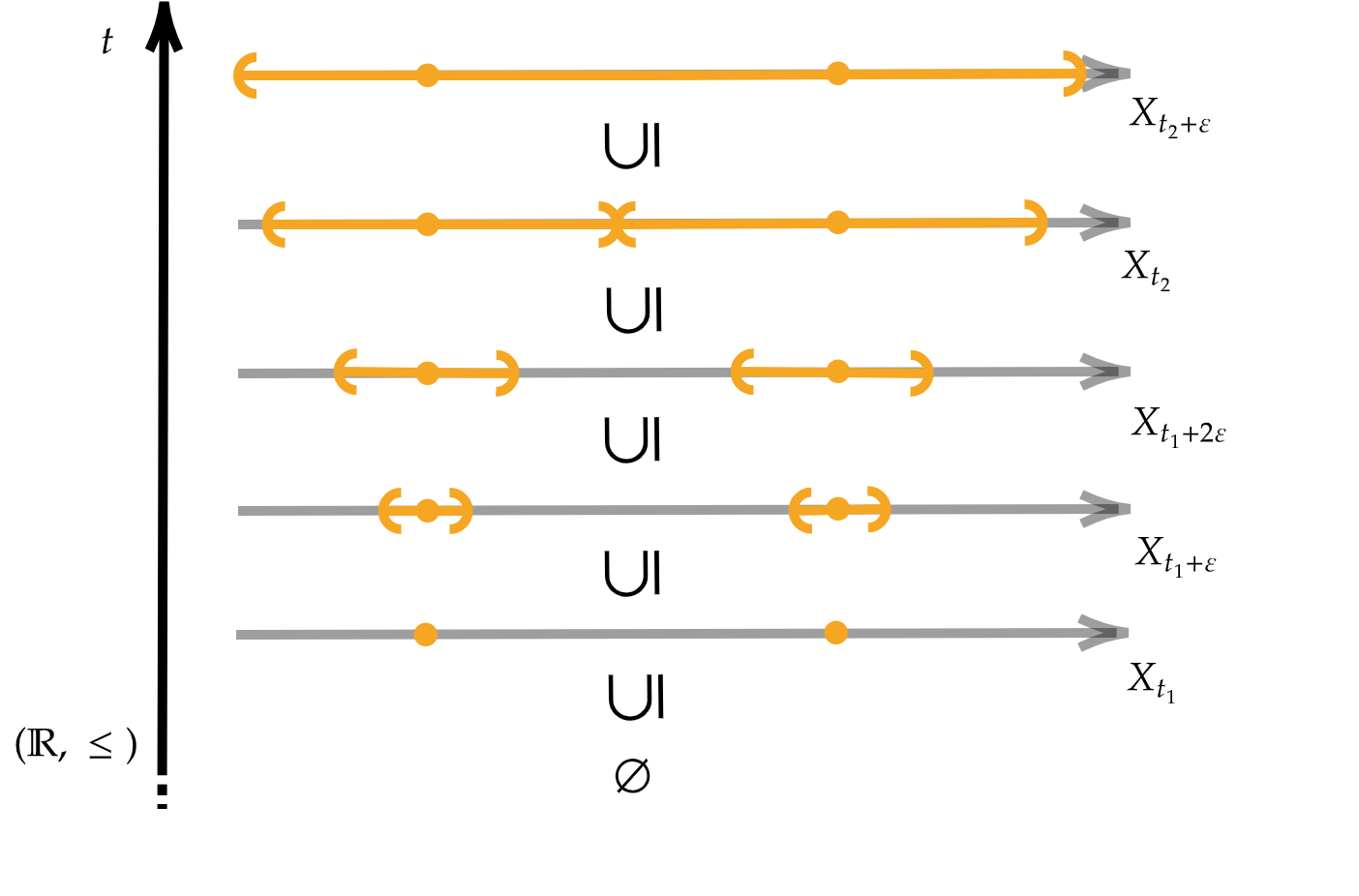}
    	\caption{A filtration $\X$.}
    	\label{fig:filtr}
	\end{subfigure}
%	\hspace{0.5 cm}
    \begin{subfigure}[c]{0.49\textwidth}
    	\centering
	    \includegraphics[width = \textwidth]{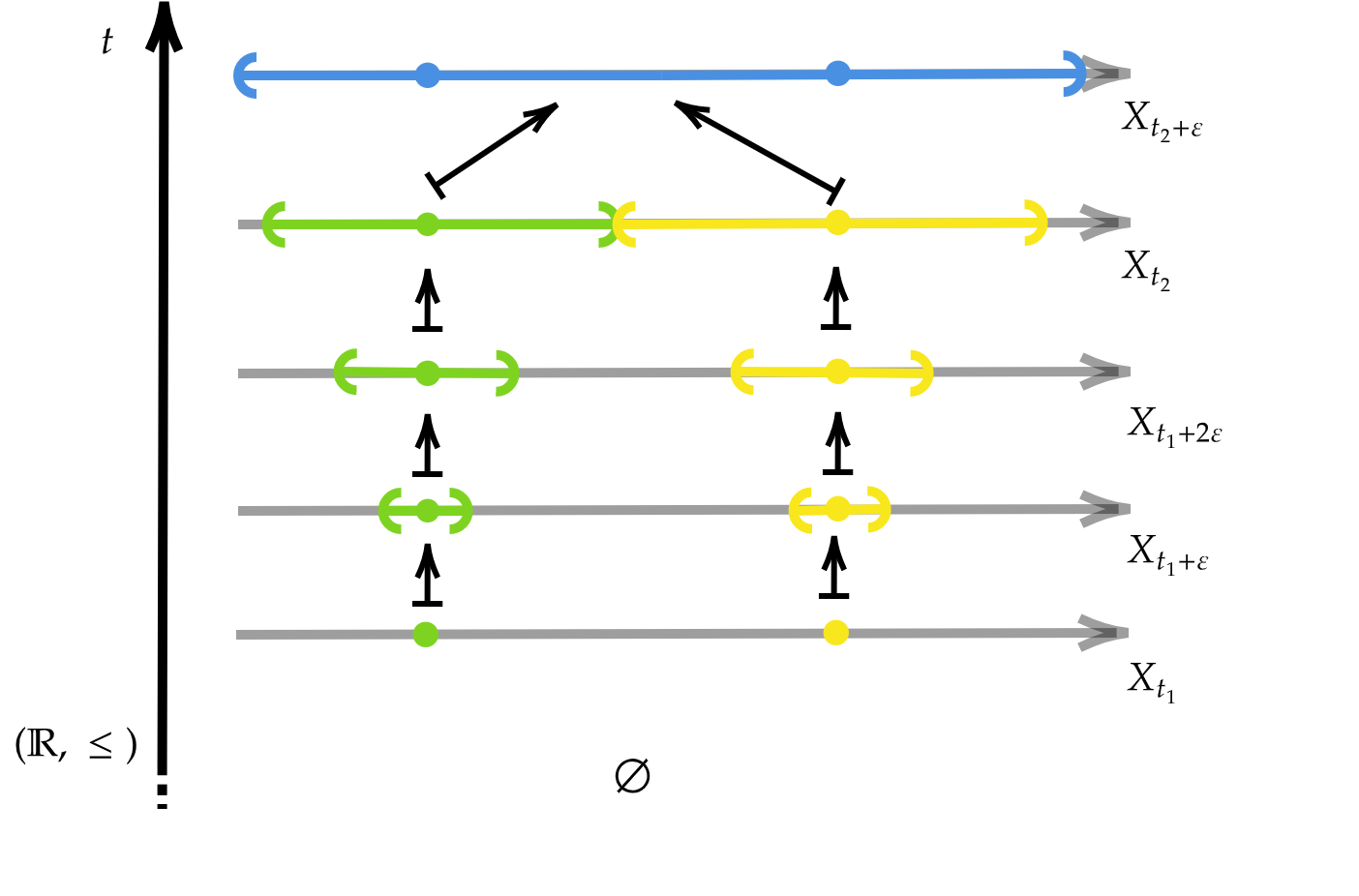}
    	\caption{The persistent set $\T$.}
    	\label{fig:abs_MT}
	\end{subfigure}
%	\hspace{0.5 cm}

    \begin{subfigure}[c]{0.49\textwidth}
    	\centering
	    \includegraphics[width = \textwidth]{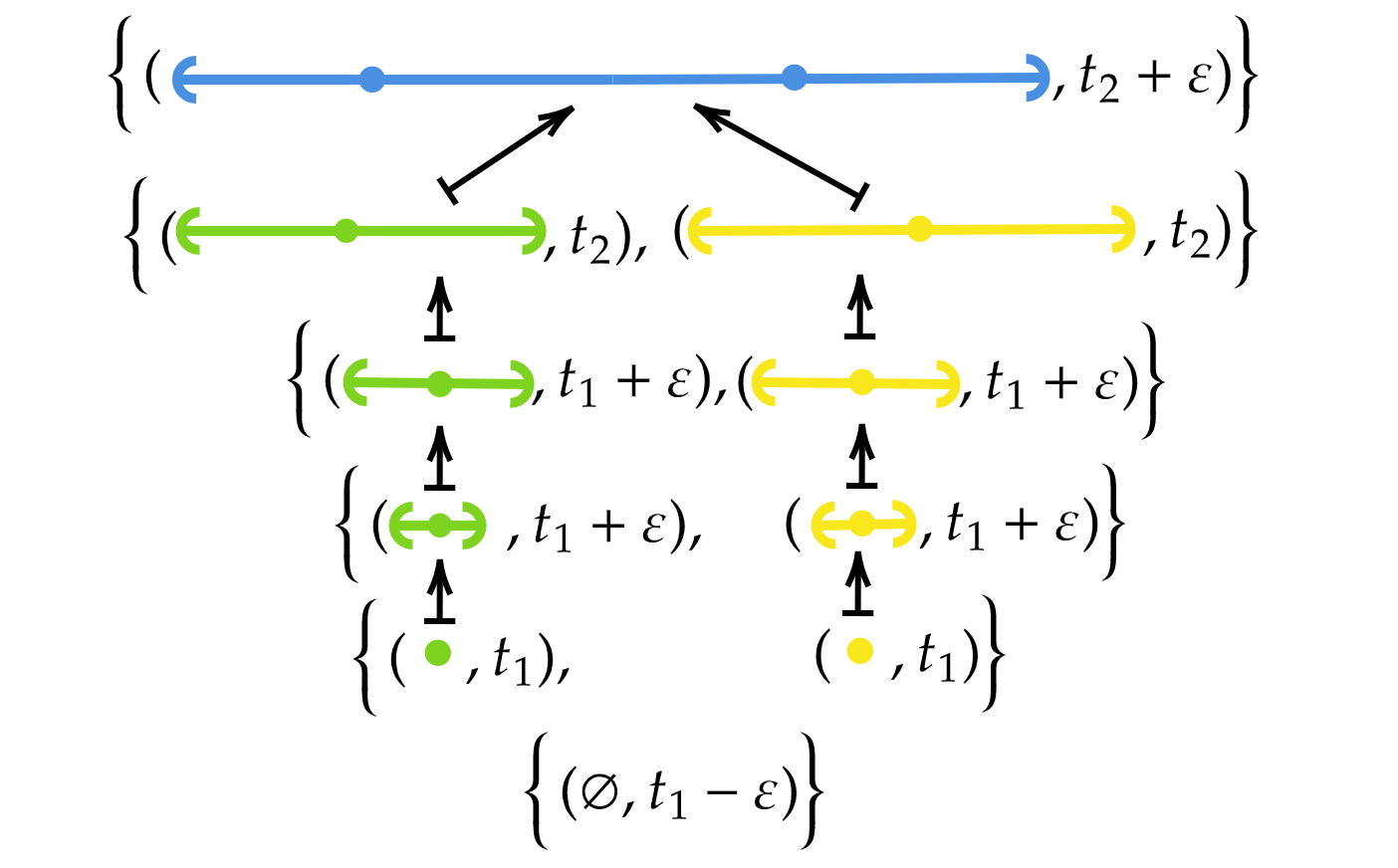}
    	\caption{The display poset $D_{\T}$.}
    	\label{fig:display}
	\end{subfigure}		
\begin{subfigure}[c]{0.49\textwidth}
    	\centering
	    \includegraphics[width = \textwidth]{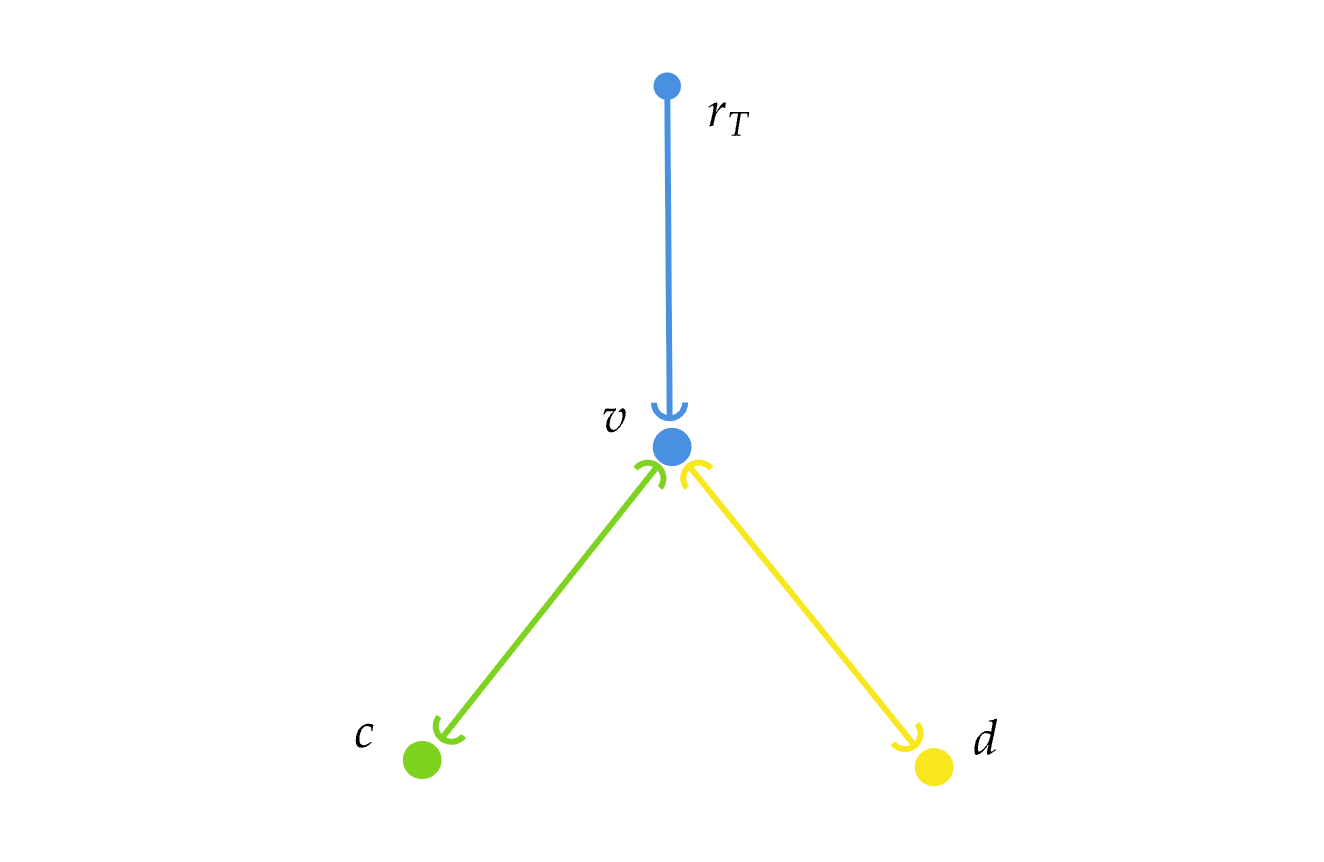}
    	\caption{The tree $\Gcal(D_{\T_{\mid C}})$.}
    	\label{fig:example}
	\end{subfigure}

\caption{An example of a filtration along with the persistent set of its path-connected components, the associated display poset and its graph representation. The colors are used throughout the plots to highlight the relations between the different objects.}
\label{fig:preliminary}
\end{figure}

The objects we introduce in the following are represented in \Cref{fig:preliminary}.

\begin{defi}[adapted from \cite{curry2021decorated}]
Given a functor $S:P\rightarrow \Sets$, with $P$ being a poset, the display poset of $S$ is defined as the set of pairs $D_{S}:=\{(p,s)\mid p\in P, s \in S(p)\}$. This is a poset with $(p,s)\leq (p',s')$ if and only if $p\leq p'$ and $S(p\leq p')(s)= s'$. For every display poset $D_{S}$ we have a monotone increasing projection
$h:D_{S}\rightarrow P$ defined by $(p,s)\mapsto p$. 
\end{defi}

The poset structure on $P$ and the one on $D_{S}$ are consistent with their respective graph representations, as made precise in the following proposition.

\begin{prop}
The map $h:D_{S}\rightarrow P$ induces a graph homomorphism (a function between vertices that also sends edges to edges) between $\Gcal(D_S)$ and $\Gcal(P)$.
\begin{proof}
By construction, any time we have $(p,s)\leq (p',s') \leq (p'',s'')$, this also implies $p\leq p'\leq p''$. Viceversa, any time we have $(p,s) \leq (p'',s'')$ and $p\leq p'\leq p''$, by functoriality there exists also $(p',s')$ such that 
$(p,s)\leq (p',s') \leq (p'',s'')$.
\end{proof}
\end{prop}

Now we restrict our attention to persistent sets.
 
\begin{defi}[\citep{patel2018generalized}]
A constructible persistent set is a functor $S: (\R,\leq)\rightarrow \Sets$ such that there is a finite collection of real numbers $\{t_1<t_2<\ldots<t_n\}$ such that:
\begin{itemize}
\item $S(t<t')=\emptyset$ for all $t<t_1$;
%\item $S(t)=\{\star\}$ for all $t>t_n$;
\item for $t,t'\in (t_i,t_{i+1})$ or $t,t'>t_n$, with $t<t'$, then $S(t<t')$ is bijective.
\end{itemize}
The set $\{t_1<t_2<\ldots<t_n\}$ is called critical set and $t_i$ are called critical values. If $S(t)$ is always a finite set, then $S$ is a finite persistent set.
\end{defi}

Note that we can always consider a minimal set of critical values \citep{pegoraro2024finitely}, which constitute the critical poset of the persistent set.
We also add a last requirement for the persistent sets we want to work with. 

\begin{defi}[\cite{pegoraro2024finitely}]
A regular persistent set $S$ is a finite constructible persistent set such that, for every critical value $t_i$ and for every $\varepsilon>0$ small enough, $S(t_i\leq t_i+\varepsilon)$ is a bijective map. 
\end{defi}

At this point we have the following immediate result, which we state without proof.

\begin{prop}
Given two regular persistent sets $S$ and $S'$, then $S\simeq S'$ if and only if:
\begin{itemize}
\item they have the same minimal set of critical values $C = \{t_1,\ldots,t_n\}$ 
\item $S_{\mid C}\simeq S'_{\mid C}$, with $S_{\mid C}:C\hookrightarrow \R\rightarrow  \Sets$.
\end{itemize}
\end{prop}

Without going into the details of the construction of a merge tree, we just show that tree-structures are the natural discrete tool to work with persistent sets. See also \Cref{rmk:merge_trees}.

\begin{prop}\label{prop:graphs_trees}
Consider a finite regular persistence set $S$ and a set of critical values $C$. If $S(t)$ has cardinality $1$ for $t$ big enough, $\Gcal(D_{S_{\mid C}})$ is a tree.
\end{prop}

\begin{rmk}\label{rmk:merge_trees}
    In the proof of \Cref{prop:graphs_trees}, we actually show that the graph $\Gcal(D_{S_{\mid C}})$ is equivalent up to degree $2$ vertices (as a directed graph) to the merge tree of $S$, upon removing from the merge tree the edge that goes to infinity. In fact, proving the actual statement of \Cref{prop:graphs_trees} reduces to observing that, within the setting considered, $\Gcal(D_{S_{\mid C}})$ is a rooted directed acyclic graph, with the $parent >child$ relationship being induced by the poset relationship on $D_{S_{\mid C}}$. 
    Note that this proposition further justifies our idea of working up to degree $2$ vertices.
\end{rmk}

\section{Edits of Dendrograms}
\label{sec:edits}

Now we present the main contributions of the paper. 

The approach we follow is to define a distance inspired by Edit Distances \citep{Tai}, but with substantial differences in the edit operations.
The philosophy of these distances is to allow certain modifications of the base object, called \emph{edits}, each being associated to a \emph{cost}, and to define the distance between two objects as the minimal cost that is needed to transform the first object into the second with a finite sequence of edits.
In this way, up to properly setting up a set of edits, one can formalize the deformation of a dendrogram modifying the local information induced by the weight function defined on the edges and the tree structure itself.  
On top of that, edit distances frequently enjoy some decomposition properties which simplify the calculations \citep{TED}, which are notoriously very heavy \citep{np_hard}. 
We also point out that, in literature, tree are often considered with weights being real numbers and avoid modeling the weight space as we do with the definition of editable spaces.

Given a space of dendrograms $(\mathcal{T}, X)$, with $(X,\odot,0_X)$ editable space, and given a dendrogram $(T,\varphi)$, we define our edits as operators $\{(T,\varphi)\}\rightarrow (\mathcal{T}, X)$. So, an edit operation is not defined on the whole $(\mathcal{T}, X)$, but only on $(T,\varphi)$, and its image can be defined up to isomorphism. Moreover, given two edits $e_i:\{(T,\varphi)\}\rightarrow (\mathcal{T}, X)$ and $e_j:e_i(\{(T,\varphi)\})\rightarrow (\mathcal{T}, X)$, we can consider their composition $e_i \circ e_j$. Any finite composition of edits is referred to as an \emph{edit path}.

The distance $d_E$ we introduce differs substantially from previously defined edit distances, as it is tailored specifically for comparing topological summaries—roughly meaning that topologically irrelevant points can be removed from a tree at no cost. 

Recall that, via the identification $E_T\cong V_T-\{r_T\}$, the upcoming transformations can be described equivalently in terms of either edges or vertices.

\begin{itemize}

\item {\em Shrinking} an edge means changing the weight value of the edge with a non zero element in $X$.
The inverse of this transformation is the shrinking which restores the original edge weight. 

\item {\em Deleting} an edge $(v_1,v_2)$ results into a new tree, with the same vertices apart from $v_1$ (the lower one), and with the parent  of the deleted vertex which gains all of its children.
%Since $E_T$ can be identified with $V_T-\{r_T\}$, sometimes we might refer to this edit also as the deletion of the vertex $v_1$, which means deleting the edge $(v_1,v_2)$.
%As mentioned before, we might also refer to this edit as the deletion of the vertex $v_1$, which indeed means deleting the edge between \(v_1\) and its parent .

The inverse of deletion is the {\em insertion} of an edge along with its child vertex.
We can insert an edge at a vertex \(v\) specifying the child of \(v\) and its children (that can be either none or any portion of the children of \(v\)) and the weight of the edge.

\item Lastly, we can eliminate a degree two vertex \(v\), that is, a parent  with an only child, connecting the two adjacent edges which arrive and depart from \(v.\) The weight of the resulting edge is the sum of the weights of the joined edges.
This transformation is the {\em ghosting} of the vertex $v$.
Its inverse transformation is called the {\em splitting} of an edge.
\end{itemize}

\begin{figure}[h!]
    \centering
	\begin{subfigure}[c]{0.49\textwidth}
    	\centering
    	\includegraphics[width = \textwidth]{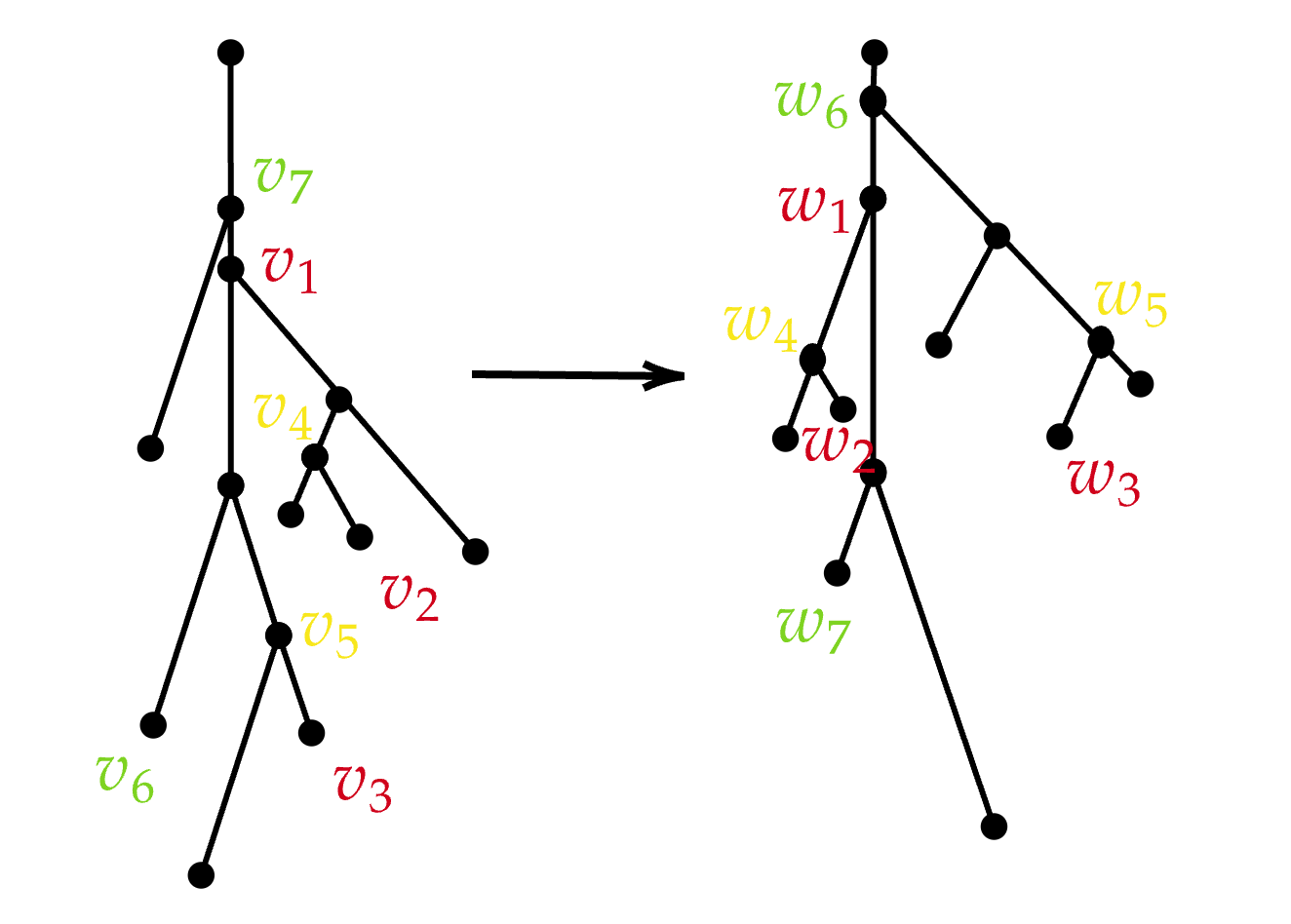}
    	\caption{Starting and target weighted trees, with highlighted the vertices involved in the edit path.}
    	\label{fig:recap}
    \end{subfigure}
    \centering
   	\begin{subfigure}[c]{0.49\textwidth}
    	\centering
    	\includegraphics[width = \textwidth]{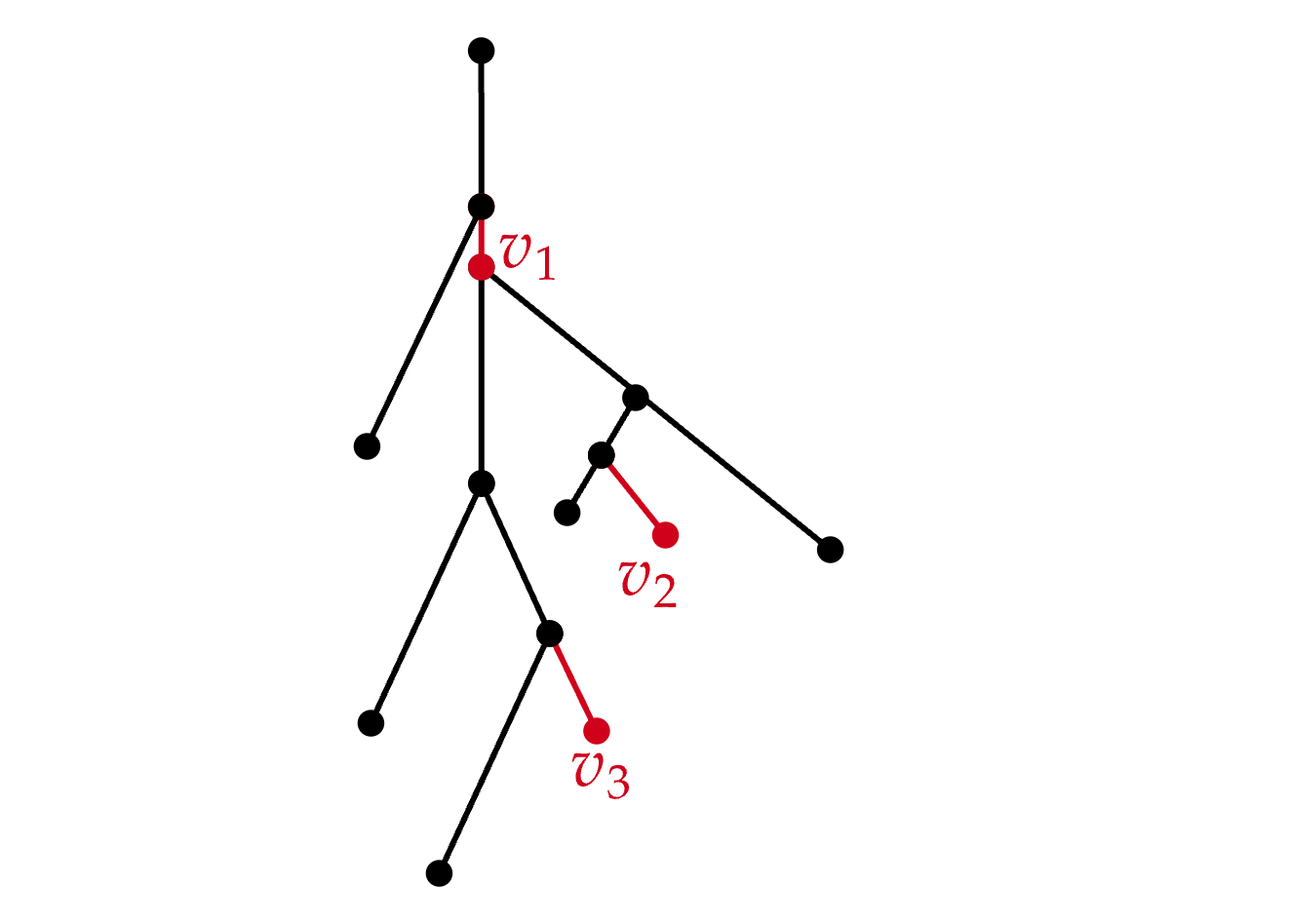}
    	\caption{Deletions of $v_1,v_2,v_3$, in red.}
    	\label{fig:deletions}
    \end{subfigure}
%    \hfill

    \centering
	\begin{subfigure}[c]{0.49\textwidth}
		\centering
		\includegraphics[width = \textwidth]{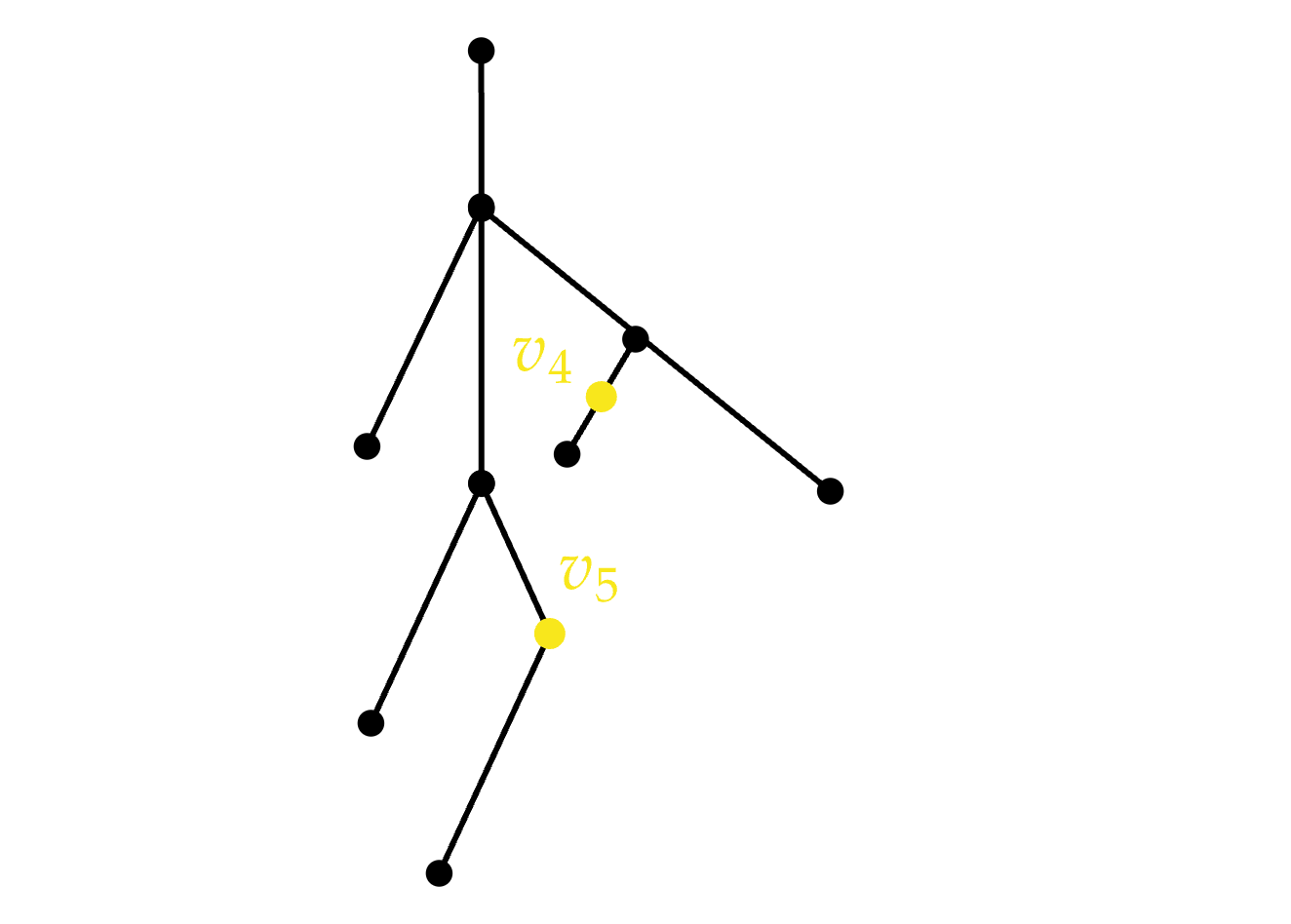}
		\caption{Ghostings of $v_4,v_5$, in yellow.}
		\label{fig:ghostings}
	\end{subfigure}
    \centering
	\begin{subfigure}[c]{0.49\textwidth}
    	\centering
    	\includegraphics[width = \textwidth]{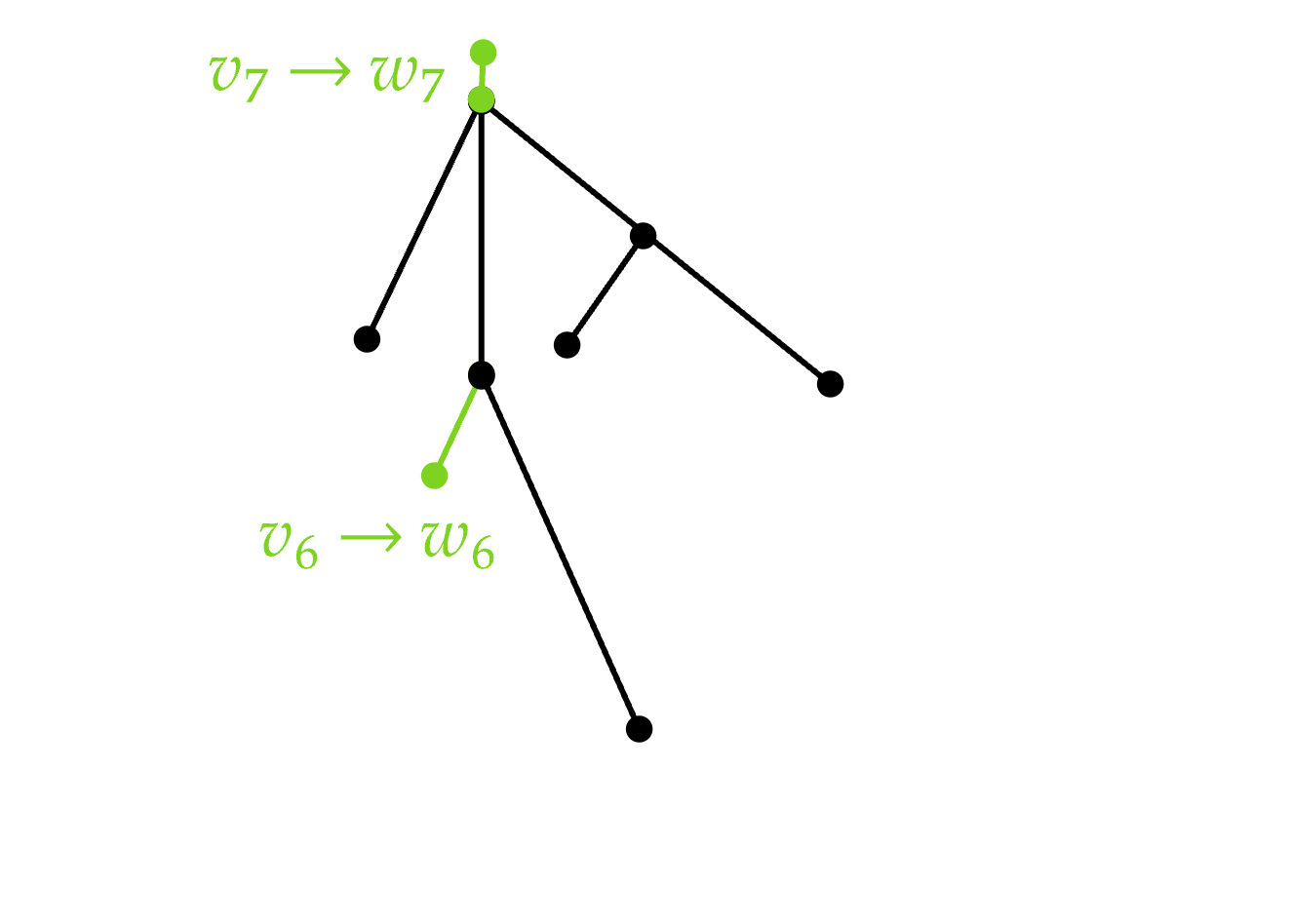}
    	\caption{Shrinkings of $v_6$ and $v_7$ to match the weight, respectively, of $w_6$ and $w_7$, in green.}
    	\label{fig:shrinkings}
    \end{subfigure}

    \centering
    \begin{subfigure}[c]{0.49\textwidth}
    	\centering
    	\includegraphics[width = \textwidth]{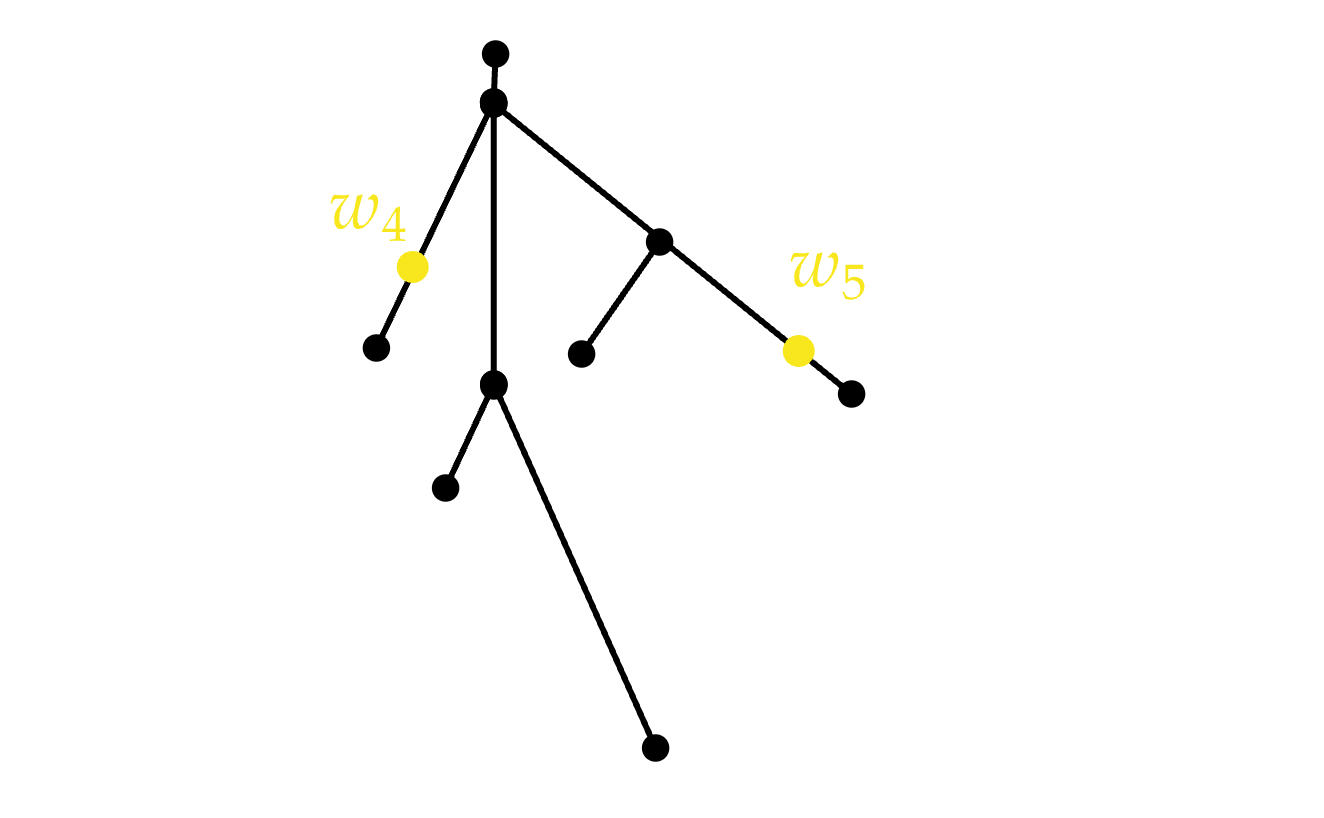}
    	\caption{Splitting edges via $w_4$ and $w_5$, in yellow.}
    	\label{fig:splittings}
    \end{subfigure}
    \centering
    \begin{subfigure}[c]{0.49\textwidth}
    	\centering
    	\includegraphics[width = \textwidth]{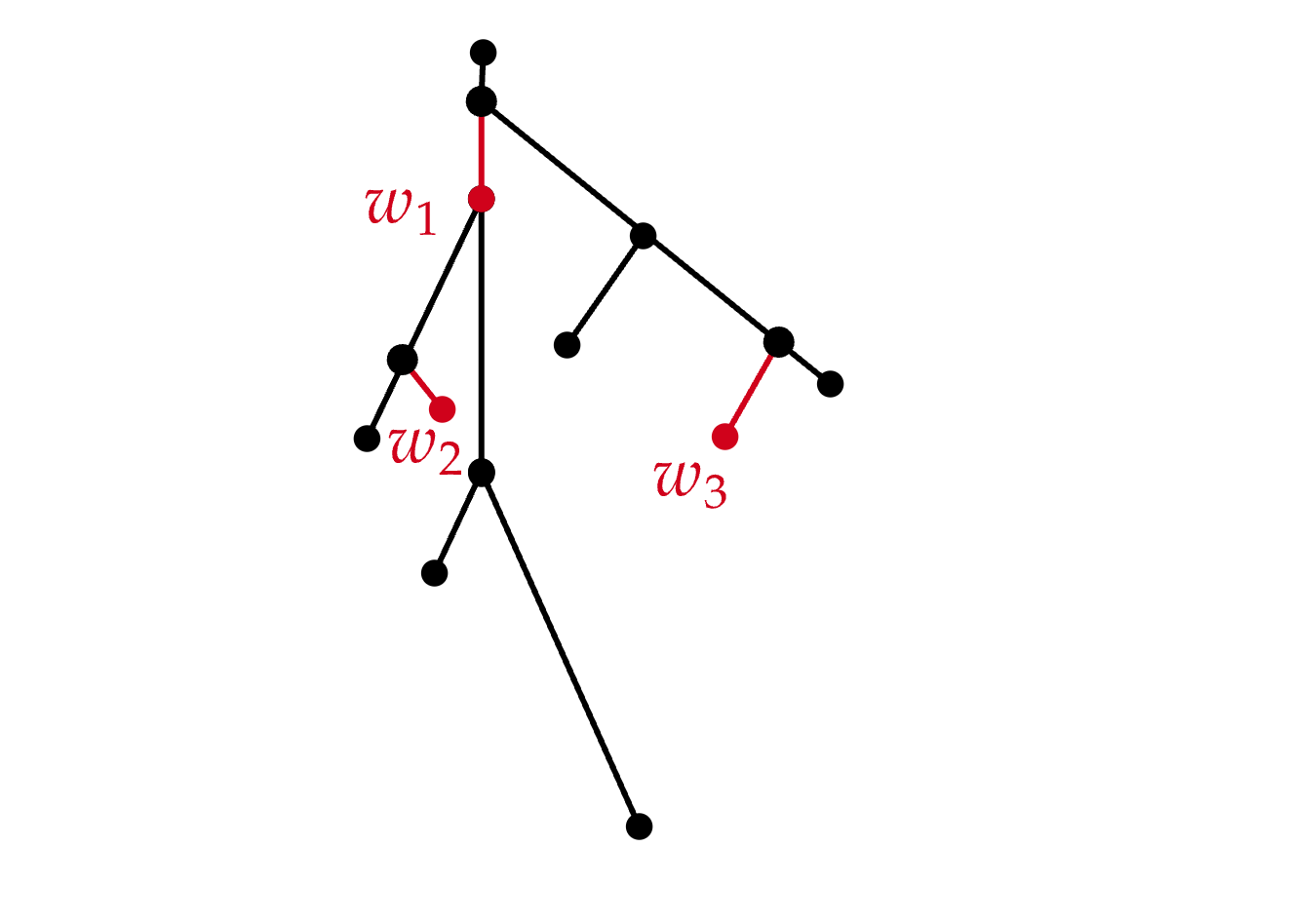}
    	\caption{Insertions of $w_1,w_2,w_3$, in red.}
    	\label{fig:insertions}
    \end{subfigure}
    \hfill
\caption{(b)$\rightarrow$(e) form an edit path made from the left weighted tree in \Cref{fig:recap} to the right one. The edit path can be represented with a mapping - \Cref{sec:mappings} - consisting of the 
pairs $(v_i,\mathfrak{D})$ for all the red vertices in \Cref{fig:deletions}, $(v_i,\mathfrak{G})$ for all the yellow vertices in \Cref{fig:ghostings}, $(v_i,w_i)$ for all the vertices associated via the green color in \Cref{fig:shrinkings}, $(\mathfrak{G},w_i)$ for all the yellow vertices in \Cref{fig:ghostings} and $(\mathfrak{D},w_i)$ for all the red vertices in \Cref{fig:insertions}. }
\label{fig:edits_TDA}
\end{figure}

  See \Cref{fig:edits_TDA} for an example of an edit path between two graphs whose weights are edges lenghts - \Cref{sec:positive_editable}.

\begin{rmk}\label{rmk:classical_edit}
Shrinking and deletions are classical edits for trees, with shrinking being usually referred to as \emph{relabeling}. The ghosting edit (and its inverse) is completely unusual
and we are not aware of any previous work employing it. Its definition is largely driven by the TDA perspective presented in the Introduction and it is fundamental for the stability results contained in \cite{pegoraro2024functional} and  \cite{pegoraro2024finitely}.
Even in the context of graph edit distances (GEDs)
\citep{ambauen2003graph, gao2010survey, zeng2009comparing, lerouge2016exact, serratosa2021redefining},
we are not aware of any existing work that adopts definitions similar to ours. While vertex splitting and merging in GEDs bear some resemblance to our notions of insertion and deletion, there no correspondence for splittings and ghostings as we define them.
\end{rmk}

We also give the following definition.

\begin{defi}\label{defi:equiv_order_2}
Dendrograms are equal up to degree $2$ vertices if they become isomorphic after applying a finite number of ghostings or splittings. We write $(T,\varphi_T)\cong_2(T',\varphi_{T'})$. 
\end{defi} 

\Cref{defi:equiv_order_2} induces an equivalence relation  which identifies the set of dendrograms inside $(\mathcal{T},X)$ that we want to treat as equal.
We call $(\mathcal{T}_2,X)$ the space of equivalence classes of dendrograms in $(\mathcal{T},X)$, equal up to degree $2$ vertices.

\subsection{Costs of Edit Operations}

Now we associate to every edit a cost, that is, a length measure in the space $(\mathcal{T},X)$. In light of this interpretation, we will often use the words \virgolette{length} and \virgolette{cost} interchangeably when referring to an edit path.

The costs of the edit operations are defined as follows:
\begin{itemize}
\item if, via shrinking, an edge goes from weight $x$ to weight $y$, then the cost of such operation is $d_X(x,y)$;
\item for any deletion/insertion of an edge with weight $x$, the cost is equal to $d_X(x,0_X)$;
\item the cost of ghosting/splitting operations is always zero. This is coherent with the properties of editable spaces: consider a ghosting merging two adjacent edges separated by a degree two vertex, with $x$ being the weight of one edge, and $y$ the cost of the other. The total weight of the two edges is the same before and after the ghosting:
$\mid d_X(x\odot y,0_X)-d_X(x,0_X)-d_X(y,0_X)\mid =0$.
\end{itemize}

%The cost of an edit path $e_0\circ\ldots\circ e_n(T)$ is given by the sum of the costs of the single edit operations.

\begin{defi}
Given two dendrograms $T$ and $T'$ in $(\mathcal{T},X)$,  define:
\begin{itemize}
\item $\Gamma(T,T')$ as the set of all finite edit paths between $T$ and $T'$;
\item $cost(\gamma)$ as the sum of the costs of the edits for any $\gamma\in\Gamma(T,T')$;
\item the multigraph edit distance as:
\[
d_E(T,T')=\inf_{\gamma\in\Gamma(T,T')} cost(\gamma)
\]
\end{itemize}
\end{defi}

By definition the triangle inequality and symmetry must hold, but, up to now, this edit distance is intractable; one would have to search for all the possible finite edit paths which connect two dendrograms in order to find the minimal ones. On top of that, having an edit which is completely \virgolette{for free}, it is not even obvious  that $d_E(T,T')>0$ for some dendrograms.
However, it is clear that $d_E$ induces a pseudo-metric on classes of dendrograms up to degree two vertices.

%\begin{defi}\label{defi:top_stable}
%A pseudo-metric on $(\mathcal{T},W)$ which induces a non trivial pseudo-metric on  $(\mathcal{T}_2,W)$ is called topologically stable.  
%\end{defi}
%
%In other words a topologically stable pseudo-metric for dendrograms is a (non trivial) pseudo-metric which identifies dendrograms which are equivalent up to degree $2$ vertices.
%The discussion on such idea continues in  
%\Cref{sec:tree_edit} after we have properly established a (pseudo) metric on dendrogram spaces.

\subsection{Mappings}
\label{sec:mappings}

We introduce a combinatorial tool, called a \emph{mapping}, which plays a central role in developing the theoretical properties of $d_E$ and in bringing it into the realm of computability. In particular, we prove that mappings can be used to parametrize a finite set of edit paths, with distinct and interpretable characteristics, which is guaranteed to contain at least an optimal edit path between two dendrograms.

\begin{defi}
A mapping between two dendrograms $T$ and $T'$ is a set 
$M\subset (E_T \cup \{\mathfrak{D},\mathfrak{G}\})\times (E_{T'} \cup \{\mathfrak{D},\mathfrak{G}\})$ 
satisfying:

\begin{itemize}
\item[(M1)] consider the projection of the Cartesian product $(E_T \cup \{\mathfrak{D},\mathfrak{G}\})\times (E_{T'} \cup \{\mathfrak{D},\mathfrak{G}\})\rightarrow (E_{T} \cup \{\mathfrak{D},\mathfrak{G}\})$; we can restrict this map to $M$ obtaining $\pi_T:M\rightarrow (E_T \cup \{\mathfrak{D},\mathfrak{G}\})$. The maps $\pi_T$ and $\pi_{T'}$ are surjective on $E_T$ and 
$E_{T'}$, i.e. $E_T\subset \im(\pi_T)$ and $E_{T'}\subset \im(\pi_{T'})$;
\item[(M2)]$\pi_T$ and $\pi_{T'}$ are injective on $M\cap (E_T\times E_{T'})$;
\item[(M3)]  given $(a,b)$ and $(c,d)\in M\cap (V_T\times V_{T'})$, $a>c$, if and only if $b>d$;
\item[(M4)] if $(a,\mathfrak{G})\in M$ (or analogously $(\mathfrak{G},a)$), 
then after applying all deletions of the form $(v,\mathfrak{D})\in M$, the vertex $a$ becomes a degree $2$ vertex. In other words: let $child(a)=\{b_1,..,b_n\}$. Then there is exactly one $i$ such that for all $j \neq i$, for all $v \in V_{sub(b_j)}$, we have $(v,\mathfrak{D})\in M$; and there is one and only one $c$ such that $c=\max\{x<b_i\mid (x,y)\in M$ for any $y \in V_{T'}\}$. 
\end{itemize}  

We call $\Mapp(T,T')$ the set of all mappings between $T$ and $T'$. 
\end{defi}

We can always have a mapping made by:
\[
\{(e,\mathfrak{D}) \mid e \in E_T\} \cup \{(\mathfrak{D},e) \mid e \in E_{T'}\}.
\]
In the following, we will see that such a mapping corresponds to the edit path which deletes all the edges of one dendrogram and inserts all the edges of the second.

We may refer to edges which appear in the pairs in $M\cap (V_T\times V_{T'})$ as the \emph{paired} or \emph{matched} edges/vertices.

\begin{rmk}
In literature usually mappings are subsets of $E_T\times E_{T'}$ and their properties are equivalent to
(M2) and (M3) (see \cite{TED, merge_farlocca} and references therein).
In this way, mappings coherently maps edges into edges, respecting the tree ordering, inducing in this way the needed \emph{relabeling} operations - see also \Cref{rmk:classical_edit}. 
The differences with our definition are caused by the introduction of the ghosting edits. Accordingly, the way in which mappings parametrize edit paths in the present work is completely novel and, for instance, it establishes some (partial) ordering between the edits.
Note that the set of mappings we consider is bigger than in usual edit distances - thus computing the distance has a higher computational cost. 
\end{rmk}

Every $M\in \Mapp(T,T')$ parametrizes a set of edit paths, with identical cost, as follows:
\begin{itemize}
\item $\gamma_{d}^T$ is made by the deletions to be done on $T$, that is, the pairs $(v,\mathfrak{D})$, executed in any order. So we obtain $T^M_d=\gamma_{d}^T(T)$, which is well defined and does not depend on the order of the deletions.
Similarly, we define $\gamma_{d}^{T'}$ as a path made by the deletions to be done on $T'$, that is, the pairs $(\mathfrak{D},w)$, executed in any order, and obtain $T'^M_d=\gamma_{d}^{T'}(T')$. 
\item One then proceeds editing $T^M_d$ by ghosting all the vertices $(v,\mathfrak{G})$ in $M$, in any order, getting a path $\gamma^T_g$ and the dendrogram $T_M:= \gamma^T_g\circ \gamma^T_d (T)$. As before, we can do an analogous procedure on $T'^M_d$, ghosting all the vertices $(\mathfrak{G},w)$ in $M$, in any order, and building a path $\gamma^{T'}_g$, along with the dendrogram $T'_M:= \gamma^{T'}_g\circ \gamma^{T'}_d (T')$.
\item Since all the remaining points in $M$ are paired, the two dendrograms $T'_M$ and $T_M$ must be isomorphic as tree structures. This is guaranteed by the properties of $M$. So one can shrink $T_M$ onto $T'_M$, and the composition of the shrinkings, executed in any order, gives an edit path $\gamma_s^T$.
\end{itemize}

By construction $\gamma_s^T\circ \gamma_g^T\circ \gamma_d^T(T)=T'_M$,
and $(\gamma_d^{T'})^{-1}\circ (\gamma_g^{T'})^{-1}\circ \gamma_s^T\circ \gamma_g^T\circ \gamma_d^T (T)=T'$.
Where the inverse of an edit path is thought as the composition of the inverses of the single edit operations, taken in the inverse order. 
 
We call $\gamma_M$ the set of all possible edit paths of the form $(\gamma_d^{T'})^{-1}\circ (\gamma_g^{T'})^{-1}\circ \gamma_s^T\circ \gamma_g^T\circ \gamma_d^T$, 
obtained by changing the order in which the edit operations are executed inside $\gamma_d$, $\gamma_g$ and $\gamma_s$.
Even if $\gamma_M$ is a set of paths, its cost is well defined:

\[
cost(M):= cost(\gamma_M)=cost(\gamma_d^T)+cost(\gamma_s^T)+cost(\gamma_d^{T'}).
\]

We prove that there exists always a mapping
$M$ such that the paths in $\gamma_M$ are optimal.

\begin{teo}[Main Theorem]
\label{teo:main_thm}
Given two dendrograms $T$ and $T'$, for every finite edit path $\gamma$, there is a mapping $M \in \Mapp(T,T')$ such that $cost(M)\leq cost(\gamma)$.
\end{teo}

Two corollaries follow immediately from the facts that \( \Mapp(T, T’) \) is finite and that any mapping with cost zero must consist entirely of ghostings and splittings.

\begin{cor}\label{cor:finite_mapp}
 For any two dendrograms $T$ and $T'$, the set $\Mapp(T,T')$ is finite. Thus, we have the following well defined pseudo-metric between dendrograms:
\[
d_E(T,T')=\inf \{cost(\gamma)\mid\gamma\in\Gamma(T,T')\}=\min\{cost(M)\mid M\in \Mapp(T,T')\}
\]
which we will refer to as the edit distance between $T$ and $T'$.
\end{cor}

\begin{cor} 
Given $T$ and $T'$ dendrograms, $d_E(T,T')=0$ if and only if $T$ and $T'$ are equal up to degree $2$ vertices. In other words $d_E$ is a metric for dendrograms considered up to degree $2$ vertices.
\end{cor}

To conclude this section, we define a particular family of mappings which help us in further restricting the search space for optimal edit paths.

\begin{defi}\label{defi:M_2_mapp}
A mapping $M\in \Mapp(T,T')$ has maximal ghostings if the following hold: $(v,\mathfrak{G})\in M$ if, and only if, $v$ is of degree $2$ after the deletions in $T$ and, similarly, $(\mathfrak{G},w)\in M$ if, and only if, $w$ is of degree $2$ after the deletions in $T'$.

A mapping $M\in \Mapp(T,T')$ has minimal deletions if the following hold: $(v,\mathfrak{D})\in M$ implies that neither $v$ nor $parent(v)$ are of degree $2$ after applying all the other deletions in $T$ and, similarly, $(\mathfrak{D},w)\in M$ implies that neither $w$ nor $parent(w)$ are of degree $2$ after applying all the other deletions in $T'$.

We collect all mappings with maximal ghostings 
in the set $M_{\Gcal}(T,T')$
and paths with minimal deletions in the set $M_{\mathcal{D}}(T,T')$. Lastly, we set: $M_2(T,T'):=M_{\Gcal}(T,T')\bigcap  M_{\mathcal{D}}(T,T')$.
\end{defi}

In other words with mappings in $M_2(T,T')$ we are always eliminating all the degree $2$ vertices which arise from deletions and we are not deleting edges which we can shrink. The following lemma then applies.  

\begin{lem}\label{lemma:M_2}
\[
\min\{cost(M)\mid  M\in \Mapp(T,T') \}=\min\{cost(M)\mid M\in M_2(T,T')\}.
\]
\end{lem}

\section{A Discussion on the Use of Different Editable Spaces}
\label{sec:examples}

In this section we discuss some key examples which regard dendrograms with values in different editable spaces.

We stress again that this is to be intended as a roadmap, the outlining of a general idea and not as a formal and coherent description on how one can treat persistent sets via dendrograms. 
Coherently, we explicitly point out which questions and problems should be assessed by future works.

\subsection{Function and Vector Valued Weights}
\label{sec:dec_MT}

We make use of the following notation: $\Y:\R\rightarrow \Top$ is a filtration of topological spaces. 
That is, $\Y$ is a functor such that $(r\leq r')\mapsto Y_{r\leq r'}:Y_r \hookrightarrow Y_{r'}$. The functor $\pi_0:\Top \rightarrow \Sets$, instead, is the functor of path connected components. We now consider persistent sets of the form $\G$. Note that every constructible persistent set can be obtained as $\G$, using the discrete topology.

We enrich a persistent set $\pi_0(\Y)$ with local information regarding $\Y$. Using functions 
$f:D_{\pi_0(\Y)}\rightarrow X$
with values in an editable space $X$.
 In particular, we propose:
\begin{enumerate}
\item diagram or Betti numbers enriched persistent sets, on the same line of thought of \cite{curry2021decorated, curry2023stability};
\item persistent sets enriched with the measure of sublevel sets. 
\end{enumerate}

Once a function $f:D_{\pi_0(\Y)}\rightarrow X$ is obtained, one then needs to obtain a dendrogram. As we explain in \Cref{sec:problems}, this step requires some attention, and thus we leave its investigation to separate works. For the sake of this discussion, we just mention that, qualitatively, this \virgolette{discretization} combines taking $\Gcal(D_{\G})$ with employing as weights of the dendrograms the restriction of $f:D_{\pi_0(\Y)}\rightarrow X$ to the edges of $\Gcal(D_{\G})$.

We illustrate this with two examples: first we consider the datum of a point cloud a then a function defined on a open subset of $\R^n$.

\begin{figure}
	\begin{subfigure}[c]{0.49\textwidth}
		\centering
		\includegraphics[width = \textwidth]{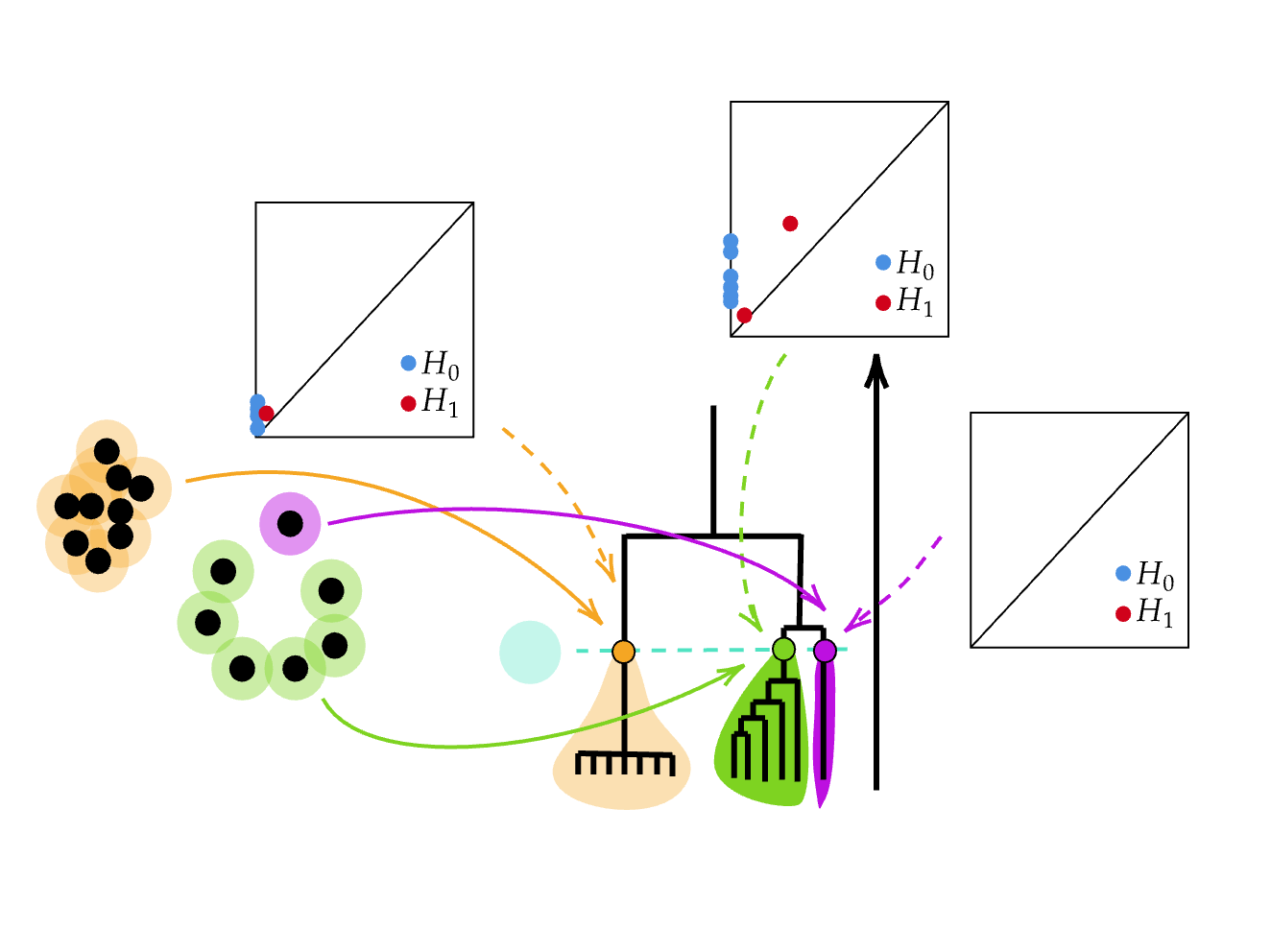}
            \captionsetup{singlelinecheck=off, margin={0.3cm, 0.1cm}}
		\caption{A point cloud (left), with the thickening giving the off-set filtration at some particular value (dashed horizontal cyan line) and the display poset of (the persistent set of) the filtration (right). At the chosen value the filtration has three path connected components, shown in different colours, each with different shapes, characterized by the (qualitative) associated persistence diagram. }
		\label{fig:PD_w_0}
	\end{subfigure}
	\begin{subfigure}[c]{0.49\textwidth}
		\centering
		\includegraphics[width = \textwidth]{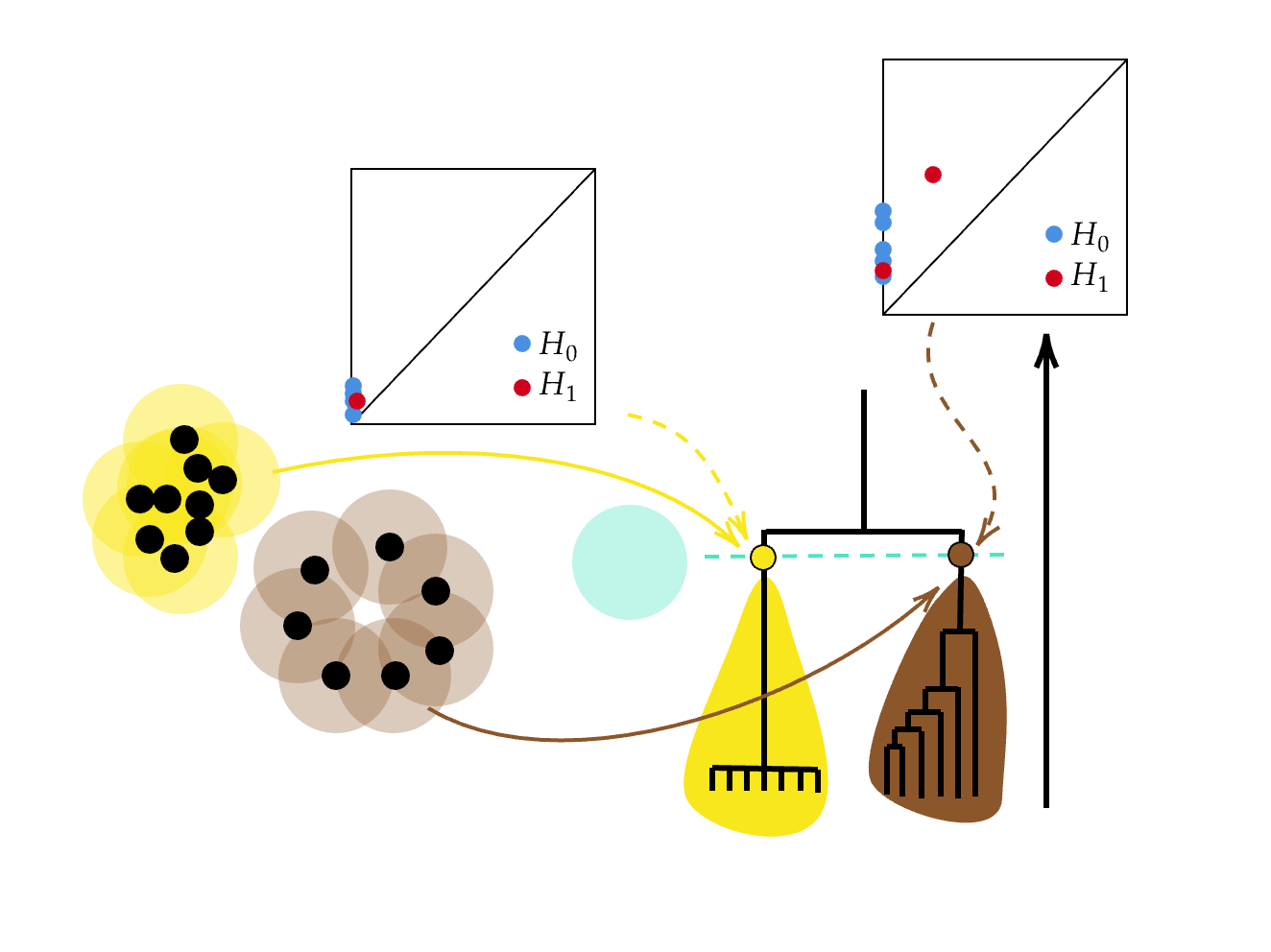}
            \captionsetup{singlelinecheck=off, margin={0.3cm, 0.1cm}}
		\caption{Same situation as in \Cref{fig:PD_w_0}, but at a different filtration value. At the chosen value - bigger than in \Cref{fig:PD_w_0}, the filtration has two path connected components, shown in different colours, each with different shapes, characterized by the (qualitative) associated persistence diagram. At this resolution, the brown cluster shows a more persistent ringed shape.}
		\label{fig:PD_w_1}
	\end{subfigure}

	\begin{subfigure}[c]{0.49\textwidth}
		\centering
		\includegraphics[width = \textwidth]{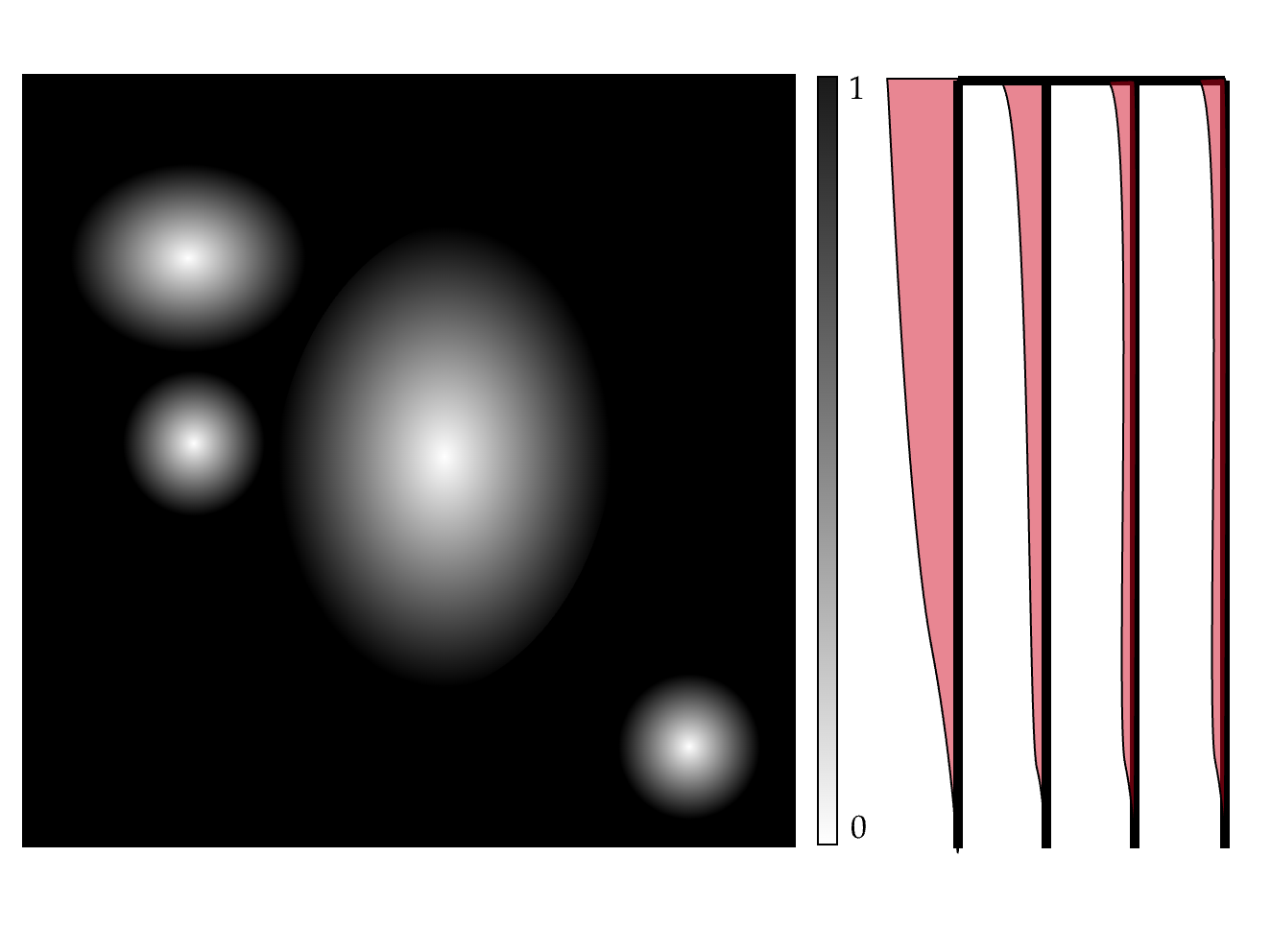}
            \captionsetup{singlelinecheck=off, margin={0.3cm, 0.1cm}}
		\caption{A real valued function defined on a square with four distinct local minima, associated to downward peaks with different shapes. The associated display poset (of the persistent set) accordingly shows four branches (black) and, in red, the (qualitative) volume of the associated downward peaks included in the sublevel sets, which increases with the value of the filtration.}
		\label{fig:vol_w_0}
	\end{subfigure}
	\begin{subfigure}[c]{0.49\textwidth}
		\centering
		\includegraphics[width = \textwidth]{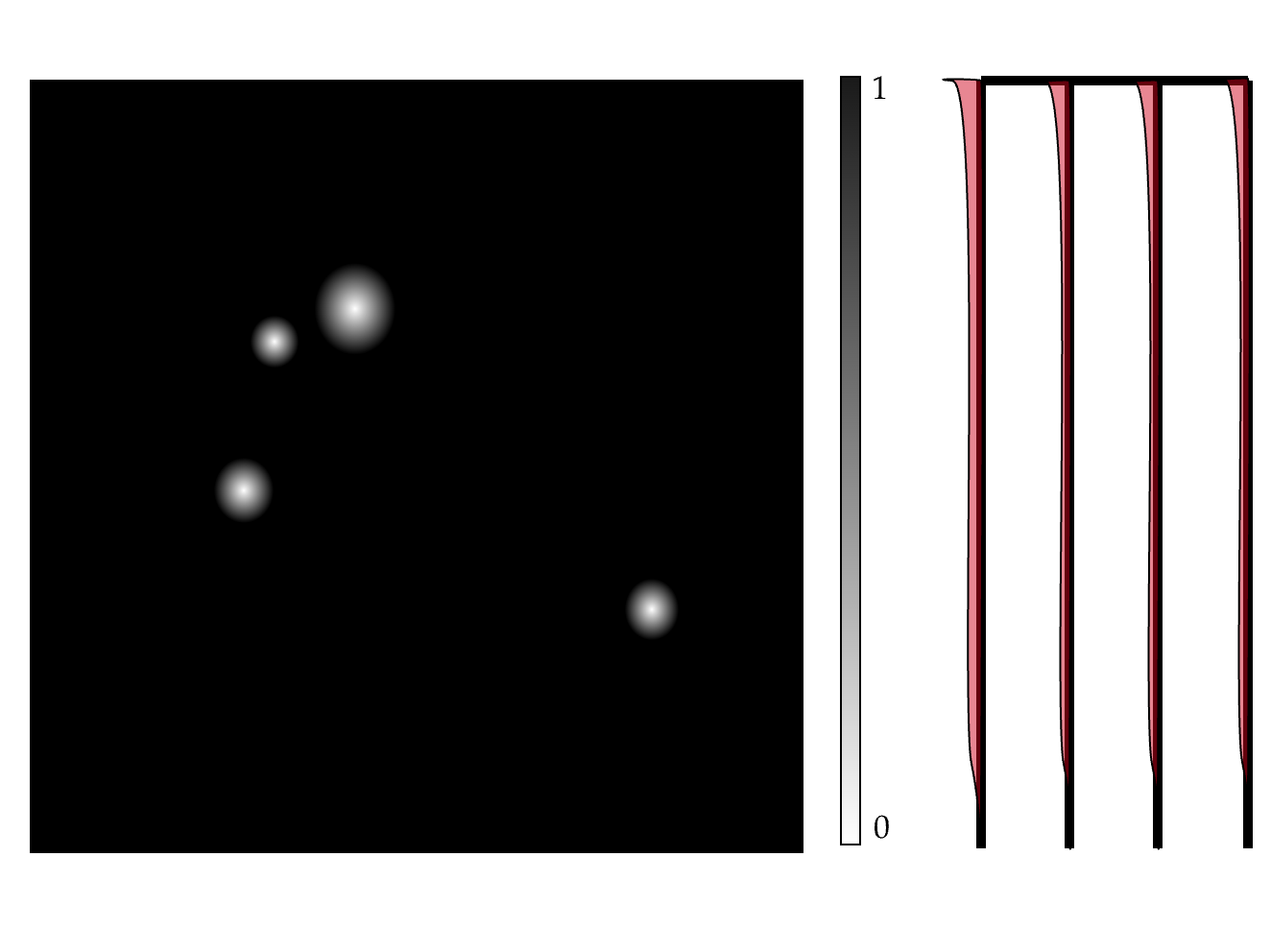}
            \captionsetup{singlelinecheck=off, margin={0.3cm, 0.1cm}}
		\caption{Same situation as in \Cref{fig:vol_w_0}, but for a different image. We clearly see that the four downward peaks have very different sizes, and that is reflected by the very different volume profiles attached to the branches of the merge tree. Note that the two images have the same  merge trees.}
		\label{fig:vol_w_1}
	\end{subfigure}
\caption{Plots related to \Cref{sec:dec_MT}.}
%\label{fig:cw_segments}
\end{figure}

Let $Y^k \subset \R^n$ be a finite set with $k$ points.
A widely employed filtration is:
\[
Y_t^k := \bigcup_{p\in Y^k} B_t(p)
\]
$B_t(p)$ being the metric ball centred in $p$, which is also called the offset filtration. 
Any point $D_{\pi_0(\Y)}$ is a pair $(t,U)$, with $t\in \R$ and $U$ being a path connected component of $Y_t^k$. Then, one
collect some topological information about $U$ and
 consider a vector-valued function $f$ with values in $X=\R^N_{\geq 0}$, for some $N<n$, defined by:
\begin{equation}\label{eq:betti}
f((t,U))=(\beta_0(U), \beta_1(U),\ldots,\beta_N(U)),    
\end{equation}
with $\beta_i(U)=\dim\text{H}_i(U,\R)$ being the $i$-th Betti number of $U$ - clearly $\beta_0(U)=1$. For an example see \cite{pegoraro2024finitelyfunc}.
An alternative and richer way to capture shape-related information about $U$ is to consider the point cloud: $U_k:=Y^k\cap U$ and then build $\text{PD}_i(U_k)$, the 
persistence diagram of the offset filtration of $U_k$ using $i$-dimensional homology.
Consequently, one can define $f$ with values in $X=D(\R^2_\Delta,\Delta)\times \ldots D(\R^2_\Delta,\Delta)$, the space of persistence diagrams with the $1$-Wasserstein distance, defined by:
\[
f((t,U))=(\text{PD}_0(U_k),\ldots,\text{PD}_N(U_k)).
\]
See \Cref{fig:PD_w_0} and 
\Cref{fig:PD_w_1} for a visual example. Observe how, at each value $t$, the merge tree gives a decomposition of the point cloud and the local shape of its clusters is captured by persistence diagrams.

Consider now a Lebesgue measurable function $g:W\rightarrow \R$, with $W$ open in $\R^n$. 
Build $Y_t=f^{-1}((-\infty,t])$, the sublevel set filtration of $f$. And consider
$D_{\pi_0(\Y)}$. This time 
another shape-related information that we can extract from a pair $(t,U)\in D_{\pi_0(\Y)}$,  is the Lebesgue measure of $U$. And thus define:  
\[
f((t,U))=\mathcal{L}(U),
\]
with $\mathcal{L}(U)$ being the Lebesgue measure of $U$. For an example of situations that can be tackled with this approach see \Cref{fig:vol_w_0} and
\Cref{fig:vol_w_1}. Note that the representation of functional data obtained with the volume-enriched persistent sets, lies in between information captured by $L_p$ metrics and purely topological information: with persistent sets one would not be able to distinguish the two functions in \Cref{fig:vol_w_0} and
\Cref{fig:vol_w_1} (they have isomorphic persistent sets), but with $L_p$ metrics one would not capture any of the shape similarities between the functions. 
We therefore believe that volume-enriched persistent sets can be a useful tool to employ when there are alignment issues with functional data i.e. all those situations where one may want to resort to some reparametrization of the domain - see, for instance, the Special Section on Time Warpings and Phase Variation
on the Electronic Journal of Statistics, Vol 8 (2). One such example can be found in \cite[Appendix,]{pegoraro2024finitelyfunc}.

Some of these ideas are formalized and in-depth explored in \cite{pegoraro2024finitelyfunc}, which focuses on representing via dendrograms objects of the form $f:D_{\pi_0 \circ \Y}\rightarrow X$, with $X$ being an editable space
and
$\pi_0 \circ Y:\R\rightarrow \Sets$ being a constructible persistent set.

\subsection{Problems to be Faced}
\label{sec:problems}

In the previous section, we outlined ideas suggesting that different persistent sets can be represented by dendrograms equipped with weight functions taking values in various editable spaces. As noted at the beginning of \Cref{sec:examples}, these ideas require formal validation before they can be reliably applied. In this section, we briefly discuss the main challenges associated with this proposed framework. Some of these challenges have been addressed in separate work, while others remain open and point to promising directions for future research.

There are three main issues that needs to be taken care of:
\begin{enumerate}
    \item Discretization: \Cref{sec:finite_pos} presents a procedure for turning a persistent set into a tree structure. This process depends crucially on the set of critical values $C$, and must be paired with the definition of an appropriate weight function, in order to obtain a dendrogram. Each of these steps—the discretization and the weight assignment—must be analyzed both individually and in combination. For example, one must verify whether the mapping from a persistent set to a dendrogram defines an injective transformation. This is not the case for the merge tree representation in \cite{merge_farlocca_2}, while injectivity is ensured in \cite{pegoraro2024finitely, pegoraro2024finitelyfunc}.
    \item Interpretability: the pipeline of representing a datum via an enriched persistent set and then with a dendrogram with values in an editable space, is in itself complex and needs to be carefully designed. Beyond the construction, it is also crucial to evaluate how well the edit distance $d_E$ captures variability between dendrograms.  For instance, if the editable space is $\R_{\geq 0}$, it is clear that editing a positively weighted tree with our edits, always gives a positively weighted tree. However, when working with functions defined on display posets, \cite{pegoraro2024finitelyfunc} shows that only a subset of dendrograms corresponds to actual functions. Therefore, additional care is needed to ensure that $d_E$ meaningfully captures the variability between such functions.
    \item Stability: closely related to interpretability is the issue of establishing stability guarantees specific to the chosen pipeline, akin to the results found in \cite{pegoraro2024functional, pegoraro2024finitelyfunc}. These results are essential to ensure that small perturbations in the data lead to correspondingly small changes in the associated dendrograms under $d_E$.
\end{enumerate}

\subsection{Stability Properties}

In \Cref{sec:problems}, we mention that establishing stability results is of pivotal importance to validate theoretically a pipeline involving dendrograms. Since some stability results have already been obtained in \cite{pegoraro2024finitely, pegoraro2024finitelyfunc}, we briefly describe them in order to give a qualitative idea of how the metrics $d_E$ behaves in different contexts.

In \cite{pegoraro2024finitely, pegoraro2024functional} $d_E$ is extended (via some non-trivial steps) to merge trees. Building on such definition, as already mentioned in the introduction, it is therein shown that:
\begin{equation}\label{eq:stability}
d_I(T_f,T_g)\leq d_E(T_f,T_g) \leq 2\cdot \text{size}(T_f) \cdot \text{size}(T_g) \parallel f-g \parallel_\infty
\end{equation}
with:
\begin{itemize}
\item $f,g$ being two tame \citep{chazal2016structure} real valued functions defined on the same path-connected topological space;
\item $T_f$ and $T_g$ are the merge trees associated to the sublevel set filtrations of $f$ and $g$;
\item $d_I$ is the interleaving distance between merge trees \citep{merge_interl};
\item $d_E$ is an adaptation of metric we define in this work, so that it induces a metric between merge trees;
\item $\text{size}(T_f)$ being the cardinality of the edge set of $T_f$. 
\end{itemize}

\Cref{eq:stability} represents a stability property in line with our expectations: 
 since our metric is a summation of all the costs of the local modifications needed to align two trees we have that the to control the distance between trees with the sup norm between functions, we need a constant that scales linearly with the size of the trees. 
 The interleaving distance, instead, is  focused on determining the maximal local modifications to turn one tree into the other and it is well known to be \emph{universally} stable, in the sense that is the biggest distance which is always bounded by the sup norm between the functions. Moreover, we have already pointed out that the $1$-Wasserstein distance between persistence diagrams satisfies an analogous relation  w.r.t. the bottleneck distances:

\[
d_B(D_f,D_g)\leq W_1(D_f,D_g)\leq (\#D_f+\#D_g)\parallel f-g\parallel_\infty.
\]
As a consequence, our metric is better at picking up differences (analogously to Wasserstein distances w.r.t bottleneck distance for persistence diagrams). On the other hand, this also means that is it more sensitive to noise.
Building on these considerations, in \cite{pegoraro2024finitely} it is argued that such stability properties may be preferable to universal stability properties, in which the operator from functions to topological representation is $1$-Lipschitz. Moreover, such comparison is framed in analogy to the bias-variance tradeoff in statistical modeling. 
We further observe that, since merge trees can be viewed as a special case of Reeb graphs, the same considerations we made for $d_I$ also apply to all the metrics discussed in \cite{bollen2022stable}—including those introduced in \cite{bauer2014measuring, di2016edit, bauer2020reeb}—due to their metric equivalence with the interleaving distance.

Lastly, \cite{pegoraro2024finitelyfunc} considers dendrograms enriched with Betti numbers of (path-connected components of) sublevel sets of functions, similarly to \Cref{eq:betti}. It is therein shown that:
\[
d_E(\varphi_f^{\Theta_p},\varphi_g^{\Theta_p}) \leq (\#\text{Crit}(f)+\#\text{Crit}(g)) 2C\parallel f-g\parallel_\infty,
\]

with:
\begin{itemize}
\item $f,g$ being two tame \citep{chazal2016structure} real valued functions defined on the same path-connected topological space;
\item $\varphi_f^{\Theta_p}$ and $\varphi_g^{\Theta_p}$ are the dendrograms associated to (enriched) the persistent sets of the sublevel set filtrations of $f$ and $g$;
\item $d_E$ is an adaptation of metric we define in this work, so that it induces a metric between functions defined on display posets of persistent sets;
\item $C>0$ is an upper bound to the values of the Betti numbers;
\item $\text{Crit}(f)$ is the set of values where the homology of the sublevel sets of $f$ changes, in any dimension. 
\end{itemize}

To conclude, we note that, if we consider the persistent homology of the sublevel set filtration of $f$ and $g$, in dimensions $i=0,\ldots,N$, and employ the persistence diagrams $D_f^i$ and $D_g^i$, we have:
\[
\sum_{i=0}^N W_1(D^i_f,D^i_g)\leq (\#\text{Crit}(f)+\#\text{Crit}(g)) \parallel f-g\parallel_\infty.
\]

\section{Computing the Edit Distance}
\label{sec:algo}

In these last section we develop an algorithm to compute the edit distance for dendrograms. 

\subsection{Decomposition Properties}
\label{sec:decomposition}

Following ideas found in \cite{TED} and exploiting the properties of editable spaces, we now prove some theoretical results which allow us to recursively split up the problem of finding an optimal mapping between two dendrograms. Then, in \Cref{sec:MILP}, we formulate each of the split-up calculations as binary optimization problems, which are then aggregated recursively in a bottom-up fashion   
in \Cref{sec:bottom_up}·

Name $T_2$ the only representative without degree $2$ vertices inside the equivalence class of $T$. One can always suppose that a dendrogram is given without degree $2$ vertices. 
Thus, for notational convenience, from now on we suppose $T=T_2$ and $T'=T'_2$.

We consider some particular subsets of $E_{T}\times E_{T'}$ which play a fundamental role in what follows.
Recall that, using $E_{T}\cong V_{T}- \{r_{T}\} $, we can induce $\pi_T:E_{T}\times E_{T'}\rightarrow V_T$.
%We define $\mathcal{C}^*(T,T')$ as follows. 

A set $M^*\subset E_{T}\times E_{T'}$ is in $\mathcal{C}^*(T,T')$ if:
\begin{itemize}
\item[(A1)] the points in $\pi_T(M^*)$ form antichains in $V_T$ (and the same for $\pi_{T'}(M^*)$ in $V_{T'}$), with respect to the partial order given by $parent >child$. This means that any two distinct vertices of $T$ (respectively of $T'$) which appear in $M^*$ are incomparable with respect to \virgolette{$>$};
\item[(A2)] the projections $\pi_T:M^*\rightarrow V_{T}$ and $\pi_{T'}:M^*\rightarrow V_{T'}$ are injective.
\end{itemize}

\begin{figure}	
    \centering
    \begin{subfigure}[c]{\textwidth}    	
    	\includegraphics[width = \textwidth]{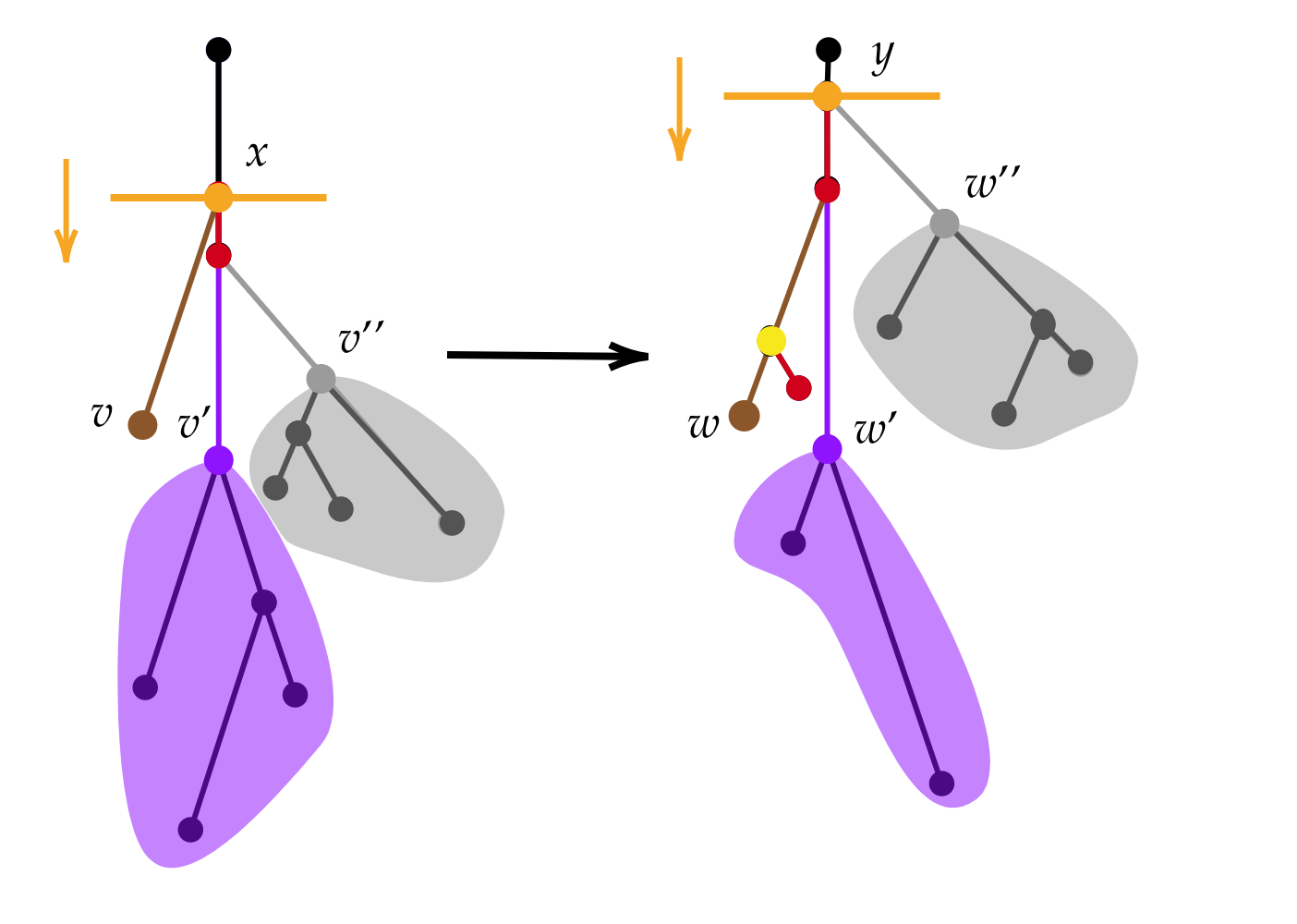}
	\end{subfigure}
%	\hspace{0.5 cm}
\caption{Given two weighted trees $T$ (left) and $T'$ (right) - which are the same of \Cref{fig:recap} - we consider $T_x=sub_T(x)$ and $T'_y=sub_{T'}(y)$ and we use \Cref{teo:decomposition} to compute 
$d_E(T_x,T'_y)$, as in the algorithm in \Cref{sec:MILP}. The set $M^*=\{(v,w),(v',w'),(v'',w'')\}$ satisfies (A1),(A2). The set $\alpha(M^*)$ is made by the deletions/insertions indicated by the red edges, the ghostings/splittings indicated by the yellow and the shrinkings given by edges of the same color, different from red. To obtain $\widetilde{T}_x$ and $\widetilde{T}'_y$ as in \Cref{prop:altro_mapping} all the black vertices covered by shaded regions must be deleted.}
\label{fig:decomposition}
\end{figure}

Consider now $M^*\in \mathcal{C}^*(T,T')$.
Starting from such set of pairs we build a set of edits which form a \virgolette{partial} mapping between $T$ and $T'$: each pair $(x,y)\in M^*$ means that we do not care of what lies below $x\in V_T$ and $y\in V_{T'}$ and we need to define edits only for the other vertices. The vertices below $x$ and $y$ will be taken care separately. See also \Cref{fig:decomposition} for a visual example.

Loosely speaking, $M^*$ is used as a \virgolette{dimensionality reduction tool}: instead of considering the problem of finding 
directly the optimal mapping between $T$ and $T'$, we split up the problem in smaller subproblems, and put the pieces together using $M^*$.
To formally do that, some other pieces of notation are needed.

Let $v\in E_T$. One can walk on the undirected graph of the tree-structure $T$ going towards any other vertex.
For any $v\in E_T$, $\zeta_v$ is the shortest
graph-path connecting $v$ to $r_{T}$. Note that this is the ordered set
$\zeta_v = \{v'\in V_T  \mid  v'>v \}$. Similarly, denote with $\zeta_x^{x'}$ the shortest path on the graph of $T$ connecting $x$ and $x'$.
Define the \emph{least common ancestor} of a set of vertices $A$ as $LCA(A):=\min \{v \in V_T \mid v\geq A\}=\min \bigcap_{a\in A} \zeta_a$.
In particular $\min \zeta_x\cap\zeta_{x'}$ is the least common ancestor between $x$ and $x'$: $LCA(x,x')=\min \zeta_x\cap\zeta_{x'}$.

By Property (A1), given $x\in V_T \cap \pi_T(M^*)$, 
%such that there is a pair $(x,y)$ in $M^*\in \mathcal{C}^*(T,T')$, 
there exist a unique $\Omega_{M^*}(x)\notin \pi_{T}(M^*)$ such that:
\[
\Omega_{M^*}(x) = \min \{LCA(x,x')\mid x'\in \pi_{T}(M^*) \text{ and }x\neq x'\}
\]  
And the same holds for $y\in V_{T'} \cap \pi_{T'}(M^*)$.  
For ease of notation we will often avoid explicit reference to $M^*$ and write directly $\Omega(x)$. 

With these bits of notation, given $M^*\in \mathcal{C}^*(T,T')$, we build the \virgolette{partial} mapping $\alpha(M^*)$: is a mapping that ignores all the vertices which lie below $x\in V_T$ and $y\in V_{T}$ if $(x,y)\in M^*$. \Cref{fig:decomposition} may help in following the upcoming paragraph. Consider $v\in V_T$:

\begin{enumerate}
\item if $(v,w)\in M^*$, then $(v,w)\in\alpha(M^*)$;
\item if there is not $x\in V_T$ such that $v<\Omega(x)$ or 
$v>\Omega(x)$, then $(v,\mathfrak{D})\in \alpha(M^*)$;
\item if there is $x\in V_T$ such that $v>\Omega(x)$ then $(v,\mathfrak{D})\in \alpha(M^*)$;
\item if there is $x\in V_T$ such that $v<\Omega(x)$:
\begin{enumerate}
\item if $v \in \zeta_x^{\Omega(x)}$ then $(v,\mathfrak{G})\in \alpha(M^*)$
\item if $v<v_i$ for some $v_i\in \zeta_x^{\Omega(x)}=\{v_0<v_1<\ldots <v_n\}$ then $(v,\mathfrak{D})\in \alpha(M^*)$;
\item if $v<x$ no edit is associated to $v$.
\end{enumerate}
\end{enumerate}

\begin{rmk}
By Properties (A1) and (A2), the conditions used to build $\alpha(M^*)$ are mutually exclusive. This means that each $v\in V_T$ satisfies one and only one of the above conditions and so $\alpha(M^*)$ is well defined.
\end{rmk}

The idea behind $\alpha(M^*)$ is that, for all pairs $(x,y)\in M^*$, we want to make the ghostings to turn the paths $\zeta_x^{\Omega(x)}$ and $\zeta_y^{\Omega(y)}$ respectively into the edges $(x,\Omega(x))$ and $(y,\Omega(y))$, and then shrink one in the other. 
As we already anticipated, $\alpha(M^*)$ takes care of all the vertices in $V_T$ and $V_{T'}$, a part from the sets $\cup_{(x,y)\in M^*}\{x'\in E_{T} \mid x'<x\}$ and $\cup_{(x,y)\in M^*}\{y'\in E_{T'} \mid y'<y\}$. For this reason we say that $\alpha(M^*)$ is a partial mapping.

We state this formally with the next proposition.

\begin{prop}\label{prop:altro_mapping}
Consider $T$ and $T'$ and $M^*\in \mathcal{C}^*(T,T')$. We obtain from such dendrograms, respectively, the dendrograms $\widetilde{T}$ and $\widetilde{T}'$ by deleting all the vertices $\bigcup_{(x,y)\in M^*}\{x'\in E_{T} \mid x'<x\}$ and $\bigcup_{(x,y)\in M^*}\{y'\in E_{T'} \mid y'<y\}$.
The set $\alpha(M^*)$ is a mapping in $M_2(\widetilde{T},\widetilde{T}')$.
\end{prop}

Now we have all the pieces we need to obtain the following key result.

\begin{teo}[Decomposition]
\label{teo:decomposition}
Given $T$, $T'$ dendrograms:

\begin{equation}
\label{eq:deco_eqn}
d_E(T,T')= \min_{M^*\in \mathcal{C}^*(T,T')} \sum_{(x,y)\in M^*}d_E(sub_T(x),sub_{T'}(y))+cost(\alpha(M^*)).
\end{equation}
\end{teo}

\subsection{Dynamical Binary Linear Programming problems}
\label{sec:MILP}

We want to use \Cref{teo:decomposition} to write a dynamical, binary linear optimization algorithm to calculate $d_E$: by translating  \Cref{teo:decomposition} into a  Binary Linear Programming (BLP) problem, we obtain a single step in a bottom-up procedure.

\subsubsection{Notation}
We are given two dendrograms $T,T'$ inside some dendrogram space $(\mathcal{T},X)$. Our objective 
is to write down \Cref{eq:deco_eqn} as a function of some binary variables.

Consider $x\in V_T$ and $y\in V_{T'}$.
Along with keeping the notation defined in  \Cref{sec:decomposition},  define  $T_x := sub_T(x)$ and $T_y := sub_{T'}(y)$, $N_x := \dim(T_x) = \#E_T $ and $N_y := \dim(T_y) = \#E_{T'}$. Moreover, given $v\in V_{T_x}$, the sequence $v_0=v<v_1<\ldots < r_T$ indicates the points in the path $\zeta_v$. Thus $v_i$ will be a vertex $v_i>v$. The same with $w\in V_{T_y}$ .

\subsubsection{Relaxing the Optimization Problem}
We would like to find $M^*\in \mathcal{C}^*(T_x,T_y)$ minimizing Equation \eqref{eq:deco_eqn} for $T_x$ and $T_y$, but this is a difficult task.
In fact, as evident in the construction of $\alpha(M^*)$, a set $M^*\in \mathcal{C}^*(T_x,T_y)$ has the role of pairing paths: if $(v,w)\in M^*$, then the paths 
$\zeta_v^{\Omega(v)}$ and $\zeta_w^{\Omega(w)}$ are shrunk one on the other by $\alpha(M^*)$. However, the points $\Omega(v)$ and $\Omega(w)$ depend on the whole set $M^*$, and not simply on the pair $(v,w)$. Modeling such global dependence gives rise to non-linear relations between paired points, and so leading to a non linear cost function, in terms of points interactions, to be minimized. For this reason we \virgolette{weaken} the last term in Equation \eqref{eq:deco_eqn},
 allowing also mappings different from $\alpha(M^*)$ to be built from $M^*$. 
In other words we minimize over $M^*\in \mathcal{C}^*(T_x,T_y)$ the following equation:
\begin{equation}
\label{eq:beta}
 \sum_{(v,w)\in M^*}d_E(sub_{T_x}(v),sub_{T_y}(w))+cost(\beta(M^*))
\end{equation}
where $\beta(M^*)$ is such that:
\begin{itemize}
\item $\beta(M^*)\in M_{\Gcal}(\widetilde{T}_{x},\widetilde{T}_{y})$ (using the notation established in  \Cref{prop:altro_mapping}, replacing $T$ and $T'$ with $T_x$ and $T_y$ respectively);
\item the set of vertices paired by $\beta(M^*)$ is exactly $M^*$: $M^*=\beta(M^*)\bigcap V_{T_x}\times V_{T_y}$.
\end{itemize}

Since, by construction $M^*=\alpha(M^*)\bigcap V_T\times V_{T'}$ and by \Cref{prop:altro_mapping}, $\alpha(M^*)\in M_2(\widetilde{T}_x,\widetilde{T}_{y})$, minimizing Equation \eqref{eq:deco_eqn} or Equation
\eqref{eq:beta} gives the same results.

\subsubsection{Setup and Variables}
%\label{sec:variables}
Suppose we already have $W_{xy}$ which is a $N_x\times N_y$ matrix such that $(W_{xy})_{v,w}= d_E(T_v,T_w)$ for all $v\in E_{T_x}$ and $w\in E_{T_y}$.
Note that:
\begin{itemize}
\item if $x$ and $y$ are leaves, $W_{xy}=0$.
\item if $v,w$ are vertices of $T_x$, $T_y$, then $W_{vw}$ is a submatrix of $W_{xy}$.
\end{itemize} 

The function to be optimized is defined on the following set of binary variables:
for every $v\in E_{T_x}$ and $w\in E_{T_y}$, for $v_i\in \zeta_v$, $v_i<r_{T_x}$, and $w_j\in \zeta_w$, $w_j<r_{T_y}$, take a binary variable $\delta^{v,w}_{i,j}$. We use $\delta$ to indicate the matrix of variables $(\delta^{v,w}_{i,j})_{v,w,i,j}$.

The mapping $\beta(M^*)$ is built according to the variables $\delta^{v,w}_{i,j}$ with value $1$:
we write a constrained optimization problem such that having $\delta^{v,w}_{i,j}=1$ means pairing the paths $\zeta_v^{v_{i+1}}$ (that is, the path which starts with $(v,v_1)$ and ends with $(v_{i},v_{i+1})$) and $\zeta_w^{w_{j+1}}$, and shrinking one in the other in the induced mapping. 

In order to pair and shrink the paths $\zeta_v^{v_{i+1}}=\{v=v_0,v_1,\ldots,v_{i+1}\}$ and $\zeta_w^{w_{j+1}}$ we need to collect some edits in the set $\beta(M^*)$ adding the following edits: 
\begin{itemize}
\item all the points $v_k\in \zeta_v^{v_{i+1}}$ with $0<k<i+1$ are ghosted, that is $(v_k,G)\in \beta(M^*)$;
\item if $v'<v_k$ for some $0<k<i+1$, then $(v',D)\in \beta(M^*)$;
\item if $v'\geq v_{i+1}$ and $v' \neq r_T$, then $(v',D)\in \beta(M^*)$ 
\item $(v,w)\in \beta(M^*)$. 
\end{itemize}

In the end, all edges which are not assigned an edit operation, are deleted. Note that analogous edits must be induced on vertices in $T_y$.
Thus, the edit $(v,w)\in \beta(M^*)$, along the edit paths induced by $\beta(M^*)$, means: shrinking the edge $(v,v_{i+1})$ onto $(w,w_{j+1})$). 
Recall that, if $\delta^{v,w}_{i,j}=1$, we do not need to define 
edits for $sub_{T_x}(v)$ and $sub_{T_y}(w)$ since, by assumption, we already know $d_E(T_v,T_w)$.

\subsubsection{Constraints}
Clearly, the set of $\delta^{v,w}_{i,j}$ with value $1$ does not always produce
$M^*\in \mathcal{C}^*(T_x,T_y)$: for instance paths could be paired multiple times.
To avoid such issues, we build a set of constraints for the variable $\delta$.

For each $v'\in V_{T_x}$ we call $\Phi(v'):=\{(v'',i)\in  V_T \times \mathbb{N}   \mid  v'=v''_i\in \zeta_{v''}^{v''_{i+1}}\}$. In an analogous way we define $\Phi(w')$ for $w'\in V_{T_y}$. 
Call $\mathcal{K}$ the set of values of $\delta$ such that for each leaf $l$ in $V_{T_x}$:
\begin{equation}
\label{eq:constr_1}
\sum_{v'\in \zeta_l}\left( \sum_{(v'',i)\in \Phi(v')}\sum_{w,j}  \delta^{v'',w}_{i,j} \right) \leq 1
\end{equation}
and for each leaf $l'$ in $V_{T_y}$:
\begin{equation}
\label{eq:constr_2}
\sum_{w'\in \zeta_{l'}}\left(\sum_{(w'',j)\in\Phi(w')}\sum_{v,i} \delta^{v,w''}_{i,j} \right)\leq 1
\end{equation}

The following proposition clarifies the properties of any value of $\delta\in\mathcal{K}$.

\begin{prop}\label{prop:constraints}
If $\delta\in\mathcal{K}$:
\begin{itemize}
\item the pairs $(v,w)$ such that  $\delta^{v,w}_{i,j}=1$ define a set $M^*\in \mathcal{C}^*(T_x,T_y)$;
\item the edits induced by all $\delta^{v,w}_{i,j}=1$ give a mapping $\beta(M^*)$ in $M_{\virgolette{G}}(\widetilde{T}_{x},\widetilde{T}_{y})$. With $\widetilde{T}_{x},\widetilde{T}_{y}$ being obtained from $T_x$ and $T_y$ as in \Cref{prop:altro_mapping}. 
\end{itemize}
\end{prop}

\begin{rmk}
If for every $\delta^{v,w}_{i,j}=1$ we have $v_{i+1}=\Omega(v)$, then $\beta(M^*)=\alpha(M^*)$.
\end{rmk}

\subsubsection{Objective Function}
Having built a mapping $\beta(M^*)$ using the binary variables, we want to define a cost functions which computes the cost of such mapping depending on 
$\delta$.

As before, consider $v\in E_{T_x}$ and interpret $\delta^{v,w}_{i,j}=1$ as pairing the paths $\zeta_v^{v_{i+1}}$ and $\zeta_w^{w_{j+1}}$;  then $v$ is paired with some $w\in E_{T_y}$ if $C(v):=\sum_{i,w,j}\delta^{v,w}_{i,j}=1$ and is ghosted if $G(v):=\sum_{\{i,v' \mid v\in \zeta_{v'}^{v'_{i}}\}}\sum_{w,j} \delta^{v',w}_{i,j}=1$.
The vertex $v$ is instead deleted if $D(v):=1-C(v)-G(v)=1$.
%Note that $C(v)\cdot G(v)=0$. 
We introduce also the following quantities, which correspond to the cost of shrinking $\zeta_v^{v_{i+1}}$ on $\zeta_w^{w_{j+1}}$:
\[
\Delta^{v,w}_{i,j}=d_X\left(\bigodot_{v'\in \zeta_v^{v_i}} \varphi_{T_x}(v'),\bigodot_{w'\in \zeta_w^{w_j}} \varphi_{T_y}(w')\right)
\]
Note that the above operations are taken inside the editable space $X$.

The function which computes the cost given by paired points is therefore: 

\[
F^C(\delta):=\sum_{v,w,i,j}\Delta^{v,w}_{i,j}\cdot \delta^{v,w}_{i,j}
\]

The contribution of deleted points is: $F^D(\delta)-F^-(\delta)$, where
\begin{equation*} 
F^D(\delta):=\sum_{v\in T_x} D(v)\cdot d(\varphi_{T_x}(v),0)+\sum_{w\in T_y} D(w)\cdot d(\varphi_{T_y}(w),0)
\end{equation*}

and

\begin{equation*}
F^-(\delta):=\sum_{v\in T_x} C(v)\cdot \mid \mid sub_{T_x}(v)\mid \mid +\sum_{w\in T_y} C(w)\cdot \mid \mid sub_{T_y}(w)\mid \mid 
\end{equation*}

where the \virgolette{norm} of a tree $T$ is  $\mid \mid T\mid \mid = \sum_{e\in E_T} d_X(\varphi(e),0)$.

Finally, one must take into account the values of $d_E(T_v,T_w)$, whenever $v$ and $w$ are paired; this information is contained in $(W_{xy})_{v,w}$:
\begin{equation*}
F^S(\delta):=\sum_{v,w}(W_{xy})_{v,w}\cdot 
 \left( \sum_{i,j}\delta^{v,w}_{i,j} \right) 
\end{equation*}

\begin{prop} 
\label{prop:opt_prob}
With the notation previously introduced:
\begin{equation}
\label{eq:problem}
d_E(T_x,T_y)=\min_{\delta\in \mathcal{K}}F^C(\delta)+F^D(\delta)-F^-(\delta)+F^S(\delta)
\end{equation}
\begin{proof}
The contribution of paired points is $F^C(\delta)$
and the contribution of deleted points is $F^D(\delta)-F^-(\delta)$. 

The cost of $\beta(M^*)$ is:
$F^\beta(\delta):= F^C(\delta)+F^D(\delta)-F^-(\delta)$.
Lastly, $F^S(\delta)$ takes into account the value of $d_E(T_v,T_w)$, if $v$ and $w$ are paired.
By  \Cref{teo:decomposition}, combined with  \Cref{prop:constraints}, the solution of the following optimization problem:
\begin{equation}
\min_{\delta\in\mathcal{K}}F^S(\delta)+F^\beta(\delta)
\end{equation} 
 is equal to $d_E(T_x,T_y)$.
\end{proof}
\end{prop}

\begin{rmk}
A solution to Problem \Cref{eq:problem} exists because the minimization domain is finite and there are admissible values; it is not unique in general.
\end{rmk}

\subsection{Bottom-Up Algorithm}
\label{sec:bottom_up}

In this section the results obtained in  \Cref{sec:decomposition} and the formulation established in \Cref{sec:MILP} are used to obtain the algorithm implemented to compute the metric $d_E$ between  dendrograms.
Some last pieces of notation are introduced in order to describe the \virgolette{bottom-up} nature of the algorithm.

Given $x\in V_T $, define $\len(x)$ to be the number of vertices in $\zeta_x$ and $\len(T)= \max_{v\in V_T}\len(v)$. Then, we set $\lvl(x)= \len(T)-\len(x)$.

The key property is that: $\lvl(x)> \lvl(v)$ for any $v\in sub(x)$.
Thus, if $W_{xy}$ is known for any $x\in \lvl_T^{-1}(n)$ and $y\in \lvl_{T'}^{-1}(m)$, then  for any $v$, $w$ in $V_T$, $V_{T'}$ such that $\lvl(v)<n$ and $\lvl(w)<m$, $W_{vw}$ is known as well.
With this notation we can write down Algorithm \ref{alg:bottomup}.

\begin{algorithm}
\SetAlgoLined
\KwResult{$d_E(T,T')$ }
 initialization: $N=\len(T)$, $M = \len(T')$, $n=m=0$\;
 \While{$n\leq N$ or $m\leq M$}{
	 \For{$(x,y)\in V_T\times V_{T'}$ such that $\lvl(x)\leq n$ and $\lvl(y)\leq m$}{
	 	 Calculate $(W_{r_{T} r_{T'}})_{x,y}$ solving Problem \eqref{eq:problem}\;
	 } 
	$n=n+1$;
	$m=m+1$;
   
 }
\Return $(W_{r_{T} r_{T'}})_{r_T,r_{T'}}$
 \caption{Bottom-Up Algorithm.}
 \label{alg:bottomup}
\end{algorithm}

An hands-on example of some iterations of Algorithm \ref{alg:bottomup} can be found in \Cref{sec:algo_ex}.
%and some simulations are carried out in \Cref{sec:simulations}.

\subsection{Computational Cost}

We end up with a result to analyze the performances of Algorithm \ref{alg:bottomup} in the case of dendrograms with binary tree structures. 
  
\begin{prop}[Computational Cost]\label{prop:complexity}
Let $T$ and  $T'$ be two dendrograms with full binary tree structures with $\dim(T)=\#E_T = N$ and $\dim(T')=M$.

Then $d_E(T,T')$ can be computed by solving $O(N\cdot M)$ BLP problems with $O(N\cdot \log(N)\cdot M\cdot \log(M))$ variables and $O(\log_2(M)+ \log_2(N))$ constraints.
\begin{proof}
In a full binary tree structure, at each level $l$ we have $2^l$ vertices. Let $L=\len(T)$ and $L'=\len(T')$. We have that, for any vertex $v\in  \lvl(l)^{-1}$, the cardinality of the path from $v$ to any of the leaves in $sub_T(v)$ is ${L-l}$ and the number of leaves in $sub_T(v)$ is $2^{L-l}$.

So, given $v\in V_T$ at level $l$ and $w\in V_{T'}$ at level $l'$, to calculate $d_E(sub_T(v),sub_{T'}(w))$ (having already $W_{vw}$) we need to solve a binary linear problem with $2^{L-l}\cdot (L-l)\cdot 2^{L'-l'}\cdot (L'-l')$ variables and $2^{L-l} + 2^{L'-l'}$ linear constraints.

Thus, to calculate $d_E(T,T')$, we need to solve $(2^{L+1}-1)\cdot (2^{L'+1}-1)$ linear binary optimization problems, each with equal or less than $2^L\cdot L \cdot 2^{L'}\cdot L'$ variables and equal or less that $2^L + 2^{L'}$ constraints. Substituting $L = \log_2(N)$ and $L' = \log_2(M)$ in these equations gives the result.
\end{proof}
\end{prop}

Note that binary dendrograms are dense (with respect to $d_E$) in any dendrogram space as long as for any $\varepsilon >0$, there is $x\in (X,\odot,0_X)$ such that $d_X(x,0)<\varepsilon$.

\begin{rmk}\label{rmk:comp_cost}
We end the manuscript with a remark about the computational cost of $d_E$ as described in 
\Cref{prop:complexity}. 
We find this result surprising to some extent. If we consider a recent implementation of the classical edit distance between unlabeled and unordered trees obtained with BLP \citep{TED}, 
then the complexity of the two metrics is very similar: the classical edit distance can be computed by solving $O(N\cdot M)$ BLP problems with $O(N\cdot  M)$ variables and $O(\log_2(M)+ \log_2(N))$ constraints - $O(N+M)$ if we count also the constraints restricting the integer variables to $\{0,1\}$, as the authors of \cite{TED} do. Thus, the introduction of the ghosting and splitting edits, with all the complexities they carry, just increases the cost by two $\log$ factors in the number of variables. 
Moreover, the condition used in \cite{zhang1996constrained} to make the classical edit distance tractable (that is, item (2) in Section 3.1. of \cite{zhang1996constrained}) is compatible with (M1)-(M4) and thus a constrained version of $d_E$ is well defined. We believe that it could be computed with a polynomial time algorithm. We leave this investigation to future works.
\end{rmk}

\section{Conclusions}
\label{sec:conclusions}

We develop a novel framework to work with 
tree that arise as topological summaries in the context of topological data analysis, possibly enriched with some abstract weights as a result of different pipelines. These kind of summaries are increasingly frequent and require some particular care to be properly handled by a metric. 
We define an edit distance employing some novel edit operations, proving that such definition induces a metric up to a certain equivalence relation between weighted trees.

By exploiting the ordering properties of trees, we produce a binary linear programming algorithm to compute such metric. The modifications that we make w.r.t. the classical edit distance between unlabeled trees make it only marginally more expensive, namely by some \emph{log} factors. For this reason we argue that polynomial time approximations should be available also for $d_E$. Still, the computational procedure we describe here is used in separate applied works, as mentioned in the introduction.
These additional modifications are essential for the  stability properties assessed in \cite{pegoraro2024functional, pegoraro2024finitelyfunc}, justifying the higher computational cost w.r.t. other edit distances for merge trees, which, instead, are unstable.

The generality of the work opens up many possible research directions. Some of which are already studied in other works, while others still await for investigations: 
\begin{itemize}
\item we think that the properties of the editable spaces can be relaxed; however, the algorithm presented in this manuscript may need to be adapted to the properties of the chosen weight space;  
\item we would like to extend this edit distance outside tree-shaped graphs, encompassing Reeb graphs \citep{reeb_2, de2016categorified} and their combinatorial approximation, namely \emph{mapper graphs} \citep{singh2007topological} are studied \citep{carriere2018statistical, carriere2018structure}, perhaps following the decomposition presented in \cite{stefanou2020tree}. We also mention that works on GEDs like 
\citep{ambauen2003graph, lerouge2016exact} could be of great help in developing such generalization;
\item as motivated by \Cref{rmk:comp_cost} we believe that polynomial time approximations can be obtained for computing $d_E$ between dendrograms;
\item further applications and case studies exploiting the generality of the weights we consider on the edges should be considered and explored.  
\end{itemize}

´

\section*{Acknowledgments}
This work was carried out as part of my PhD Thesis, under the supervision of Professor Piercesare Secchi.
I also acknowledge the support of the Wallenberg AI, Autonomous Systems and Software Program (WASP), and of the SciLifeLab and Wallenberg National Program for Data-Driven Life Science (DDLS), which fund the project: Topological Data Analysis of Functional Genome to Find Covariation Signatures.

\newpage

\appendix

\section*{Outline of the Appendix}

\Cref{sec:algo_ex} contains an hands-on example of some iterations of the algorithm presented in \Cref{sec:bottom_up}. %\Cref{sec:simulations} contains some numerical simulations. 
Lastly, 
\Cref{sec:proofs_TDA} contains the proofs of the results in the paper.




\section{Algorithmic Example}
\label{sec:algo_ex}

Here we present in details the first steps of the  Algorithm \ref{alg:bottomup}, used to calculate the distance between two merge trees.

We consider the following pair of weighted trees. Let $(T,w_T)$ be the tree given by: $V_T=\{a,b,c,d,r_T\}$, $E_T=\{(a,d),(b,d),(d,r_T),(c,r_T)\}$ and $w_T(a)=w_T(b)=w_T(d)=1$, $w_T(c)=5$; the tree $(T',w_{T'})$ instead, is defined by: $V_{T'}=\{a',b',c',d',r_{T'}\}$, $E_{T'}=\{(a',d'),(b',d'),(d',r_{T'}),(c',r_{T'})\}$ and 
$w_{T'}(a)=1$, $w_{T'}(b)=w_{T'}(c)=2$ and $w_{T'}(d)=3$. We employ the notation $w_T$ for the weight function instead of $\varphi_T$ 
to indicate the common situation of a tree with positive real weights.

\subsubsection*{Step: $n=m=0$}
This step is trivial since we only have pairs between leaves, like $(a,a')$, which have trivial subtrees and thus $d_E(sub_T(a),sub_{T'}(a'))=0$.

\subsubsection*{Step: $n=m=1$}
The points $x\in V_T$ with $\lvl_T(x)\leq 1$ are $\{a,b,c,d\}$ and 
the points $y\in V_{T'}$ with $\lvl_{T'}(y)\leq 1$ are $\{a',b',c',d'\}$. Thus the pairs $(x,y)$ which are considered are: $(d,d')$, $(d,a')$, $(d,b')$, $(d,c')$ and $(a,d')$, $(b,d')$, $(c,d')$. The pairs between leaves, like $(a,a')$ have already been considered.

\subparagraph*{Pair: $(d,d')$}
Let $T_d=sub_{T}(d)$ and $T_{d'}=sub_{T'}(d')$. The set of internal vertices are respectively $E_{T_{d}}=\{a,b\}$ and $E_{T_{d'}}=\{a',b'\}$. For each vertex $v<\text{root}$ in each subtree, where \virgolette{$\text{root}$} stands for $d$ or $d'$, roots of $T_{d}$ and $T_{d'}$ respectively, we have $\zeta_v=\{v_0=v,v_1=\text{root}\}$.  
Thus, the binary variables we need to consider, are the following: 
$\delta_{0,0}^{a,a'}$,
$\delta_{0,0}^{a,b'}$, $\delta_{0,0}^{b,a'}$ and $\delta_{0,0}^{b,b'}$.
The quantities $\Delta_{i,j}^{v,w}$ are given by:
$\Delta_{0,0}^{a,a'}=0$,
$\Delta_{0,0}^{a,b'}=1$, $\Delta_{0,0}^{b,a'}=0$ and $\Delta_{0,0}^{b,b'}=1$.
Thus:
\[
F^C(\delta)= 0\cdot \delta_{0,0}^{a,a'}+\delta_{0,0}^{a,b'}+0\cdot\delta_{0,0}^{b,a'}+\delta_{0,0}^{b,b'}
\]
While: 
\[
F^D(\delta)= (1-\delta_{0,0}^{a,a'}-\delta_{0,0}^{a,b'})\cdot 1 
+ (1-\delta_{0,0}^{b,a'}-\delta_{0,0}^{b,b'})\cdot 1 +
(1-\delta_{0,0}^{a,a'}-\delta_{0,0}^{b,a'})\cdot 1 +
(1-\delta_{0,0}^{a,b'}-\delta_{0,0}^{b,b'})\cdot 2
\]
and:
\[
F^-(\delta)= (\delta_{0,0}^{a,a'}+\delta_{0,0}^{a,b'})\cdot 0 
+ (\delta_{0,0}^{b,a'}+\delta_{0,0}^{b,b'})\cdot 0 +
(\delta_{0,0}^{a,a'}+\delta_{0,0}^{b,a'})\cdot 0 +
(\delta_{0,0}^{a,b'}+\delta_{0,0}^{b,b'})\cdot 0
\] 
and:
\[
F^S(\delta)= \delta_{0,0}^{a,a'}\cdot 0 +\delta_{0,0}^{a,b'}\cdot 0+\delta_{0,0}^{b,a'}\cdot 0+\delta_{0,0}^{b,b'}\cdot 0
\]

Lastly the constraints are:
\[
\delta_{0,0}^{a,a'}+\delta_{0,0}^{a,b'}\leq 1;\text{ }
\delta_{0,0}^{b,a'}+\delta_{0,0}^{b,b'}\leq 1;\text{ }
\delta_{0,0}^{a,a'}+\delta_{0,0}^{b,a'}\leq 1;\text{ }
\delta_{0,0}^{a,b'}+\delta_{0,0}^{b,b'}\leq 1
\]
A solution is given by $\delta_{0,0}^{a,a'}=\delta_{0,0}^{b,b'}=1$ and $\delta_{0,0}^{a,b'}=\delta_{0,0}^{b,a'}=0$, which entails $F^C(\delta)=1$, $F^D(\delta)=0$, $F^-(\delta)=0$ and $F^S(\delta)=0$ 
and $d_E(T_{d},T_{d'})=1$.

\subparagraph*{Pair: $(d,a')$}
Obviously: $d_E(sub_T(d),sub_{T'}(a'))=\mid \mid sub_T(d)\mid \mid $.
All the pairs featuring a leaf and an internal vertex (that is, a vertex which is not a leaf), such as $(d,b')$, $(a,d')$ etc. behave similarly.

\subsubsection*{Step: $n=m=2$}
The points $x\in V_T$ with $\lvl_T(x)\leq 2$ are $\{a,b,c,d,r_T\}$ and 
the points $y\in V_{T'}$ with $\lvl_{T'}(y)\leq 2$ are $\{a',b',c',d',r_{T'}\}$. Thus the pairs $(x,y)$ which are considered are $(d,r_{T'})$, $(r_T,d')$, $(r_T,r_{T'})$ and then the trivial ones: 
$(r_T,a')$, $(r_T,b')$, $(r_T,c')$ and $(a,r_{T'})$, $(b,r_{T'})$, $(c,r_{T'})$. Some pairs have already been considered and thus are not repeated.

\subparagraph*{Pair: $(d,r_{T'})$}
Let $T_d=sub_{T}(d)$ and $T'=sub_{T'}(r_{T'})$. The set of internal vertices are respectively $E_{T_{d}}=\{a,b\}$ and $E_{T_{d'}}=\{a',b',c',d'\}$.  
Thus, the binary variables we need to consider, are the following: 
$\delta_{0,0}^{a,a'}$,
$\delta_{0,1}^{a,a'}$,
$\delta_{0,0}^{a,b'}$, 
$\delta_{0,1}^{a,b'}$, 
$\delta_{0,0}^{a,c'}$,
$\delta_{0,0}^{a,d'}$,
$\delta_{0,0}^{b,a'}$,
$\delta_{0,1}^{b,a'}$,
$\delta_{0,0}^{b,b'}$,
$\delta_{0,1}^{b,b'}$
$\delta_{0,0}^{b,c'}$,
and
$\delta_{0,0}^{b,d'}$.

The quantities $\Delta_{i,j}^{v,w}$ are given by:
$\Delta_{0,0}^{a,a'}=0$,
$\Delta_{0,1}^{a,a'}=3$,
$\Delta_{0,0}^{a,b'}=1$, 
$\Delta_{0,1}^{a,b'}=4$, 
$\Delta_{0,0}^{a,c'}=1$,
$\Delta_{0,0}^{a,d'}=2$,
$\Delta_{0,0}^{b,a'}=0$,
$\Delta_{0,1}^{b,a'}=3$,
$\Delta_{0,0}^{b,b'}=1$,
$\Delta_{0,1}^{b,b'}=4$,
$\Delta_{0,0}^{b,c'}=1$ and 
$\Delta_{0,0}^{b,d'}=2$.
The function $F^C(\delta)$ is easily obtained by summing over $\delta_{i,j}^{v,w}\cdot \Delta_{i,j}^{v,w}$.

While: 
\[
F^D(\delta)= (1-\delta_{0,0}^{a,a'}-\delta_{0,1}^{a,a'}-\delta_{0,0}^{a,b'}-\delta_{0,1}^{a,b'}-\delta_{0,0}^{a,c'}-\delta_{0,0}^{a,d'})\cdot 1 
+ \ldots +
(1-\delta_{0,0}^{a,d'}-\delta_{0,0}^{b,d'})\cdot 3
\]
and:
\[
F^-(\delta)= (\delta_{0,0}^{a,a'}+\delta_{0,1}^{a,a'}+\delta_{0,0}^{a,b'}+\delta_{0,1}^{a,b'}+\delta_{0,0}^{a,c'}+\delta_{0,0}^{a,d'})\cdot 0
+ \ldots +
(\delta_{0,0}^{a,d'}+\delta_{0,0}^{b,d'})\cdot 3
\] 
and:
\[
F^S(\delta)= (\delta_{0,0}^{a,a'}+\delta_{0,1}^{a,a'})\cdot 0 +
(\delta_{0,0}^{a,b'}+\delta_{0,1}^{a,b'})\cdot 0 
+ \ldots +
\delta_{0,0}^{a,d'}\cdot 3 +
\delta_{0,0}^{b,d'}\cdot 3 
\]

Lastly the constraints are:
\[
\delta_{0,0}^{a,a'}+\delta_{0,1}^{a,a'}+
\delta_{0,0}^{a,b'}+\delta_{0,1}^{a,b'}+
\delta_{0,0}^{a,c'}+\delta_{0,0}^{a,d'}\leq 1
\]
\[
\delta_{0,0}^{b,a'}+\delta_{0,1}^{b,a'}+
\delta_{0,0}^{b,b'}+\delta_{0,1}^{b,b'}+
\delta_{0,0}^{b,c'}+\delta_{0,0}^{b,d'}\leq 1
\]
\[
\delta_{0,0}^{a,a'}+\delta_{0,1}^{a,a'}+
\delta_{0,0}^{b,a'}+\delta_{0,1}^{b,a'}+
\delta_{0,0}^{a,d'}+\delta_{0,0}^{b,d'}\leq 1
\]
\[
\delta_{0,0}^{a,b'}+\delta_{0,1}^{a,b'}+
\delta_{0,0}^{b,b'}+\delta_{0,1}^{b,b'}+
\delta_{0,0}^{a,d'}+\delta_{0,0}^{b,d'}\leq 1
\]
\[
\delta_{0,0}^{a,c'}+\delta_{0,0}^{b,c'}\leq 1
\]
In this case there are many minimizing solutions. One is given by:
$\delta_{0,1}^{a,a'}=\delta_{0,0}^{b,c'}=1$ and all other variables equal to $0$. This value of $\delta$ is feasible since the variables $\delta_{0,1}^{a,a'}$ and $\delta_{0,0}^{b,c'}$ never appear in the same constraint.
This value of $\delta$ entails $F^C(\delta)=3+1$, $F^D(\delta)=2$, $F^-(\delta)=0$ and $F^S(\delta)=0$, 
and thus $d_E(T_{d},T')=6$.

Another solution can be obtained with:
$\delta_{0,0}^{a,d'}=\delta_{0,0}^{b,c'}=1$ and all other variables equal to $0$. Also this value of $\delta$ is feasible since the variables $\delta_{0,0}^{a,d'}$ and $\delta_{0,0}^{b,c'}$ never appear in the same constraint.
This value of $\delta$ entails $F^C(\delta)=2+1$, $F^D(\delta)=w_{T'}(a')+w_{T'}(b')=1+2$, $F^-(\delta)=\mid \mid sub_{T'}(d')\mid \mid =3$ and $F^S(\delta)=d_E(sub_T(a),sub_{T'}(d'))=\mid \mid sub_{T'}(d')\mid \mid =3$, 
and thus $d_E(T_{d},T')=3+3-3+3=6$.

\subparagraph*{Pair: $(r_T,d')$}
This and the other pairs are left to the reader.

\section{Proofs}\label{sec:proofs_TDA}

\medskip\noindent
\underline{\textit{Proof of}  \Cref{prop:wass_editable}.} 

We indicate with $\text{Lip}(\R^2_\Delta)$ the set of functions $k:\R^2_\Delta\rightarrow \R$ which are Lipschitz on 
$\R^2_\Delta$ endued with the $1$-norm and such that $k(\Delta)=0$.

By \cite{bubenik2022universality} we have:
\[
W_1([D],[D'])=W_1(\sum_{x\in [D]}\{x\}, \sum_{y\in [D']}\{y\})= \sup_{k\in\text{Lip}(\R^2_\Delta)}(\sum_{x\in [D]}k(x)- \sum_{y\in [D']}k(y)).
\]
At this point $(P3)$ and $(P4)$ follow easily by, respectively:
\begin{itemize}
\item  replacing $[D']$ with the empty diagram and noting:
\[
\sup_{k\in\text{Lip}(\R^2_\Delta)}(\sum_{x\in [D]}k(x))=\sum_{x\in [D]}\min_{y\in \Delta}\parallel x-y\parallel_1;
\]
\item observing that when we have $W_1([D+D''],[D'+D''])$ the contributions by the elements in $[D'']$ when computing:
\[
\sum_{x\in [D+D'']}k(x)- \sum_{y\in [D'+D'']}k(y)
\] 
appear in both sums, erasing themselves out.
\end{itemize}
\hfill$ \blacksquare $

\bigskip\noindent
\underline{\textit{Proof of}  \Cref{prop:graphs_trees}.}

\smallskip\noindent

For the construction of the merge tree associated to a persistent set see \cite{pegoraro2024finitelyfunc}.

Given a persistent set $S$ as in the hypotheses of the statement.
To simplify the notation we make a slight abuse of notation and call $S$ also the restriction of $S$ to its critical set $C$.

 Let $D_S$ be its display poset and $(T_S,f_S:V\rightarrow\R)$ with $T_S =(V,E)$ its merge tree. By construction 
$V-\{r_{T_S}\}$ - with $r_{T_S}$ being the root of the merge tree and $f_S(r_{T_S})=+\infty$ - can be sent into $D_S$: every vertex $v$ in $V-\{r_{T_S}\}$ is identified by a critical value $t_i=f_S(v)$ and an element $v\in S(t_i)$. Thus we have $V-\{r_{T_S}\}\hookrightarrow D_S$. Moreover, this embedding is an embedding of posets if we consider $T_S$ with the relation  generated by $child<parent $.

Consider now $(t_j,p)\in  D_S-V-\{r_{T_S}\}$:
this happens if there are $t_i<t_j$ and $t_k>t_j$ such that for every $\varepsilon_i>0$ such that $t_j-\varepsilon_i> t_i$ and for every $\varepsilon_k>0$ such that $t_j+\varepsilon_k < t_k$, we have
 $\#S(t_j\leq t_j+\varepsilon_k)^{-1}(S(t_j\leq t_j+\varepsilon_k)(p))=1$ and 
$\#S(t_j-\varepsilon_i\leq t_j)^{-1}(p)=1$. 
Moreover we can choose $t_j$ and $t_k$ so that:
\begin{itemize}
\item  $\max \{(t,s)\in D_S \mid (t,s) < (t_j,p)\}=(t_i,S(t_i\leq t_j)^{-1}(p))$;
\item  $\min \{(t,s)\in D_S \mid (t,s)> (t_j,p)\}=(t_k,S(t_i\leq t_k)(p))$
\end{itemize}
This implies that for any edge $e$ connecting $(t,s)\rightarrow (t_j,p)$ in the DAG representing $D_S$, we have that if $(t,s)\neq(t_i,S(t_i\leq t_j)^{-1}(p))$, then we also have an edge $(t,s)\rightarrow (t_i,S(t_i\leq t_j)^{-1}(p))$. But this implies that we can remove $e$ from the DAG without changing its transitive closure. That is, $e$ does not appear in $\Gcal(D_S)$. The same for any edge    
$(t_j,p)\rightarrow (t,s)$ with 
$(t,s)\neq(t_k,S(t_i\leq t_k)(p))$.
Thus $(t_j,p)$ induces a degree $2$ vertex in $\Gcal(D_S)$. 

To conclude the proof we just need to prove that $\Gcal(D_S)$ is a tree. In fact, if this is the case, then we have a poset map $V-\{r_T\}\hookrightarrow V_{\Gcal(D_S)}$ so that all the vertices not in the image of this map are of degree $2$. I.e. the map is an isomorphism of posets onto the image, which is the poset obtained by removing all degree $2$ vertices from $\Gcal(D_S)$. 

We know that $\Gcal(D_S)$ is a DAG which has, by construction one maximal element. The only thing we need to prove is that for every $x=(t_i,p)\in V_{\Gcal(D_S)}$ we only have one directed edge of the form $e=(x,y)$. But this is easily seen: $y=S(t_i\leq t_{i+1})(p)$. By functoriality, for every other $z=(t,s)\geq z$, we also have $z\geq y$. Thus we are done. 

\hfill$ \blacksquare $

\bigskip\noindent
\underline{\textit{Proof of}  \Cref{lemma:M_2}.}

\smallskip\noindent
Any degree $2$ vertex which is not ghosted is paired with another degree $2$ vertex. Ghosting both of them does not increase the cost of the mapping.

\hfill$ \blacksquare $

\medskip\noindent
\underline{\textit{Proof of}  \Cref{teo:main_thm}.} 

\smallskip\noindent
To lighten the notation we use the following symbols:
\begin{itemize}
\item the edit induced by $(v,\mathfrak{D})$ is called $v_d$ and $v^{-1}_d$ stands for $(\mathfrak{D},v)$.
\item the edit induced by $(v,\mathfrak{G})$ is called $v_g$ and $v^{-1}_g$ stands for $(\mathfrak{G},v)$.
\item the edit induced by $(v,v')$ is called $v_{\varphi,\varphi'}$ with $\varphi$ being the original weight function, and $\varphi'$ the weight function after the shrinking.
\end{itemize}

We know that the set of finite edit paths between two dendrograms is nonempty.

Suppose that $\gamma$ is a finite edit path.
This means that $\gamma$ is the composition of a finite set of edits. We indicate such ordered composition with $\gamma=\prod_{i=0}^N e_i$ with $e_i$ edit operation.
We would like to change the order of the edit operations without raising the cost and changing the endpoints of the edit path. This is not always possible. However we can work it around in the useful cases using properties (P1)-(P4). 
In particular, we would like to know  when we can commute a generic edit $e_i$ in the following situations:
\begin{itemize}
\item  $v_d \circ e_i$ and $e_i\circ v^{-1}_d$
\item $v_g\circ e_i$ and $e_i\circ v^{-1}_g$.
\end{itemize}

Moreover we want to reduce the edit path to max one edit for any vertex of $T$ and $T'$.

We divide the upcoming part of the proof in subsections, each devoted to different combinations of edits.

\subsubsection*{Inverse Operations}

We point out the following fact, which will be used often in the proof. Suppose we have a pair of edits $e_i \circ e_j$ that can be replaced with the pair $e'_i \circ e'_j$. Now consider $e^{-1}_j \circ e^{-1}_i$:
\begin{equation}\label{eq:inverse}
e^{-1}_j \circ e^{-1}_i = (e_i \circ e_j)^{-1}=(e'_i \circ e'_j)^{-1}=e'^{-1}_j \circ e'^{-1}_i.
\end{equation}

\subsubsection*{$v_d$ and $v_d^{-1}$}

When we delete or insert one vertex, we are modifying the tree 
structure at the level of its parent and its children. Therefore, we are only taking into account operations involving the parent of the considered vertex, on the vertex himself or on the children of the deleted/inserted vertex.

\begin{itemize}
\item $v_d \circ v'_g$: if $v$ is a child of $v'$, can be safely replaced with $v_d\circ v'_d$. Instead of ghosting the parent and then deleting the whole edge, we can delete both edges one by one; conserving the length of the path (P3). If, instead, $v$ is parent of $v'$, we can safely commute the operations.
\item $v_d\circ v_g^{-1}$: the two edits can be replaced with $v'_{\varphi,\varphi'}$ with $v'$ parent of $v$ (after the insertion) and $\varphi'$ properly defined not to raise the cost of the path. In fact we are inserting $v$ on an edge and then deleting it. This can obviously be achieved by shrinking the original edge (without changing the path length - (P4)).
\item $v_d\circ v'^{-1}_g$: if $v'$ is parent of $v$ after the splitting, we can again replace the two edits with shrinking: instead of inserting a point in an edge, and deleting then the edge below, we can directly shrink the original edge (P4). If, instead, $v'$ is inserted below $v$, we have the same situation, but seen from the point of view of the vertex which is the child of $v'$ after the splitting.
\item $v_d\circ v_{\varphi,\varphi'}$: the two edits can be replaced by $v_d$ potentially diminishing the length of the path, but surely not raising it (P1).
\item $v'_g\circ v_d^{-1}$: if $v'$ is the parent of $v$, this edit can be replaced with just $v'_{\varphi,\varphi'}$ with appropriate weights: we are inserting an edge under a vertex which (in this case) becomes of degree  two and is ghosted. We can directly modify the edge without changing the length of the path (P4). If, instead, $v'$ becomes a child of $v$ after the insertion, we can simply shrink $v'$ to obtain the same result without raising the cost (P4).
\item $v'^{-1}_g\circ v_d^{-1}$: if $v'$ is on the edge inserted with $v_d^{-1}$, we cannot commute these two edits,  but we can replace them with two insertions: instead of inserting an edge and then splitting it, we can directly insert two smaller edges; without changing the cost of the path (P3). In all other situations, we can commute the two operations.
\item $v_{\varphi,\varphi'}\circ v_d^{-1}$: the two edits can be replaced with an insertion directly with weight $\varphi'$, possibly shortening the path (P1).
\item $v'^{-1}_d\circ v_d$: suppose that the parent of $v'$, after its insertion, is the parent of $v$ before we apply both edits. If the children of $v'$ are different from the children of $v$, this operations do not commute. If the children are the same, they can be changed with a shrinking of $v$, reducing the length of the path by at most $cost(v'^{-1}_d)+cost(v_d)$ (P1). 
\end{itemize}

\subsubsection*{$v_g$ and $v_g^{-1}$}

Like in the previous case, we only take into account transformations concerning the parent and the child of the added/ghosted degree  two vertex.

\begin{itemize}
\item $v_g\circ v'_g$: always commute (P2).
\item $v_g\circ v'^{-1}_g$: always commute, provided that we define carefully the splitting $v'^{-1}_g$ (P2).
\item  $v_g\circ v_{\varphi,\varphi'}$: with such operations,  we are shrinking a vertex before ghosting it. However, we can achieve the same result, without increasing the path length, by ghosting the vertex at first, and then shrinking its child (P1)-(P4).
\item $v'_{\varphi,\varphi'}\circ v_g^{-1}$: if either $v'=v$, or $v$ is the parent of $v'$, we can replace those edits with an appropriate shrinking of the vertex which is the child of $v$ after the splitting and an appropriate insertion of $v'$, without changing the length of the path (P3)-(P4).
\item  $v_g \circ v'_d$: if $v$ is the parent of $v'$, such edits do not commute and can't be replaced by a similar operation which inverts ghosting and deletion. 
\end{itemize}

\subsubsection*{$v_{w,w'}$}
\begin{itemize}
\item $v_{\varphi',\varphi''}\circ v_{\varphi,\varphi'}$: those edit can be replaced by $v_{\varphi,\varphi''}$ which is either conserving or shortening the path (P1).
\item  $v_{\varphi,\varphi'}^{-1}$: we can proceed as in the previous point, since: $v_{\varphi,\varphi'}^{-1}=v_{\varphi',\varphi}$.
\end{itemize}

Thanks to these properties, we can take a given path $\gamma=\prod_{i=0,...,N} e_i$ and modify the edit operations in order to obtain the following situation:

\begin{itemize}
\item the first operations are all in the form $v_d$; this can be achieved because $ v_d\circ -$ can be always rearranged, potentially by changing the path (as shown before) and shortening it. Of course, there can be only one deletion for each vertex of $T$;
\item then we have all the edits in the form $v_g$; since $v_g\circ -$ is exchangeable any time but when we have $v_g \circ v'_d$, this is not a problem. Observe that all degree  two vertices which were not deleted can be ghosted (at most one time); 
\item in the same way we can put last all the paths in the form $v_d^{-1}$ and before them $v_g^{-1}$. All the new vertices appearing with the insertion of edges and the splitting of edges with degree  two vertices are all nodes which remain in $T'$ and which are not further edited;
\item in the middle we are left with the shrinking paths. Since we can substitute $v_{\varphi,\varphi'}\circ v_{\varphi',\varphi''}$ with $v_{\varphi,\varphi''}$, we can obtain at most one transformation on each vertex.
\end{itemize}

Thus, we can find at least one edit path of the form:
\[
\overline{\gamma}= (\gamma_d^{T'})^{-1} \circ (\gamma_g^{T'})^{-1} \gamma_s^T\circ\gamma_g^{T}\circ\gamma_d^T 
\]

with: 
\begin{itemize}
\item $\gamma_d^T=\prod v_d$;
\item $\gamma_g^T=\prod v_g$;
\item $\gamma_s^T=\prod v_{\varphi,\varphi'} $;
\item $(\gamma_g^{T'})^{-1}=\prod v_g^{-1} $;
\item $(\gamma_d^{T'})^{-1}=\prod v_d^{-1}$.
\end{itemize}

As for the paths $\gamma_M$ described in \Cref{sec:mappings}, all the permutations of the edits inside each $\gamma_d^T,\gamma_g^T,\ldots,$ does not change its endpoints and its cost.
Thus, $\overline{\gamma}$ is such that $\gamma(T)=\overline{\gamma}(T)=T'$ (up to isomorphism) and $cost(\overline{\gamma})\leq cost(\gamma)$.

The key point is that $\overline{\gamma}$ can be easily realized as a mapping made by the following elements:
\begin{itemize}
\item $(v,\mathfrak{D})$, $\forall v_d \in \gamma_d^T$;
\item $(v,\mathfrak{G})$, $\forall v_g \in \gamma^T_g$;
\item $(v,v')$, $\forall v_{\varphi,\varphi'} \in \gamma_s^T$, where $v'$ is the edge in $T'$ (more precisely, in $\gamma_g^{T'} \circ \gamma_d^{T'} (T')$) associated to the shrinking $v_{\varphi,\varphi'}$;
\item $(\mathfrak{G},v)$, $\forall v^{-1}_g \in (\gamma_g^{T'})^{-1}$;
\item $(\mathfrak{D},v)$, $\forall v^{-1}_d \in (\gamma_d^{T'})^{-1}$.
\end{itemize}

 \hfill$ \blacksquare $

\bigskip\noindent
\underline{\textit{Proof of}  \Cref{prop:altro_mapping}.}

\smallskip\noindent

Condition (M2) coincide with condition (A2). Condition (M3) is clearly satisfied because of the antichain condition (A1).
Consider a vertex $v\in E_T$. 
The only case in which $v$ is not edited is when $v<x$ with $x\in v_T\cap \pi_T(M^*)$. However, in this case, $v$ does not appear in $\widetilde{T}$,
and thus (M1) is satisfied. Moreover, after the deletions, all degree $2$ vertices are ghosted, and (M4) follows.

\hfill$ \blacksquare $

%\paragraph{Proof of Lemma~\Cref{lemma:proj_conv}}
\bigskip\noindent
\underline{\textit{Proof of}  \Cref{teo:decomposition}.} 

\smallskip\noindent
Let $M\in M_2(T,T')$ such that $d_E(T,T')=cost(M)$.

We note that $parent >child$ induces a partial order relation  also on the pairs given by paired points in $M$: $(x,y)>(v,w)$ if $x>v$ and $y>w$. In fact, by property (M3),  $x>v$ if and only if $y>w$.
So we can select $(x_i,y_i)$, the maxima with respect to this partial order relation . 
%Such procedure is well defined because, again by property (M3), maxima in $T$ are paired with maxima in $T'$. 
Thus, we obtain $(x_0,y_0)$,...,$(x_n,y_n)$ which form an antichain (both in $V_T$ and $V_{T'}$).

Clearly $M^*=\{(x_0,y_0)$,...,$(x_n,y_n)\}\in \mathcal{C}^*(T,T')$. Now we build $\alpha(M^*)$ and compare the cost of its edits with the ones in $M$.
Let $\bar{x}=LCA(x_i,x_j)$. Since $\bar{x}>x_i,x_j$, it is not paired in $M$. Since $x_i$ and $x_j$ are paired, $\bar{x}$ cannot be ghosted, so it is deleted in $M$. Any point $x$ above $\bar{x}$ is deleted for the same reasons. So the edits above $\bar{x}$ are shared between $\alpha(M^*)$ and $M$.

In $\alpha(M^*)$ we ghost any point between $\bar{x}$ and $x_i$ (and the same for $x_j$) and this is not certain to happen in $M$ (some points could be deleted). Nevertheless, even in the worst case, these ghostings are guaranteed not to increase the distance.
For instance, suppose $x_i<x<\bar{x}$ is deleted in $M$ and ghosted by $\alpha(M^*)$, then: 
\[
d(x_i\odot x,y_i)\leq d(x_i\odot x,y_i \odot x)+d(y_i \odot x,y_i)= d(x_i,y_i)+d(x,0)
\]
by properties (P1)-(P4).
Since $\alpha(M^*) \in M_2(\widetilde{T},\widetilde{T}')$ by  \Cref{prop:altro_mapping}, we have:
\[
\sum_{(x,y)\in M^*}d_E(sub_T(x),sub_{T'}(y))+cost(\alpha(M^*)) \leq cost(M)
\]

\smallskip

Now we prove the other inequality.

\smallskip

Consider $M^*$ which realizes the minimum of the right side of Equation \eqref{eq:deco_eqn}, and $M_i$ which realizes $d_E(sub(x_i),sub(y_i))$ with $(x_i,y_i)\in M^*$.	
We build a mapping $M$ collecting edits in the following way: for every $x'\in E_T$ if $x'\in sub(x_i)$, we take the edit associated to it from $M_i$, otherwise we know that it is edited by  $\alpha(M^*)$, and we take it from there; the set of these assignments gives $M\in M_2(T,T')$ whose cost is exactly 
$\sum_{(x_i,y_i)\in M^*} cost(M_i)+cost(\alpha(M^*))$. 
This gives the second inequality.

%Thus:
%
%\[
%\sum_{(x,y)\in M^*}d_E(sub_T(x),sub_{T'}(y))+cost(\alpha(M^*)) \geq cost(M)
%\]
 \hfill$ \blacksquare $

\bigskip\noindent
\underline{\textit{Proof of}  \Cref{prop:constraints}.}

\smallskip\noindent
Having fixed a leaf $l$, the constraint in Equation \eqref{eq:constr_1} allows for at most one path $\zeta_{v}^{v_{i+1}}\subset \zeta_l$ to be kept after the edits induced by all the variables equal to $1$. 
Moreover if $(v'',i) \in \Phi(v)\cap \Phi(v')$, then $v=v''_i = v'$. Thus, variables are added at most one time in \Cref{eq:constr_1} and \Cref{eq:constr_2}.
Which means that for any $a\in V_{T_x}$, we are forcing that the vertex $a$ can be an internal vertex or lower extreme of at most one path $\zeta_v^{v_{i+1}}$ such that $\delta_{i,j}^{v,w}=1$. In other words if two \virgolette{kept} paths $\zeta_v^{v_{i+1}}$ and $\zeta_{v'}^{v'_{i'+1}}$ (i.e. with $\delta_{i,j}^{v,w}=\delta_{i',j}^{v',w}=1$) intersect each other, it means that they just share the upper extreme $v_{i+1}=v'_{i'+1}$.
These facts together imply that (if the constraints are satisfied) the edits induced on $T_x$
by $\delta^{v,w}_{i,j}=1$ and $\delta^{v',w'}_{i',j'}=1$ are always compatible: if $v''\in V_{sub_T(v_{i})}$ then it is not touched by (the edits induced by) $\delta^{v',w'}_{i',j'}=1$ (and the same exchanging the role of $v'$ and $v$), if $v''$ is equal or above $v_{i+1}$ and/or $v'_{i'+1}$, then it is deleted in any case. 
Lastly, by noticing that if $\delta^{v,w}_{i,j}\in\Phi(v')$ then $\delta^{v,w'}_{i,j'}\in\Phi(v')$ for all other possible $w'$ and $j'$, we see that every path $\zeta_v^{v_i}$ is paired with at most one path $\zeta_w^{w_j}$, and viceversa.

As a consequence, for any vertex $v'$ in any of the tree structures, at most one point on the path $\zeta_{v'}$ is paired in $M^*$, guaranteeing the antichain condition.
Moreover, any point of $T_x$ which is in $\pi_{T_x}(M^*)$ is assigned to one and only one point of $T'_y$ and viceversa. The edits induced by $\delta=1$ clearly satisfy properties (M2)-(M4). Passing to $\widetilde{T}_x$ and $\widetilde{T}_{y}$, also (M1) is satisfied.

\hfill$ \blacksquare $

%\MPnote{mettere questione curvatura di alexandrov un po' meglio!}

\bibliography{references}

\begin{thebibliography}{60}
\providecommand{\natexlab}[1]{#1}
\providecommand{\url}[1]{\texttt{#1}}
\expandafter\ifx\csname urlstyle\endcsname\relax
  \providecommand{\doi}[1]{doi: #1}\else
  \providecommand{\doi}{doi: \begingroup \urlstyle{rm}\Url}\fi

\bibitem[Aho et~al.(1972)Aho, Garey, and Ullman]{transitive_reduction}
A.~V. Aho, M.~R. Garey, and J.~D. Ullman.
\newblock The transitive reduction of a directed graph.
\newblock \emph{SIAM Journal on Computing}, 1\penalty0 (2):\penalty0 131--137,
  1972.
\newblock \doi{10.1137/0201008}.

\bibitem[Ambauen et~al.(2003)Ambauen, Fischer, and Bunke]{ambauen2003graph}
R~Ambauen, Stefan Fischer, and Horst Bunke.
\newblock Graph edit distance with node splitting and merging, and its
  application to diatom identification.
\newblock In \emph{International Workshop on Graph-Based Representations in
  Pattern Recognition}, pages 95--106. Springer, 2003.

\bibitem[Bauer et~al.(2014)Bauer, Ge, and Wang]{bauer2014measuring}
Ulrich Bauer, Xiaoyin Ge, and Yusu Wang.
\newblock Measuring distance between reeb graphs.
\newblock In \emph{Proceedings of the thirtieth annual symposium on
  Computational geometry}, pages 464--473, 2014.

\bibitem[Bauer et~al.(2020)Bauer, Landi, and M{\'e}moli]{bauer2020reeb}
Ulrich Bauer, Claudia Landi, and Facundo M{\'e}moli.
\newblock The reeb graph edit distance is universal.
\newblock \emph{Foundations of Computational Mathematics}, pages 1--24, 2020.

\bibitem[Biasotti et~al.(2008)Biasotti, Giorgi, Spagnuolo, and
  Falcidieno]{reeb_2}
Silvia Biasotti, Daniela Giorgi, Michela Spagnuolo, and Bianca Falcidieno.
\newblock Reeb graphs for shape analysis and applications.
\newblock \emph{Theoretical Computer Science}, 392\penalty0 (1-3):\penalty0
  5--22, 2008.

\bibitem[Bille(2005)]{survey_ted}
Philip Bille.
\newblock A survey on tree edit distance and related problems.
\newblock \emph{Theoretical Computer Science}, 337\penalty0 (1-3):\penalty0 217
  -- 239, 2005.

\bibitem[Billera et~al.(2001)Billera, Holmes, and
  Vogtmann]{billera2001geometry}
Louis~J Billera, Susan~P Holmes, and Karen Vogtmann.
\newblock Geometry of the space of phylogenetic trees.
\newblock \emph{Advances in Applied Mathematics}, 27\penalty0 (4):\penalty0
  733--767, 2001.

\bibitem[Bollen(2022)]{bollen2022stable}
Brian Bollen.
\newblock \emph{Stable, Discriminative Distances on Reeb Graphs and Merge
  Trees}.
\newblock PhD thesis, The University of Arizona, 2022.

\bibitem[Botnan and Lesnick(2022)]{botnan2022introduction}
Magnus~Bakke Botnan and Michael Lesnick.
\newblock An introduction to multiparameter persistence.
\newblock \emph{arXiv:2203.14289}, 2022.

\bibitem[Bubenik and Elchesen(2022)]{bubenik2022universality}
Peter Bubenik and Alex Elchesen.
\newblock Universality of persistence diagrams and the bottleneck and
  wasserstein distances.
\newblock \emph{Computational Geometry}, 105:\penalty0 101882, 2022.

\bibitem[Cardona et~al.(2021)Cardona, Curry, Lam, and
  Lesnick]{cardona2021universal}
Robert Cardona, Justin Curry, Tung Lam, and Michael Lesnick.
\newblock The universal $\ell^p$-metric on merge trees.
\newblock \emph{arXiv}, 2112.12165 [cs.CG], 2021.

\bibitem[Carlsson and De~Silva(2010)]{carlsson2010zigzag}
Gunnar Carlsson and Vin De~Silva.
\newblock Zigzag persistence.
\newblock \emph{Foundations of computational mathematics}, 10:\penalty0
  367--405, 2010.

\bibitem[Carlsson and M\'emoli(2013)]{carlsson2013classify}
Gunnar Carlsson and Facundo M\'emoli.
\newblock Classifying clustering schemes.
\newblock \emph{Foundations of Computational Mathematics}, 13, 2013.

\bibitem[Carlsson and M{\'e}moli(2013)]{carlsson2013classifying}
Gunnar Carlsson and Facundo M{\'e}moli.
\newblock Classifying clustering schemes.
\newblock \emph{Foundations of Computational Mathematics}, 13:\penalty0
  221--252, 2013.

\bibitem[Carriere and Oudot(2018)]{carriere2018structure}
Mathieu Carriere and Steve Oudot.
\newblock Structure and stability of the one-dimensional mapper.
\newblock \emph{Foundations of Computational Mathematics}, 18:\penalty0
  1333--1396, 2018.

\bibitem[Carriere et~al.(2018)Carriere, Michel, and
  Oudot]{carriere2018statistical}
Mathieu Carriere, Bertrand Michel, and Steve Oudot.
\newblock Statistical analysis and parameter selection for mapper.
\newblock \emph{The Journal of Machine Learning Research}, 19\penalty0
  (1):\penalty0 478--516, 2018.

\bibitem[Cavinato et~al.(2022)Cavinato, Pegoraro, Ragni, Sollini, Erba, and
  Ieva]{cavinato2022imaging}
Lara Cavinato, Matteo Pegoraro, Alessandra Ragni, Martina Sollini, Paola~Anna
  Erba, and Francesca Ieva.
\newblock Author correction: Imaging-based representation and stratification of
  intra-tumor heterogeneity via tree-edit distance.
\newblock \emph{Scientific Reports}, 12\penalty0 (1):\penalty0 22540, 2022.

\bibitem[Chazal et~al.(2016)Chazal, De~Silva, Glisse, and
  Oudot]{chazal2016structure}
Fr{\'e}d{\'e}ric Chazal, Vin De~Silva, Marc Glisse, and Steve Oudot.
\newblock \emph{The structure and stability of persistence modules}.
\newblock Springer, 2016.

\bibitem[Curry(2018)]{curry2018fiber}
Justin Curry.
\newblock The fiber of the persistence map for functions on the interval.
\newblock \emph{Journal of Applied and Computational Topology}, 2, 2018.

\bibitem[Curry et~al.(2021)Curry, DeSha, Garin, Hess, Kanari, and
  Mallery]{curry2021trees}
Justin Curry, Jordan DeSha, Ad{\'e}lie Garin, Kathryn Hess, Lida Kanari, and
  Brendan Mallery.
\newblock From trees to barcodes and back again ii: Combinatorial and
  probabilistic aspects of a topological inverse problem.
\newblock \emph{arXiv}, 2107.11212, 2021.

\bibitem[Curry et~al.(2022)Curry, Hang, Mio, Needham, and
  Okutan]{curry2021decorated}
Justin Curry, Haibin Hang, Washington Mio, Tom Needham, and Osman~Berat Okutan.
\newblock Decorated merge trees for persistent topology.
\newblock \emph{Journal of Applied and Computational Topology}, 2022.

\bibitem[Curry et~al.(2023{\natexlab{a}})Curry, Mio, Needham, Okutan, and
  Russold]{curry2023convergence}
Justin Curry, Washington Mio, Tom Needham, Osman~Berat Okutan, and Florian
  Russold.
\newblock Convergence of leray cosheaves for decorated mapper graphs.
\newblock \emph{arXiv:2303.00130}, 2023{\natexlab{a}}.

\bibitem[Curry et~al.(2023{\natexlab{b}})Curry, Mio, Needham, Okutan, and
  Russold]{curry2023stability}
Justin Curry, Washington Mio, Tom Needham, Osman~Berat Okutan, and Florian
  Russold.
\newblock Stability and approximations for decorated reeb spaces.
\newblock \emph{arXiv preprint arXiv:2312.01982}, 2023{\natexlab{b}}.

\bibitem[De~Silva et~al.(2016)De~Silva, Munch, and Patel]{de2016categorified}
Vin De~Silva, Elizabeth Munch, and Amit Patel.
\newblock Categorified reeb graphs.
\newblock \emph{Discrete \& Computational Geometry}, 55\penalty0 (4):\penalty0
  854--906, 2016.

\bibitem[Di~Fabio and Landi(2016)]{di2016edit}
Barbara Di~Fabio and Claudia Landi.
\newblock The edit distance for reeb graphs of surfaces.
\newblock \emph{Discrete \& Computational Geometry}, 55\penalty0 (2):\penalty0
  423--461, 2016.

\bibitem[Divol and Lacombe(2021)]{divol2021understanding}
Vincent Divol and Th{\'e}o Lacombe.
\newblock Understanding the topology and the geometry of the space of
  persistence diagrams via optimal partial transport.
\newblock \emph{Journal of Applied and Computational Topology}, 5:\penalty0
  1--53, 2021.

\bibitem[Edelsbrunner and Harer(2008)]{PH_survey}
Herbert Edelsbrunner and John Harer.
\newblock Persistent homology---a survey.
\newblock In \emph{Surveys on discrete and computational geometry}, volume 453
  of \emph{Contemporary Mathematics}, pages 257--282. American Mathematical
  Society, Providence, RI, 2008.

\bibitem[Edelsbrunner and Harer(2010)]{edelsbrunner2022computational}
Herbert Edelsbrunner and John~L Harer.
\newblock \emph{Computational topology: an introduction}.
\newblock American Mathematical Society, 2010.

\bibitem[Felsenstein and Felenstein(2004)]{phylo}
Joseph Felsenstein and Joseph Felenstein.
\newblock \emph{Inferring phylogenies}, volume~2.
\newblock Sinauer associates Sunderland, MA, 2004.

\bibitem[Gao et~al.(2010)Gao, Xiao, Tao, and Li]{gao2010survey}
Xinbo Gao, Bing Xiao, Dacheng Tao, and Xuelong Li.
\newblock A survey of graph edit distance.
\newblock \emph{Pattern Analysis and applications}, 13:\penalty0 113--129,
  2010.

\bibitem[Garba et~al.(2021)Garba, Nye, Lueg, and
  Huckemann]{garba2021information}
Maryam~Kashia Garba, Tom~MW Nye, Jonas Lueg, and Stephan~F Huckemann.
\newblock Information geometry for phylogenetic trees.
\newblock \emph{Journal of Mathematical Biology}, 82\penalty0 (3):\penalty0
  1--39, 2021.

\bibitem[Gasparovic et~al.(2019)Gasparovic, Munch, Oudot, Turner, Wang, and
  Wang]{merge_intrins}
Ellen Gasparovic, E.~Munch, S.~Oudot, Katharine Turner, B.~Wang, and Yusu Wang.
\newblock Intrinsic interleaving distance for merge trees.
\newblock \emph{ArXiv}, 1908.00063v1[cs.CG], 2019.

\bibitem[Hatcher(2000)]{hatcher}
Allen Hatcher.
\newblock \emph{Algebraic topology}.
\newblock Cambridge University Press, Cambrige, 2000.

\bibitem[Hein et~al.(1995)Hein, Jiang, Wang, and Zhang]{np_hard}
Jotun Hein, Tao Jiang, Lusheng Wang, and Kaizhong Zhang.
\newblock On the complexity of comparing evolutionary trees.
\newblock In \emph{Combinatorial Pattern Matching}, pages 177--190, Berlin,
  Heidelberg, 1995. Springer Berlin Heidelberg.

\bibitem[Hong et~al.(2017)Hong, Kobayashi, and Yamamoto]{TED}
Eunpyeong Hong, Yasuaki Kobayashi, and Akihiro Yamamoto.
\newblock Improved methods for computing distances between unordered trees
  using integer programming.
\newblock In \emph{International Conference on Combinatorial Optimization and
  Applications}, pages 45--60. Springer, 2017.

\bibitem[Lerouge et~al.(2016)Lerouge, Abu-Aisheh, Adam, and
  H{\'e}roux]{lerouge2016exact}
Julien Lerouge, Zeina Abu-Aisheh, S{\'e}bastien~Romain Adam, and Pierre
  H{\'e}roux.
\newblock Exact graph edit distance computation using a binary linear program.
\newblock In \emph{Structural, Syntactic, and Statistical Pattern Recognition:
  Joint IAPR International Workshop, S+ SSPR 2016, M{\'e}rida, Mexico, November
  29-December 2, 2016, Proceedings}, volume 10029, page 485. Springer, 2016.

\bibitem[Levenshtein et~al.(1966)]{levenshtein1966binary}
Vladimir~I Levenshtein et~al.
\newblock Binary codes capable of correcting deletions, insertions, and
  reversals.
\newblock In \emph{Soviet physics doklady}, volume~10, pages 707--710. Soviet
  Union, 1966.

\bibitem[Lueg et~al.(2024)Lueg, Garba, Nye, and Huckemann]{lueg2024foundations}
Jonas Lueg, Maryam~K Garba, Tom~MW Nye, and Stephan~F Huckemann.
\newblock Foundations of the wald space for phylogenetic trees.
\newblock \emph{Journal of the London Mathematical Society}, 109\penalty0
  (5):\penalty0 e12893, 2024.

\bibitem[M{\'e}moli et~al.(2022)M{\'e}moli, Wan, and
  Wang]{memoli2022persistent}
Facundo M{\'e}moli, Zhengchao Wan, and Yusu Wang.
\newblock Persistent laplacians: Properties, algorithms and implications.
\newblock \emph{SIAM Journal on Mathematics of Data Science}, 4\penalty0
  (2):\penalty0 858--884, 2022.

\bibitem[Morozov et~al.(2013)Morozov, Beketayev, and Weber]{merge_interl}
Dmitriy Morozov, Kenes Beketayev, and Gunther Weber.
\newblock Interleaving distance between merge trees.
\newblock \emph{Discrete \& Computational Geometry}, 49:\penalty0 22--45, 2013.

\bibitem[Murtagh and Contreras(2017)]{dendro_1}
Fionn Murtagh and Pedro Contreras.
\newblock Algorithms for hierarchical clustering: An overview, ii.
\newblock \emph{Wiley Interdisciplinary Reviews: Data Mining and Knowledge
  Discovery}, 7\penalty0 (6):\penalty0 e1219, 2017.

\bibitem[Patel(2018)]{patel2018generalized}
Amit Patel.
\newblock Generalized persistence diagrams.
\newblock \emph{Journal of Applied and Computational Topology}, 1, 2018.

\bibitem[Pegoraro(2024{\natexlab{a}})]{pegoraro2024finitely}
Matteo Pegoraro.
\newblock A finitely stable edit distance for merge trees.
\newblock \emph{arXiv preprint arXiv:2111.02738v5}, 2024{\natexlab{a}}.

\bibitem[Pegoraro(2024{\natexlab{b}})]{pegoraro2024finitelyfunc}
Matteo Pegoraro.
\newblock A finitely stable edit distance for functions defined on merge trees.
\newblock \emph{arXiv}, 2108.13108v7 [math.CO], 2024{\natexlab{b}}.

\bibitem[Pegoraro and Secchi(2024)]{pegoraro2024functional}
Matteo Pegoraro and Piercesare Secchi.
\newblock Functional data representation with merge trees.
\newblock \emph{arXiv}, 2108.13147v6 [stat.ME], 2024.

\bibitem[Pont et~al.(2022)Pont, Vidal, Delon, and Tierny]{merge_wass}
Mathieu Pont, Jules Vidal, Julie Delon, and Julien Tierny.
\newblock Wasserstein distances, geodesics and barycenters of merge trees.
\newblock \emph{IEEE Transactions on Visualization and Computer Graphics},
  28\penalty0 (1):\penalty0 291--301, 2022.

\bibitem[Serratosa(2021)]{serratosa2021redefining}
Francesc Serratosa.
\newblock Redefining the graph edit distance.
\newblock \emph{SN Computer Science}, 2\penalty0 (6):\penalty0 438, 2021.

\bibitem[Singh et~al.(2007)Singh, M{\'e}moli, Carlsson,
  et~al.]{singh2007topological}
Gurjeet Singh, Facundo M{\'e}moli, Gunnar~E Carlsson, et~al.
\newblock Topological methods for the analysis of high dimensional data sets
  and 3d object recognition.
\newblock \emph{PBG@ Eurographics}, 2:\penalty0 091--100, 2007.

\bibitem[Skraba and Turner(2020)]{skraba2020wasserstein}
Primoz Skraba and Katharine Turner.
\newblock Wasserstein stability for persistence diagrams.
\newblock \emph{arXiv preprint arXiv:2006.16824}, 2020.

\bibitem[{Sridharamurthy} et~al.(2020){Sridharamurthy}, {Masood},
  {Kamakshidasan}, and {Natarajan}]{merge_farlocca}
R.~{Sridharamurthy}, T.~B. {Masood}, A.~{Kamakshidasan}, and V.~{Natarajan}.
\newblock Edit distance between merge trees.
\newblock \emph{IEEE Transactions on Visualization and Computer Graphics},
  26\penalty0 (3):\penalty0 1518--1531, 2020.

\bibitem[Sridharamurthy and Natarajan(2021)]{sridharamurthy2021comparative}
Raghavendra Sridharamurthy and Vijay Natarajan.
\newblock Comparative analysis of merge trees using local tree edit distance.
\newblock \emph{IEEE Transactions on Visualization and Computer Graphics},
  2021.

\bibitem[Stefanou(2020)]{stefanou2020tree}
Anastasios Stefanou.
\newblock Tree decomposition of reeb graphs, parametrized complexity, and
  applications to phylogenetics.
\newblock \emph{Journal of Applied and Computational Topology}, 4\penalty0
  (2):\penalty0 281--308, 2020.

\bibitem[Tai(1979)]{Tai}
Kuo-Chung Tai.
\newblock The tree-to-tree correction problem.
\newblock \emph{Journal of the ACM}, 26:\penalty0 422–433, 1979.

\bibitem[Touli(2020)]{merge_frechet}
Elena~Farahbakhsh Touli.
\newblock Frechet-like distances between two merge trees.
\newblock \emph{ArXiv}, 2004.10747v1[cs.CC], 2020.

\bibitem[Ushizima et~al.(2012)Ushizima, Morozov, Weber, Bianchi, Sethian, and
  Bethel]{ushizima2012augmented}
Daniela Ushizima, Dmitriy Morozov, Gunther~H Weber, Andrea~GC Bianchi, James~A
  Sethian, and E~Wes Bethel.
\newblock Augmented topological descriptors of pore networks for material
  science.
\newblock \emph{IEEE transactions on visualization and computer graphics},
  18\penalty0 (12):\penalty0 2041--2050, 2012.

\bibitem[Wetzels et~al.(2022)Wetzels, Leitte, and Garth]{merge_farlocca_2}
Florian Wetzels, Heike Leitte, and Christoph Garth.
\newblock Branch decomposition-independent edit distances for merge trees.
\newblock \emph{Comput. Graph. Forum}, 41\penalty0 (3):\penalty0 367--378,
  2022.

\bibitem[Xu and Tian(2015)]{dendro_2}
Dongkuan Xu and Yingjie Tian.
\newblock A comprehensive survey of clustering algorithms.
\newblock \emph{Annals of Data Science}, 2, 2015.

\bibitem[Yan et~al.(2019)Yan, Wang, Munch, Gasparovic, and
  Wang]{yan2019structural}
Lin Yan, Yusu Wang, Elizabeth Munch, Ellen Gasparovic, and Bei Wang.
\newblock A structural average of labeled merge trees for uncertainty
  visualization.
\newblock \emph{IEEE transactions on visualization and computer graphics},
  26\penalty0 (1):\penalty0 832--842, 2019.

\bibitem[Zeng et~al.(2009)Zeng, Tung, Wang, Feng, and Zhou]{zeng2009comparing}
Zhiping Zeng, Anthony~KH Tung, Jianyong Wang, Jianhua Feng, and Lizhu Zhou.
\newblock Comparing stars: On approximating graph edit distance.
\newblock \emph{Proceedings of the VLDB Endowment}, 2\penalty0 (1):\penalty0
  25--36, 2009.

\bibitem[Zhang(1996)]{zhang1996constrained}
Kaizhong Zhang.
\newblock A constrained edit distance between unordered labeled trees.
\newblock \emph{Algorithmica}, 15\penalty0 (3):\penalty0 205--222, 1996.

\end{thebibliography}

\end{document}